\documentclass[11pt]{article}
\usepackage[hidelinks]{hyperref}

\usepackage[english,activeacute]{babel}
\usepackage[utf8]{inputenc}
\usepackage[Lenny]{fncychap}
\usepackage{float}
\usepackage{bbm}

\usepackage{amssymb,amsmath,amsthm} 
\usepackage{graphics,graphicx} 
\usepackage{subfig}
\DeclareGraphicsExtensions{.pdf,.png,.jpg,.eps}
\usepackage{lscape}
\usepackage{amsfonts,mathrsfs}
\usepackage{multicol}
\usepackage{appendix}
\usepackage[usenames,dvipsnames]{color}
\definecolor{armygreen}{rgb}{0.29, 0.33, 0.13}
\usepackage{makeidx} 
\usepackage{appendix} 

\newcommand{\RR}{\mathbb{R}}

\newcommand{\Lc}{\mathcal{L}}
\newcommand{\B}{\mathcal{B}}

\newtheorem{Obs}{Remark}[section]

\numberwithin{figure}{section}
\numberwithin{table}{section}
\numberwithin{equation}{section}

\numberwithin{equation}{section}

\makeindex 
\begin{document}
\title{RBFs methods for null control problems of the Stokes system with Dirichlet and Navier--slip boundary conditions.}
\author{
Pedro Gonz\'alez--Casanova\footnotemark[1],\, Louis Breton\footnotemark[2],\,\,\,and\,Cristhian Montoya\footnotemark[3]
} 

\footnotetext[1]{Instituto  de Matem\'aticas, Universidad Nacional Aut\'onoma de M\'exico, M\'exico.  {\tt casanova@matem.unam.mx}}

\footnotetext[2]{Facultad de Ciencias,Universidad Nacional Aut\'onoma de M\'exico,  M\'exico.  {\tt louis.breton@ciencias.unam.mx}}

\footnotetext[3]{Departamento de Matem\'atica, Universidad T\'ecnica Federico Santa Mar\'ia, Casilla 110-V, Valpara\'iso, Chile.  {\tt cristhian.montoya@usm.cl}}

\date{}
\maketitle
\abstract{The purpose of this article is to introduce radial basis function, (RBFs), methods for solving null control problems for the Stokes system with few internal scalar 
	controls and Dirichlet or Navier--slip boundary conditions. To the best of our knowledge, it has not been reported in the literature any numerical approximation through 
	RBFs to solve the direct Stokes problem with Navier--slip boundary conditions. In this paper we fill this gap to show its application for solving the null control problem
	for the Stokes system. 
	To achieve this goal, we introduce two radial basis function solvers, one global and the other local, to discretized the primal and adjoint systems related to the control
	 problem. Both techniques are based on divergence-free global RBFs. Stability analysis for these methods is performed in terms of the spectral radius of the corresponding
	Gram matrices. By using a conjugate gradient algorithm, adapted to the radial basis function setting, the control problem is solved. Several test problems in two dimensions
	are numerically solved by these RBFs methods to test their feasibility. The solutions to these problems are also obtained by finite elements techniques, (FE), 
	to compare their relative performance.
	}\\

\section{\normalsize{Introduction}} 
	The main goal of this article is to introduce radial basis function, (RBFs), methods to solve null control problems for the Stokes system with few internal scalar controls 
	and Dirichlet or Navier--slip boundary conditions. In other to solve this problem, we first present fast and efficient radial basis functions (RBFs) algorithms for the Stokes
	 equations with Dirichlet or Navier--slip boundary conditions, and to solve from a numerical point of view the null control problem associated to the Stokes system with 
	 such boundary conditions. The control function is given as an external source acting over a subdomain and having few scalar controls, for example, 
	the control function can be $\textbf{v}=(v_1,0)$, where $v_1$ is a scalar field supported in a small subdomain.  First, we mention some recent discretization 
	schemes linked to the direct Stokes problem with Dirichlet conditions. The notion of the null control problem for the Stokes system will be mentioned later on.
	
	Stationary and evolutionary Stokes equations have been recently solved using analytically divergence-free and curl-free approximation spaces, 	
	\cite{wendland2009divergence} and \cite{fuselier2016high}. These schemes involves matrix-valued, positive definite kernels, as they have been introduced in
	 \cite{amodei1991vector} and \cite{narcowich1994generalized}. 
	
	In \cite{wendland2009divergence}, using the well known fact that the vector field $\mathbf{u}$ has a vector potential $\mathbf{w}$, i.e., that 
	$\mathbf{u} = \nabla\times \mathbf{w}$, the author define the potential as the curl of a scalar compact support radial basis function 
	$\psi$, i.e., $\mathbf{w} = \Delta\times \psi$. The field $\mathbf{u}$ is then approximated by a linear combination of the divergency-free matrix-valued 
	kernel $\mathbf{\Phi}_{Div} = \nabla\times \Delta\times \psi =\{ - \Delta I + \nabla \nabla^T \} \psi(\mathbf{x})$, where $\Delta$ is the Laplacian and $I$  the identity matrix. 
	A basic element in this approach, is that the author utilizes a vector composed by the velocity and the pressure, i.e., the vector $(\mathbf{u}, p)\in \RR^{d+1}$. Thus, 
	by forming an Ansatz build from a combination of a divergency free RBF kernel $\mathbf{\Phi}_{Div}$ and a scalar RBF kernel $\psi$, he is able to discretize the 
	stationary Stokes problem such that incompressibility is analytically satisfied.  Boundary conditions are incorporated to the discrete scheme, which is proved to be 
	invertible. {\it It is important to note that both the velocity and the pressure are coupled in the algebraic system of this approach}.

Within the same spirit, in \cite{fuselier2016high}, they built both a divergency free $\mathbf{\Phi}_{Div}$ and a rotational free $\mathbf{\Psi}_{Rot} = -\nabla \nabla^T\psi(\mathbf{x})$ matrix-valued kernels. Two major differences with respect to the first work, is that in \cite{fuselier2016high}, the radial basis function $\psi$ is now global and that they solve the evolutionary instead the stationary Stokes problem. In fact, based in the Leray projection, which is based in the Hodge decomposition theorem,  they first decouple the Stokes equation into a forced diffusion-type equation for $\mathbf{u}$ and then using the ortogonal complement of the Leray projector they obtain an auxiliary equation for the gradient of the pressure in terms of the velocity $\mathbf{u}$. These two equations are then discretized, for the space variable, by collocation using the RBFs kernel $\mathbf{\Phi}_{Div}$ for the velocity and $\psi$ for the pressure. Using the method of lines with a proper time integrator they solve the evolutionary Stokes problem.
	
We stress that in \cite{fuselier2016high}, unlike the current standard Leray-type projection methods
\cite{brown2001accurate}, the method proposed by the authors allows to build both tangential and normal boundary conditions into the discrete approximation to the Leray projection. This is in the spirit of RBFs global Hermite or symmetric methods for the solution of scalar PDEs, \cite{fasshauer2007meshfree}.  However, we emphasize that this scheme imposes non-slip boundary conditions. This is a consequence of the fact that the Leray projector applied to the stationary Stokes problem, eliminates the pressure from corresponding velocity equation,  which in turn means that we cannot impose slip boundary conditions using the decouple procedure as it is currently proposed by Fuselier, \cite{fuselier2016high}.
	
A first limitation of the method given in \cite{fuselier2016high} is that it incurs two dimension-dependent costs: an initialization cost of $O(N^3)$, and a cost of $O(N^2)$ per time-step, where $N$ is the total number of nodes. This is a consequence of the use of global RBFs to approximate the Leray projector. A second limitation, which is well known for RBF global methods, is that the condition number of the Gram matrix increases very fast as the number of nodes grows, thus limiting the size of the problems that can be solved.
\vskip 0.4cm
\noindent	
We note that the method introduced by  Kein and Wendland, \cite{Keim2016AHA}, is based upon the Leray projection, which converts the Stokes equation into a vector-valued heat equation. 
As mention above, our approach, heavily relies on Wendland method that incorporates the pressure into the solution.  However, a difference is that we use global instead compact support RBFs to build the incompressible vector. Also in our approach we do not use the Leray projection to transform the system into a vector-valued heat equation. Finally we incorporate Navier--slip boundary conditions into the direct solver of the stokes problem.
\vskip 0.4cm
\noindent
The approach to solve the direct Stokes problem that we introduce in this paper has two proposes. First we use Wendland formulation \cite{wendland2009divergence}, for global RBFs, to solve the stationary and evolutionary Stokes problems for the velocity and pressure simultaneously. Thus, we can perform the discretization for Navier--slip boundary conditions, and second we aim to suppress the former limitations of Fuselier algorithm. Specifically, we use a  RBFs  Local Hermite Method  instead the global Hermite or symmetric collocation technique.  By using extended precession it possible to solve problems having a large number of data. This is possible, whenever the value of the fill distance and shape parameter have values such that the condition number of the local matrices, in extended precision, are numerically well posed. Also the numerical complexity of the algorithm, (CPU time), is reduced with respect to global methods. It is worth noting that, although here we aim to investigate the extended precision approach, several alternatives to the ill-conditioned problem of RBFs collocation methods have been recently formulated, see for example \cite{Lehto2017, Fornberg2013, Fornberg2015} and references therein for a comprehensive review on this subject.  We also recall that in a different context, namely for interpolation of vector fields, both global techniques and LHI methods have been investigated by one of the authors, (see \cite{Cervantes2013},  \cite{Cervantes2018}).
\vskip 0.4cm
\noindent
Once the direct methods are formulated, our approximation strategy to solve the null control problem for the Stokes system is based in applying the conjugate gradient method (CGM) for PDE's, supported on the classical formulation for control problems under the form of appropriate convex quadratic optimization problems \cite{glowinski1994exact, glowinski2008exact}. In this work, numerical experiments are carried out in two dimensions using  distributed controls with either all scalar components or few scalar controls and taking into account the two types of boundary conditions, Dirichlet and Navier--slip. Moreover, we use finite elements to obtain numerical solutions of these problems and compare its performance against the RBF results.
\vskip 0.4cm
\noindent
Following in the context of the null control problem, the numerical setting for the RBF algorithm has thus two parts:
	\begin{itemize}
    	\item A fast and efficient RBFs solver for the evolutionary Stokes problems.
		\item A splitting iterative algorithm, namely the conjugate method (CGM), to solve the coupled system of equations (optimal system).
	\end{itemize}
	\vskip 0.4cm
It is worth mentioning that only a few articles have appeared in radial basis function methods in control theory. Up to our knowledge, only two works on optimal control by 
RBF techniques has been reported in this field, \cite{Pearson2013} and \cite{Casanova-Zavaleta}.
\vskip 0.4cm
	
	The paper is organized as follows. In Section \ref{section_problem_formulation} we give the continuous description of the null control problem. Section 
	\ref{section_LHI_scalar} we recall some algorithmic elements of the LHI method. The divergence free global method and stability results are developed in Section
	\ref{section.DFGM}. In Section \ref{section_LHI_Stokes}, we introduce a RBF--LHI approach for solving the Stokes 
	system. Finally,  in Section \ref{section_control} we conclude the paper by applying these RBFs methods to the null controllability 
	problem with few scalar controls and either Dirichlet or Navier--slip conditions.   Conclusions and final remarks are also included.
	
\vskip 0.4 cm

\section{\normalsize{Problem formulation}}\label{section_problem_formulation}

In this section we introduce the notation and the continuos setting of the Stokes control problem that will be numerically solved in this article.
\vskip 0.4 cm
Let us first introduce some notation. Let $\Omega$ be a nonempty bounded connected open subset of $\mathbb{R}^d$ ($d=2$ or $d=3$) of class $C^{\infty}$. Let $T>0$ and let $\omega\subset\Omega$ be a (small) nonempty open subset which is the control domain.
Let further $Q:=\Omega\times(0,T)$,\, $\Sigma:=\partial\Omega\times(0,T)$ and by $\boldsymbol{\nu}(x)$ the outward unit normal vector to $\Omega$ at the point $x\in\partial\Omega$. 
	
Moreover, let
   \begin{equation*}
	H:=\{ \mathbf{u}\in L^{2}(\Omega)^{d}: \nabla\cdot \mathbf{u}=0 \,\,\text{in}\,\, \Omega,\,\,\, 
	\mathbf{u}\cdot \boldsymbol{\nu}=0 \,\,
	\text{on}\,\,\partial\Omega \}
    \end{equation*}
and  
    \begin{equation*}
      V:=\{ \mathbf{u}\in H_{0}^{1}(\Omega)^{d}: \nabla\cdot \mathbf{u}=0 \,\,\text{in}\,\, \Omega\}. 
\end{equation*}

\vskip 0.4 cm
The continuous null control problem for the Stokes system with either Dirichlet  or  Navier--slip boundary conditions, that we are interested in, is defined as follows:
\vskip 0.4 cm
Given an initial data $y_0$, we are looking for a control function $\mathbf{v}=\mathbf{v}(x,t)$ acting in $\omega\times(0,T)$ with $supp\,\mathbf{v}\subset\omega\times(0,T)$  such that the solution of the problem 
\begin{equation}\label{intro.Stokes.system}
    	\left\{
    	\begin{array}{lll}
        	\mathbf{y}_{t}-\mu\nabla\cdot (D\mathbf{y})+\nabla p=\mathbf{v}1_{\omega} &\text{ in }& Q,\\
        	\nabla \cdot \mathbf{y}=0&  \text{ in }& Q, \\
        	\mbox{+BC}&\text{ on }&\Sigma,\\
        	\mathbf{y}(\cdot,0)=\mathbf{y}_0(\cdot) & \text{ in }&\Omega,
   	 	\end{array}
    	\right.
\end{equation} 
satisfies 
\begin{equation}\label{intro.finalstate.zero}
		\mathbf{y}(\cdot, T)=0\quad\mbox{in}\quad\Omega.	
\end{equation}
In \eqref{intro.Stokes.system}, $D$ denotes the symmetrized gradient of $\mathbf{y}$, $\mu>0$ is the viscosity coefficient and $p$ is the pressure . 
\vskip 0.4 cm
We focus our work in two types of boundary conditions (BC) on $\Sigma$, namely:
	\begin{equation}\label{intro.BC.Stokes}
		\underbrace{\mathbf{y}=0}_{(a)\,\, \mbox{Dirichlet}}\quad\mbox{or}\quad
		\underbrace{\mathbf{y}\cdot \boldsymbol{\nu}=0,\,\,\, (\sigma(\mathbf{y},p)\cdot \boldsymbol{\nu})_{tg}=0,}_{(b)
		\,\,\mbox{Navier--slip}}	
	\end{equation}
where $\sigma(\mathbf{y},p):=-pId+2\mu D\mathbf{y}$ is the stress tensor and $tg$ stands for the tangential component of the corresponding vector field, i.e., $$\mathbf{y}_{tg}=\mathbf{y}-(\mathbf{y}\cdot \boldsymbol{\nu})\boldsymbol{\nu}.$$
	
	From a Physical point of view, Navier--slip boundary condition arises from the interaction between wall and fluid and when the temperature is high. These behavior involves a movement on the boundary (slip), loosing energy, which do not penetrate the boundary (impermeable boundary), among other factors. The use of these slip conditions allow to describe phenomena observed in nature and remove un--physical singularities, see for instance \cite{cebeci2012analysis}, \cite{he2009numerical} and references therein for more details.
	Now, from a mathematical point of view, such boundary conditions say that the tangential component of the stress tensor is proportionality to the tangential component of the velocity \cite{navier1823memoire}: $$\mathbf{y}\cdot \boldsymbol{\nu}=0,\,\,\, (\sigma(\mathbf{y},p)\cdot \boldsymbol{\nu})_{tg}+k\mathbf{y}_{tg}=0,$$ where $k$ is a function that measures the local viscous coupling fluid--solid. We highlight that this proportionality factor can depend on  the velocity as well as on the pressure, which complicate both the theoretical analysis and numerical solutions. In our numerical experiments we only deal the case \eqref{intro.BC.Stokes}--(b) and invite to interested reader to see \cite{cebeci2012analysis, navier1823memoire} for a complete discussion on this subject.
\vskip 0.4 cm

We now characterize the control problem in terms of the optimal or minimum value of a quadratic convex funcional in $(L^2(Q))^2$ in the sense of \cite{garciamontoyaosses17}. 

Namely, for $\mathbf{y}_0\in H$, we aim to obtain the control $\mathbf{v}$  with one vanishing component ($j$th component, $j\in \{1,2\}$) such that it minimizes the functional $J$ defined by
\begin{equation}\label{functional.Stokes1}
J(\mathbf{v}):=\displaystyle\frac{1}{2}\iint\limits_{\omega\times(0,T)}|\mathbf{v}|^{2}\,dx\,dt
+\displaystyle\frac{1}{2c_1}\|\mathbf{y}(.,T)\|^{2}_{L^2(\Omega)}\,dx
+\displaystyle\frac{1}{2c_2}\iint\limits_{\omega\times(0,T)}|v_j|^{2}\,dx\,dt,
\end{equation}
	where $\mathbf{y}$ is solution of the Stokes system \eqref{intro.Stokes.system},  $c_1, c_2$ are arbitrary positive
	numbers associated respectively to the final condition $\mathbf{y}(\cdot, T)=0$ and control function $\mathbf{v}$.

	As it is well known, this problem can be restated in terms of the Lagrangian formulation. Let  
	${\cal L}:[H_{0}^{1}(\Omega)^2]^2\times [L^{2}_0(\Omega)]^2\times L^2(Q)^2\rightarrow\mathbb{R}$
	defined by:
	\begin{equation}\label{lagragian.stokes1}
	\begin{array}{lll}
	{\cal L}(\mathbf{y},\mathbf{w},q,p,\mathbf{v}):=&\displaystyle\frac{1}
	{2}\iint\limits_{\omega\times(0,T)}|\mathbf{v}|^{2}\,dx\,dt
	+\displaystyle\frac{1}{2c_1}\|\mathbf{y}(.,T)\|^{2}_{L^2(\Omega)}\,dx
	+\displaystyle\frac{1}{2c_2}\iint\limits_{\omega\times(0,T)}|v_j|^{2}\,dx\,dt\\
 		&\quad +\displaystyle\iint\limits_{Q}\mathbf{y}_{t}\cdot \mathbf{w}-\mu\nabla\cdot (D\mathbf{y})\cdot\mathbf{w}
 		-div(\mathbf{w}) p-div(\mathbf{w})q\,dx\,dt\\
 		&\quad-\displaystyle\iint\limits_{\omega\times(0,T)}\mathbf{v}\cdot \mathbf{w}\,dx\,dt+
 		\displaystyle\int\limits_{\Omega}
 		(\mathbf{y}(\cdot,0)-\mathbf{y}_{0})\cdot
 		\mathbf{w}(\cdot,0)\,dx,
	\end{array}
	\end{equation}
	where $\mathbf{w}$ is the Lagrange multiplier. The optimal value of the control $\mathbf{v}$ can be obtained in 
	the following way.  It is easy to verify that the Fr\'echet derivate of $J$ with respect to $\mathbf{v}$ is:
	\begin{equation}\label{Stokes1.optimality.condition}
		\frac{\partial J}{\partial \mathbf{v}}(\mathbf{v})=v_i-w_i\,\,\mbox{if}\,\, i\neq j\quad \mbox{and}\quad
		\frac{\partial J}{\partial \mathbf{v}}(\mathbf{v})=\displaystyle\frac{1}{2c_2}v_j-w_j,\quad\mbox{in}\,\,
	\omega\times(0,T), 
	\end{equation}
	where $\mathbf{w}\in V$ is the solution of the adjoint system of \eqref{intro.Stokes.system}:
	\begin{equation}\label{Stokes.adjointsystem1}
    	\left\{
    	\begin{array}{lll}
        	-\mathbf{w}_{t}-\mu\nabla\cdot D\mathbf{w}+\nabla q=0 &\text{ in }& Q,\\
        	\nabla \cdot \mathbf{w}=0&  \text{ in }& Q, \\
        	\mbox{+BC}&\text{ on }&\Sigma,\\
        	\mathbf{w}(\cdot,T)=-\displaystyle\frac{1}{c_1}\mathbf{y}(\cdot,T) & \text{ in }&\Omega.
   	 	\end{array}
    	\right.
	\end{equation}

	In \cite{garciamontoyaosses17} the authors proved that for every $c_1>0,\, c_2>0$, there 
	exists a unique minimal control $\mathbf{v}$ associated to \eqref{functional.Stokes1} which is characterized by 
	the optimality system given by \eqref{intro.Stokes.system}, \eqref{Stokes1.optimality.condition},  
	\eqref{Stokes.adjointsystem1} with Dirichlet boundary conditions \eqref{intro.BC.Stokes}--(a). 

	We stress that the Lagrangian formulation \eqref{lagragian.stokes1} is the continuous setting for the mixed FE 
	formulation that will be used in the subsection of numerical results.  	
 	\vskip 0.4 cm
	From an abstract point of view, the control theory for the Stokes system 
	with internal controls has been studied intensively for the 
    	mathematical community. The interested reader can be see 
    	\cite{fursikov1999exact,FCGIP04,corona,imanuvilov,guerrero2006local,carreno2013local,
    	coron2014local,guerrero2018local} and references therein for a complete description.
    	Meanwhile, from a computational point of view, we only know the recent paper by
	Fernandez--Cara et. al, \cite{fernandez2017numerical}, whose numerical experiments are developed 
	in two dimension for the heat, Stokes and Navier--Stokes with Dirichlet boundary conditions. 
	The implemented methodology in
	\cite{fernandez2017numerical} for the fluid equations is based in the so called Fursikov--Imanuvilov 
	formulation \cite{fursikov1996controllability} and Lagrangian approximation throughout mixed finite
	elements. Another approximation scheme to the null control problem is given in
	\cite{fernandez2015theoretical} for a turbulence model and also using Dirichlet boundary conditions. 
	
	\vskip 0.4cm
	We close this section by pointing out that, as far as we know, it does not exists a numerical approximation through 
	RBFs for the Stokes problem with Navier--slip boundary conditions. In the following sections we fill this gap and 
	show its application for solving the null control problem for the Stokes system with few scalar controls. 

\section{\normalsize{Numerical method: notation and preliminary remarks}}\label{section_LHI_scalar}
	In what follows, and for the sake of completeness, we first briefly recall the scalar LHI method,  introduced by Stevens, et al. \cite{stevens2011local}. (see \cite{Casanova-Zavaleta} for a similar notation). This scalar setting will be used later in this paper to formulate the generalized vectorial technique to solve the Stokes problem. 
	
	In the LHI scalar approach, we aim to obtain the RBF approximation of the analytic solution
	$u$ of a linear steady partial differential well posed problem
	\begin{equation}\label{equation:solctrsa}
    	\left\{
    	\begin{array}{lll}
        	\Lc u(x)=f(x) &\text{ in }& \Omega,\\
        	\B u(x)=g(x)&  \text{ on }& \partial\Omega,
        \end{array}
    	\right.
	\end{equation} 
	\begin{figure}[ht!]   
	\begin{center}
  		\includegraphics[scale=0.45]{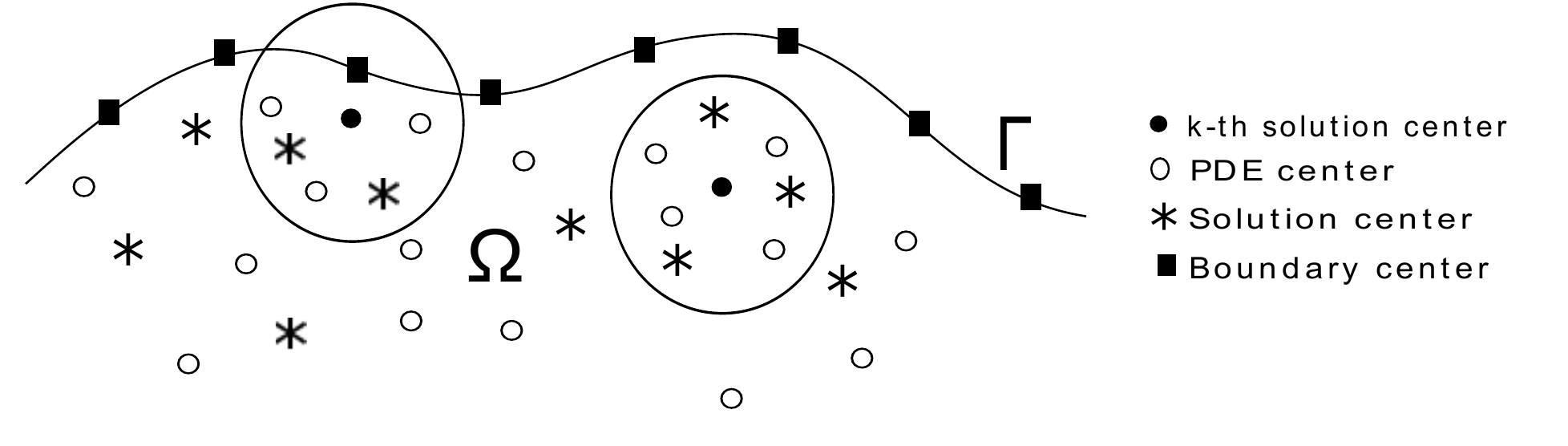}
  		\caption{ Centers and local subdomains for the LHI method}
			\label{fig:localsd} 
	\end{center}
	\end{figure}
	where $\Omega\subset \mathbb{R}^d$ represents the spatial domain, the right--hand sides 
	$f: \Omega\to \mathbb{R}$ and $g: \partial\Omega\to \mathbb{R}$ are given,   $\Lc $ and $\B$ are 
	linear partial differential operators in the domain $\Omega$ and on the contour $\partial\Omega$,  
	which are locally approximated in the following way.
	First, let $\Omega_{n} \subset {\overline{\Omega}}$ be a set of $n_t$ total number of scattered nodes. Consider now  
	the following subsets of $\Omega_{n}$. Let $\Omega_{c} \subset \Omega_{n}$ be a subset of $n_c$ nodes called 
	centers, see Figure \ref{fig:localsd}. 
	
	Let now $D^k$ be a disk, of variable radius, with center at the $k^{th}$ node of $\Omega_{c}$, 
	(recall that $\Omega_{c}$ is the set of centers of the disks), and consider the set of $n^{(k)}$ 
	fixed number of nodes $\Omega_n\cap D^k$,  see again  Figure \ref{fig:localsd}. To perform the local discretization 
	we introduce the following notation:
	\begin{equation*}
	\begin{array}{ll}
	& \Omega_{sc} =\{x_1^{(k)},\dots,x_{n_{sc}^k}  \} \subset \Omega_c\,\, \mbox{be a set of $n_{sc}^k$ nodes
	(called solution centers)}.\\	
	
	& \partial \Omega_{fc} =\{  x_{n_{sc}^k +1},\ldots,x_{n_{sc}^k +
	n_{fc}^k} \} \subset \partial \Omega\,\,\mbox{the boundary nodes}.\\
	
	&\Omega_{pdec} = \{x_{n_{sc}^k + n_{fc}^k+ 1},\ldots,x_{n_{sc}^k + n_{fc}^k + n_{pdec}^k} \} \subset \Omega 
	\,\, \mbox{ is a set of interior nodes related $\Lc$}.
	\end{array}
	\end{equation*}
	By simplicity, we shall denote the node $x_1^{(k)}$ as the center of the disk $D^k$ for each $k$. Then, 
	every disk $D^k$  
	define a local sub-system as follows:
	\begin{align}
  		\Lc u(x)&=f(x)  \hspace{1cm} x   \in \Omega_{pdec},  \label{equation:pdectrs} \\
  		\B u(x)&=g(x)  \hspace{1cm} x  \in \partial \Omega_{fc}, \label{equation:bdyctrs}\\
  		 u(x_i)&=h_i   \hspace{1.3cm} x_i \in  \Omega_{sc}   \label{equation:solctrs},
	\end{align}
	where $h_i$ are the unknown values.
	  	
	This procedure generates a set of local linear systems given by
	\begin{equation}\label{eq:sl}
  		A^{(k)}\beta^{(k)}= d^{(k)},
	\end{equation}
	\noindent
	which is obtained by substituting in (\ref{equation:pdectrs})--(\ref{equation:solctrs}) the 
	following radial ansatz for the local domains displayed (in circles) in Figure \ref{fig:localsd}
	\begin{equation}\label{lhia}
	\begin{array}{ll}
		\hat{u}^{(k)}(x)&=\sum\limits_{j=1}^{n_{sc}^k} \beta_j^{(k)} \phi_j( r ) 
		+\sum\limits_{j=n_{sc}^k+1}^{n_{sc}^k+n_{fc}^k} \beta^{(k)}_j\B^{\xi}\phi_j(r) 
		+\sum\limits_{j=n_{sc}^k+n_{fc}^k+1}^{n_{sc}^k+n_{fc}^k+n_{pdec}^k}\beta_j^{(k)} \mathcal{L}^{\xi
		}\phi_j(r)
		+ p_k^m,
	\end{array}	 
	\end{equation}
 	where $k$ is the local system index, $n_{sc}^k$ the number of
	solution centers in the local system, $n_{pdec}^k$ denotes the number of PDE centers in the local system and
	$n_{fc}^k$ is the number of boundary centers in the local system. Besides,  $\Lc^{\xi} \phi_j(r) := \Lc
	\phi( \| x-\xi \|)|_{\xi=\xi_j}$, $\B^{\xi} \phi_j(r) := \B \phi( \| x-\xi \|)|_{\xi=\xi_j} $, being 
	$\phi(r)$ the 
	inverse multi-quadric and $p_k^m$ a polynomial in $\mathbb{R}^d$ of degree $m$, which is an element of
	the null space of \eqref{equation:solctrsa}. The momentum condition is also required in this step, see 
	\cite{stevens2011local}. Thus, the local linear system (\ref{eq:sl}) can be expressed in vectorial notation 
	by defining $A^{(k)}$ (called Gram matrix) and the right hand vector $d^{(k)}$ as follows
	\begin{displaymath}
 	 A^{(k)} = \left[
    \begin{matrix}
      \Phi_{ij} &\B^\xi[\Phi_{ij} ] & \Lc^\xi[\Phi_{ij}]     & P_{ij}   \\
       \B^{x}[\Phi_{ij}] & \B^{x}\B^{\xi}[\Phi_{ij}] & \B^{x}\Lc^{\xi}[\Phi_{ij}]   & \B^{x}[P_{ij}] \\
      \Lc^{x}[\Phi_{ij}]  & \Lc^{x} \B^{\mathbf{\xi}}[\Phi_{ij}]  & \Lc^{x} \Lc^{\mathbf{\xi}}[\Phi_{ij}]
        &   \Lc^{x} [P_{ij}]\\
       P_{ji} &  \B^{\xi} [P_{ji}]  & \Lc^\xi  [P_{ji}]    & 0 \\ 
    \end{matrix}
    \right] \quad \textnormal{ and } \quad
    d^{(k)} = 
    \left[
      \begin{matrix}
        h_i\\
         g_i\\
        f_i \\
        0 \\
      \end{matrix}
      \right],
	\end{displaymath}
	which is well known to be invertible, see \cite{wendland2004scattered}. Thus, we have that 
	$\beta^{(k)}=(A^{(k)})^{-1}d^{(k)}$, and using \eqref{lhia},  $\hat{u}^{(k)}(x)$ can be rewritten by 
	
	\begin{equation}\label{eq:lapp}
    	\hat{u}^{(k)}(x)=H(x)\beta^{(k)}=H(x)( A^{(k)})^{-1}d^{(k)}=W^{(k)} d^{(k)},
	\end{equation}
	where
	\begin{displaymath}
  		H^{(k)}(x)=\left[
    	\phi(\|x-\xi\|) \quad \B^{\xi} \phi(x-\xi) \quad  \Lc^{\xi}\phi(\|x-\xi\|)  \quad p^m(x)    \right]
    	\quad \mbox{and} \quad W^{(k)}:=H(x)( A^{(k)})^{-1}.
	\end{displaymath}
	$W^{(k)}$ is known as the vector of weights, \cite{wendland2004scattered}. Now, if $\cal{J}$ is an arbitrary
	differential operator given, we can compute its value at $\hat{u}^{(k)} $ by
	
	\begin{displaymath}
	\label{eq:system4}
		{\cal{J}} \hat{u}^{(k)} (x)={\cal{J}} H(x)( A^{(k)})^{-1} d^{(k)}={\cal{J}} (W^{(k)})(x)d^{(k)}. 
	\end{displaymath}
	Let $u_c = \left [ u(x_1^{(k)})\right ]_{k=1}^{n_c^{(k)}}$ be the set of values of the 
	exact solution at the centers of each disk $D^k$, which are unknown values (they belong to 
	the vector $d^{(k)}$). In order to obtain these unknown values, we consider the following system of equations
	
	\begin{equation}
  	\label{eq:system}
		h \left (x_1^{(k)}\right ) = W^{(k)}_{\cal{J}}\left (x_1^{(k)}\right )  d^{(k)},\quad 
		k = 1,\ldots ,n_c,
	\end{equation}
	\noindent
	where $W^{(k)}_{\cal{J}} =  {\cal{J}}  (W^{(k)})$ and $(h, {\cal J})$ is defined by $(f, {\cal L})$.

	We shall denote by $S u_c = b$ the linear system  \eqref{eq:system}, whose variables are the
	values at the solution centers $\Omega_{c}$. Recall that each row of the matrix $S$ is composed by zero 
	elements except for the  weights, (centers), corresponding to each disk $D^k$, that is, each row of size 
	$n_c$, 
	has only $n_c^{(k)}$ elements different form zero. Moreover, since $n_c \gg  n_c^{(k)}$, the matrix $S$ 
	is sparse. 
	To built the matrix $S$ efficiently, i.e., to compute the weights, we can solve the following equations
	\begin{equation}\label{eq:calculo de pesos}
		A^{(k)}W^{(k)}_{\cal{J}}\left (x_1^{(k)}\right )  = {\cal{J}} H^{(k)}(x_1^{(k)}), \quad 
		k = 1,\ldots ,n_c.
	\end{equation}
	Since the the matrix $S$ is sparse standard solvers and preconditioners can be used. Besides, it is worth pointing out that by using the method of lines and a proper numerical time integrator, non stationary linear PDEs problems can be solve by the LHI method. To review exhaustively the LHI method, the interested reader can see 
\cite{stevens2011local, fasshauer2007meshfree} and references therein.

\section{\normalsize{Divergence free global method for evolutionary problems}}\label{section.DFGM}
	In this section we present two  methods for solving the time dependency for the unsteady Stokes system. 
	The spatial dependency follows divergence--free RBF approximation like in \cite{wendland2009divergence}.
	Let us define $L(\mathbf{y},p)=-\mu\Delta \mathbf{y}+\nabla p$ and consider the system
	\begin{equation}\label{Stokes.system}
		\left\{ \begin{array}{lll}
		\mathbf{y}_{t}+L(\mathbf{y},p)=\mathbf{f} & \text{ in } & Q,\\
		\nabla\cdot\mathbf{y}=0 & \text{ in } & Q,\\
		B\mathbf{y}=\mathbf{g} & \text{ on } & \Sigma,\\
		\mathbf{y}(\cdot,0)=\mathbf{y}_{0}(\cdot) & \text{ in } & Q,
	\end{array}\right.
	\end{equation} 	
	where $B$ represents only one boundary operador describing \eqref{intro.BC.Stokes}.
	
\subsection{\normalsize{Time backward scheme: global RBF collocation}}\label{subsection.operator.discre}

	The goal is to create a modified PDE--operator via a finite difference approximation for the time derivative, 
	where the left--hand side is unknown and right--hand side is known. The time scheme follows some ideas by 
	\cite{stevens2011local}. In order to  illustrate this method we use 
	backward finite differences (BFD) methods that are appropriate for the Stokes equations. 
	The formulation at any time step for the system \eqref{Stokes.system} is given as follows:
	\begin{equation}\label{bdf.Stokes.system}
	\left\{ \begin{array}{lll}
		\mathbf{y}^{n+s}+\Delta t\,\beta_{s}\,L(\mathbf{y}^{n+s},p^{n+s})
		=\Delta t\,\beta_{s}\,\mathbf{f}^{n+s}
		+\sum\limits_{k=0}^{s-1}\sigma_{k}\mathbf{y}^{n+k} & \text{ in } & Q,\\
		\nabla\cdot\mathbf{y}^{n+s}=0 & \text{ in } & Q,\\
		B\mathbf{y}^{n+s}=\mathbf{g}^{n+s} & \text{ on } & \Sigma,
	\end{array}\right.
	\end{equation} 
	\noindent
	where $\beta_s,\sigma_k$ are known parameters within BFD methods. Thus, at each step we solve the 
	following modified PDE 
	\begin{equation}\label{DFG.modifiedPDE}
	\left\{ \begin{array}{lll}
		\bar{L}(\mathbf{y}^{n+s},p^{n+s})=\boldsymbol{F}^{n+s} & \text{ in } & Q,\\
		\nabla\cdot\mathbf{y}^{n+s}=0 & \text{ in } & Q,\\
		B\mathbf{y}^{n+s}=\mathbf{g}^{n+s} & \text{ on } & \Sigma,
	\end{array}\right.
	\end{equation} 
	where the source term is computed using previous steps and $\bar{L},\,\boldsymbol{F}^{n+s}$ are defined by
	\begin{equation*}
		\bar{L}(\mathbf{y}^{n+s},p^{n+s}):=\mathbf{y}^{n+s}+\Delta t\,\beta_{s}\,L(\mathbf{y}^{n+s},p^{n+s}),\quad\quad
		\boldsymbol{F}^{n+s}:=\
		\Delta t\,\beta_{s}\,\mathbf{f}^{n+s}+\sum_{k=0}^{s-1}\sigma_{k}\mathbf{y}^{n+k}.
	\end{equation*} 
	
	\noindent The field $\mathbf{y}^{n+s}$ is then approximate by a linear combination of the divergence-free 
	matrix--valued kernel 
	$\mathbf{\Phi}_{Div}=\nabla\times\Delta\times\psi=\{-\Delta I+\nabla\nabla^{T}\}\psi$,
	where $\Delta$ is the Laplacian, $I$ is the identity matrix and $\psi:\mathbb{R}^d\to \mathbb{R}:$ is 
	a positive definitive function. Following \cite{wendland2009divergence}, the velocity pressure
	vector $(\mathbf{y}^{n+1},p^{n+1})$ can be approximate by a linear combination of:
	\[
	\mathbf{\Phi}=\left[{\begin{array}{cc}
	\mathbf{\Phi}_{Div} & 0\\
	0 & \phi
	\end{array}}\right]:\mathbb{R}^{d}\rightarrow\mathbb{R}^{(d+1)\times(d+1)},
	\]
	where $\phi:\mathbb{R}^d\to \mathbb{R}$ is a positive definitive function.
	
	On the other hand, using the generalized interpolation collocation method, see \cite{wendland2004scattered},  
	The RBF ansatz, using Navier--slip boundary conditions (see \eqref{intro.BC.Stokes}--(b)), is given by 
	\begin{equation}\label{Anzatdist.Stokes.system}
		(\hat{\mathbf y}^{n+s},\hat{p}^{n+s})(\boldsymbol{x})=\sum_{i=1}^{d}\sum_{j=1}^{n_{b}}\left(B_{ij}
		\right)\Phi(\boldsymbol{x-z})\alpha_{ij}^{B}
		+\sum_{i=1}^{d}\sum_{j=1}^{n_{in}}\left(\overline{L_{ij}}\right)\Phi(\boldsymbol{x-z})\alpha_{ij}^{L},
	\end{equation}
	where $n_b,\, n_{in}$ are the total numbers of boundary and interior nodes respectively, 
	$\alpha=(\alpha_{ij}^B,\alpha_{ij}^L)\in\mathbb{R}^{d(n_{b}+n_{in})}$ and 
	$B_{ij}$,$\overline{L_{ij}}$ are vector-valued functionals defined as follows:
	\begin{equation}\label{DFG.ansatz.slip.functionals1}
		\overline{L_{ij}}=\delta_{\boldsymbol{x_{j}^{in}}}e_{i}+\beta_{s}\Delta t(\delta_{\boldsymbol{x_{j}^{in}}}\circ-\mu\Delta 
		e_{i}+\delta_{\boldsymbol{x_{j}^{in}}}\circ\partial_{i}e_{d+1})
	\end{equation}
	and $B_{ij}$ is
	\begin{equation}\label{DFG.ansatz.slip.functionals2}
		\delta_{\boldsymbol{x_{j}^{b}}} \sum_{k=1}^d e_k\nu_k \mbox{ for } i=1 \quad\mbox{and}\quad 
		\delta_{\boldsymbol{x_{j}^{b}}}\circ \mu \sigma(\bar{\boldsymbol{e}}_d,e_{d+1})\boldsymbol{\nu}.\tau_{i-1}, \mbox{ for } i=2...d 
	\end{equation}
	where $\boldsymbol{\nu}$ the outward unit normal vector to $\Omega$, $\bar{\boldsymbol{e}}_d = (e_1,...,e_d)$ and ($\tau_1,...,\tau_{d-1}$) is orthonormal basis 
	of the tangent space.
	  
	\noindent Now, in order to simplify the notation, we shall proceed our analysis in 2D by defining the vector value
	function $(\phi^{B_{1}},\phi^{B_{2}}, \phi^{\overline{L_{1}}},\phi^{\overline{L_{2}}})$
	such that the ansatz \eqref{Anzatdist.Stokes.system} can be rewritten in the form
	\begin{eqnarray}\label{eq:compact_anzat}
		(\hat{\mathbf y}^{n+s},\hat{p}^{n+s})(\boldsymbol{x})=&
		 &\hspace{-0,5cm} \sum_{j=1}^{n_{b}}\phi^{B_{1}}(\boldsymbol{x}-\boldsymbol{x_{j}^{b}})\alpha_{1j}^{B}
	 	+\sum_{j=1}^{n_{b}}\phi^{B_{2}}(\boldsymbol{x}-\boldsymbol{x_{j}^{b}})\alpha_{2j}^{B}\\
 		&&\quad + \sum_{j=1}^{n_{in}}\phi^{\overline{L_{1}}}(\boldsymbol{x}-\boldsymbol{x_{j}^{in}})\alpha_{1j}^{L} 
 		+\sum_{j=1}^{n_{in}}\phi^{\overline{L_{2}}}(\boldsymbol{x}-\boldsymbol{x_{j}^{in}})\alpha_{2j}^{L}. \nonumber
	\end{eqnarray}

	\noindent Once that the ansatz has been replaced in \eqref{DFG.modifiedPDE}, we obtain the discrete system
	\begin{align}\label{eq:global scheme}
		H_{+}\bar{\alpha{}}^{n+s}=\left(\begin{array}{c}
		\boldsymbol{F}^{n+s}\\
		\mathbf{g}^{n+s}
	\end{array}\right),
	\end{align}
	where the collocation matrix $H_{+}$ is given by:
	\[
	H_{+}=\left(\begin{array}{ccccc}
		\overline{L_{1}}\phi^{B_{1}} & \overline{L_{1}}\phi^{B_{2}} & 
		\overline{L_{1}}\phi^{\overline{L_{1}}} & \overline{L_{1}}\phi^{\overline{L_{2}}}\\
		\overline{L_{2}}\phi^{B_{1}} & \overline{L_{2}}\phi^{B_{2}} & 
		\overline{L_{2}}\phi^{\overline{L_{1}}} & \overline{L_{2}}\phi^{\overline{L_{2}}}\\
		B_{1}\phi^{B_{1}} & B_{1}\phi^{B_{2}} & 
		B_{1}\phi^{\overline{L_{1}}} & B_{1}\phi^{\overline{L_{2}}}\\
		
		B_{2}\phi^{B_{1}} & B_{2}\phi^{B_{2}} & 
		B_{2}\phi^{\overline{L_{1}}} & B_{2}\phi^{\overline{L_{2}}}
		
		\end{array}\right).
	\]
	\begin{Obs}\label{remark.timescheme.glogalscheme}
		In the case of Dirichlet boundary conditions, the vector--value functional $B_{ij}$ defined in 
		\eqref{DFG.ansatz.slip.functionals2} is replaced by $B_{ij}=\delta_{\boldsymbol{x_{j}}} e_i$, and 
		the collocation matrix $H_+$ is similar to the previous one.  
	\end{Obs}
	

	\subsection{\normalsize{Stability analysis}}\label{subsection.stability.global1}
	 In this subsection, we proceed to present the stability analysis associated to the previous scheme,  by using a  
	 matrix method similar to the procedure developed in \cite{chinchapatnam2006unsymmetric}. First, using equation (\ref{eq:compact_anzat}) we 
	 define the interpolation matrix $A\in\mathbb{R}^{(n_{in}+n_{b})\times(n_{in}+n_{b})}$ such that
	 
	\begin{equation}\label{eq:interporlation matrix}
		A\bar{\alpha}^{n+s}=\left(\begin{array}{c}
		M_{\phi_{in}}\\
		M_{\phi_{b}}
	\end{array}\right)\bar{\alpha}^{n+s}
	=\left(y^{n+s}(x_{1}^{in}),\dots,y^{n+s}(x_{n_{in}}^{in}),y^{n+s}(x_{1}^{b}),\dots ,y^{n+s}(x_{n_{b}}^{b})\right)^{T},
	\end{equation}
	where $M_{\phi_{in}}\in\mathbb{R}^{n_{in}\times(n_{in}+n_{b})}$, $M_{\phi_{b}}\in\mathbb{R}^{n_{b}\times(n_{in}+n_{b})}$.
	
	Following \cite{chinchapatnam2006unsymmetric}, we define 
	$H_{-}:=\left(\begin{array}{c}
	M_{\phi_{in}}\\
		0
	\end{array}\right)$ 
	and putting together \eqref{eq:interporlation matrix}  and \eqref{eq:global scheme} obtain

	\begin{equation*}
		H_{+}\bar{\alpha}^{n+s}=H_{-}\sum_{k=0}^{s-1}\sigma_{k}\bar{\alpha}^{n+k}+\left(\begin{array}{c}
		\Delta t\,\boldsymbol{f}^{n+s}\\
		\boldsymbol{g}^{n+s}
	\end{array}\right).
	\end{equation*}

	\noindent
	Thus, it follows that
	\begin{equation*}
		\boldsymbol{y}^{n+s}	=AH_{+}^{-1}H_{-}A^{-1}\sum_{k=0}^{s-1}\sigma_{k}\boldsymbol{y}^{n+s}+AH_{+}^{-1}
		\left(\begin{array}{c}
		\Delta t\,\boldsymbol{f}^{n+s}\\
		\boldsymbol{g}^{n+s}
		\end{array}\right).
	\end{equation*}

	Denoting by $\mathbf{y}^{n}$ the exact solution and by $\hat{\mathbf{y}}^{n}$ 
	the numerically computed solution, the error 
	$\mathbf{e}^{n}=\mathbf{y}^{n}-\mathbf{\hat{y}}^{n}$ satisfies the equation
	 
	\begin{equation*}
		\mathbf{e}^{n+s}=K\sum\limits_{k=0}^{s-1}\sigma_{k}\mathbf{e}^{n+k}+E_{n+s},
	\end{equation*}
	where $E_{n+s}$ is the local error in the scheme \eqref{eq:global scheme}
	and $K=AH_{+}^{-1}H_{-}A^{-1}$. Besides, since $E_{n+s}$ is small and therefore bounded, 
	the error analysis can be analyzed using the equation
	\begin{equation*}\label{diference eq}
		\mathbf{e}^{n+s}=K\sum\limits_{k=0}^{s-1}\sigma_{k}\mathbf{e}^{n+k}.
	\end{equation*}
	By assuming that $K$ is diagonalizable, i.e.,  $K=D^{-1}\Lambda D$,  we can define $\mathbf{z}^{n}:=D\mathbf{e}^{n}$ 
	and therefore \eqref{diference eq} is equivalent to 
	\begin{equation*}\label{diference eq2}
		\mathbf{z}^{n+s}=\Lambda \sum\limits_{k=0}^{s-1}\sigma_{k}\mathbf{z}^{n+k}.
	\end{equation*}
	\noindent
	Since $\Lambda$ is a diagonal matrix, for every $j=1,\dots, d(n_b+n_{in})$ we have that 
	$z^{n+s}_{j}=\sum\limits_{k=0}^{s-1}\lambda_{j}\sigma_{k}z^{n+k}_{j}$, 
	and whose solution is given by $	z^{n}_{j}=\sum\limits_{k=0}^{s-1}C_{k}^{j}r_{k}^{n}$,
	where $C_{k}^{j}$ are arbitrary complex constants and $r_k$ are the roots of the associated polynomial 
	of the finite diference equation. 
		
	\noindent Finally, since $\|\mathbf{e}^{n}\|$ goes to zero iff $\|\mathbf{z}^{n}\|$ tends to zero, 
	the method will be stable as long as the eigenvalues of $K$ belong to the stability region of
	\begin{equation}\label{eq:region.stability}
		\pi(r,\lambda)=r^s - \sum_{k=0}^{s-1}\lambda\sigma_{k} r^{i}.
	\end{equation}
	
	As a consequence of the boundary locus technique \cite{lambert1991numerical}, we can deduce the following stability regions: 
	
	\begin{figure}[H]\label{fig.stencils}
	\begin{center}
		\subfloat[]{\includegraphics[scale=0.3]
		{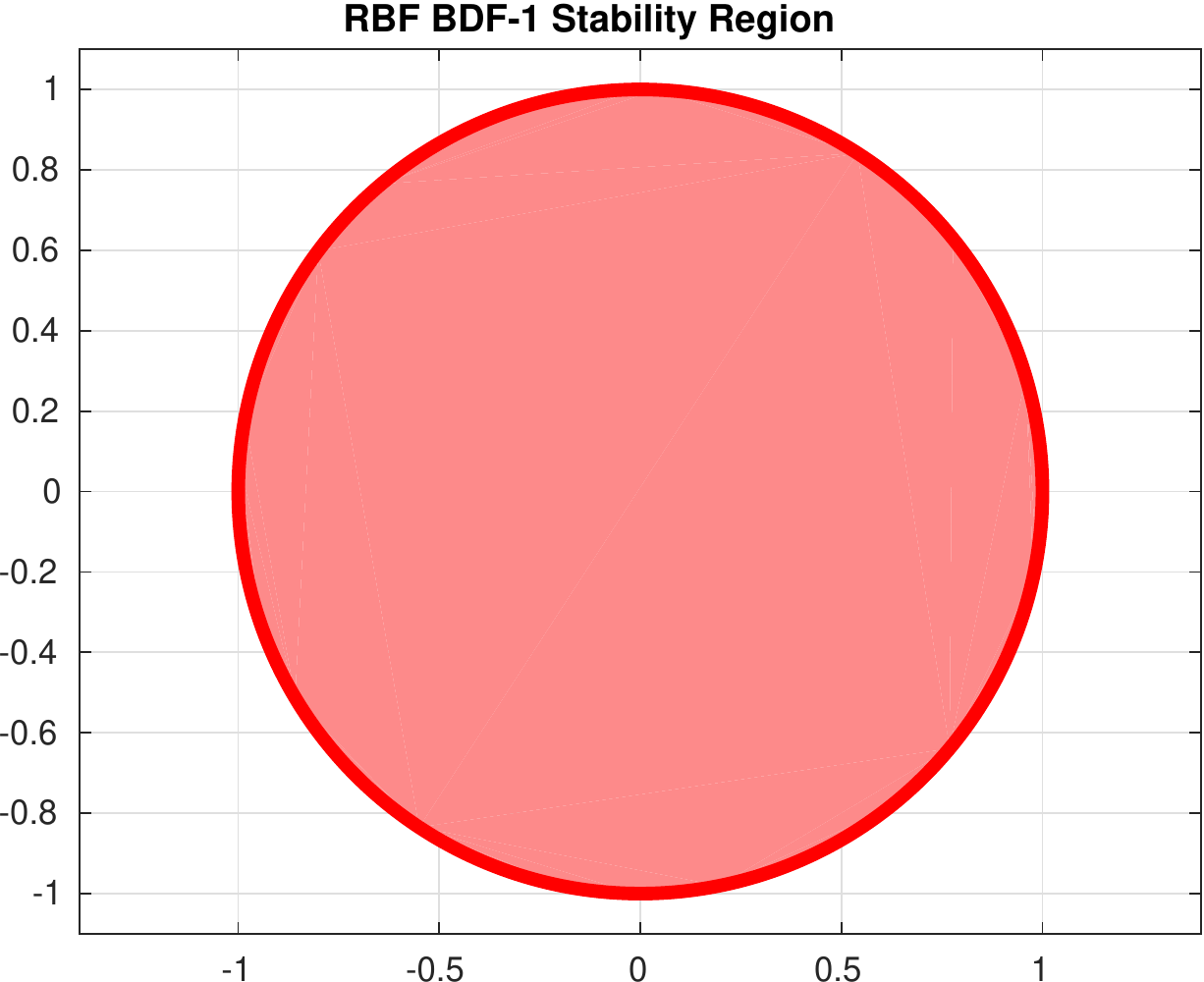}}
		\hspace{1.0cm} 
		\subfloat[]{\includegraphics[scale=0.3]
		{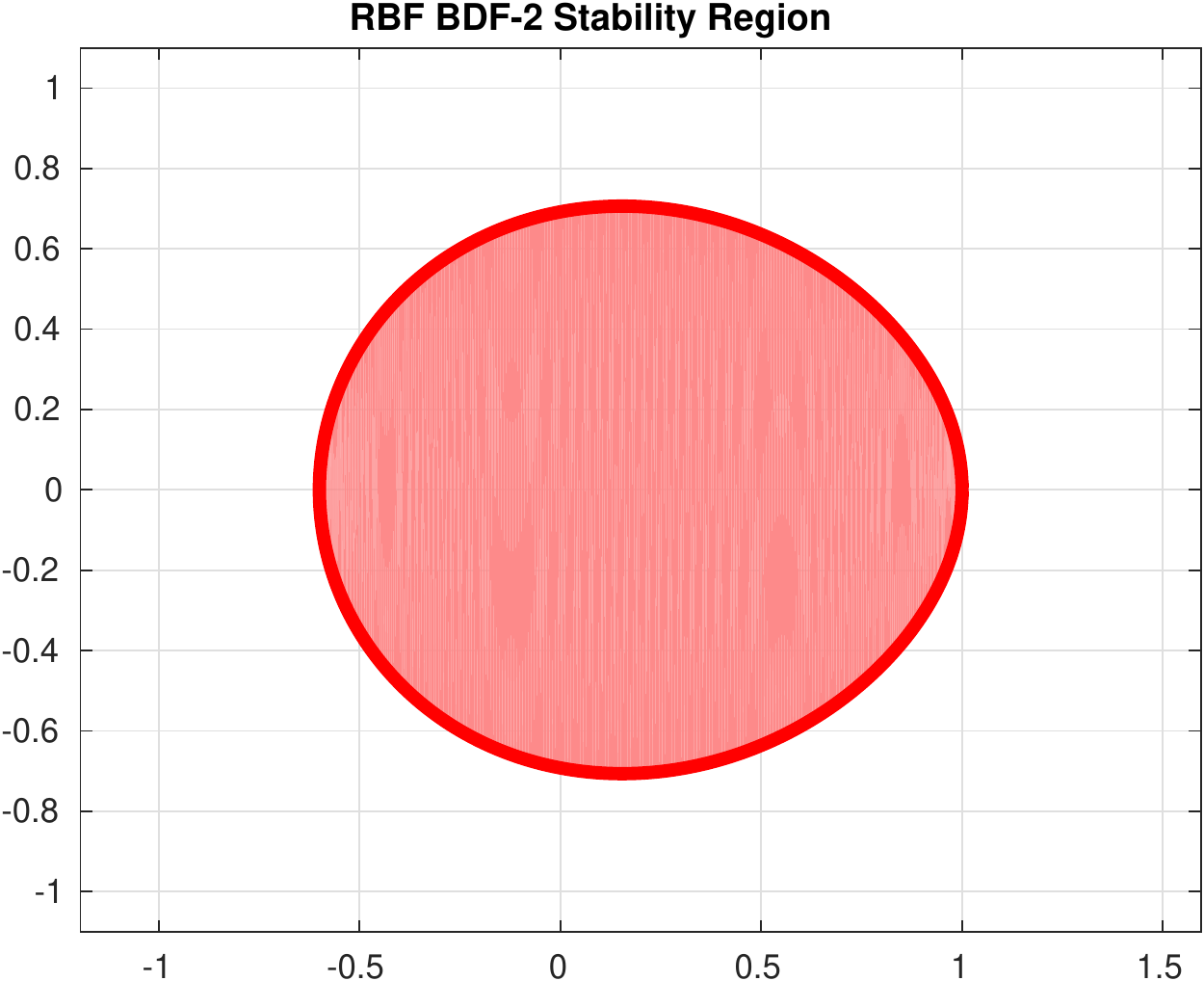}}
		\hspace{1.0cm} 
		\subfloat[]{\includegraphics[scale=0.3]
		{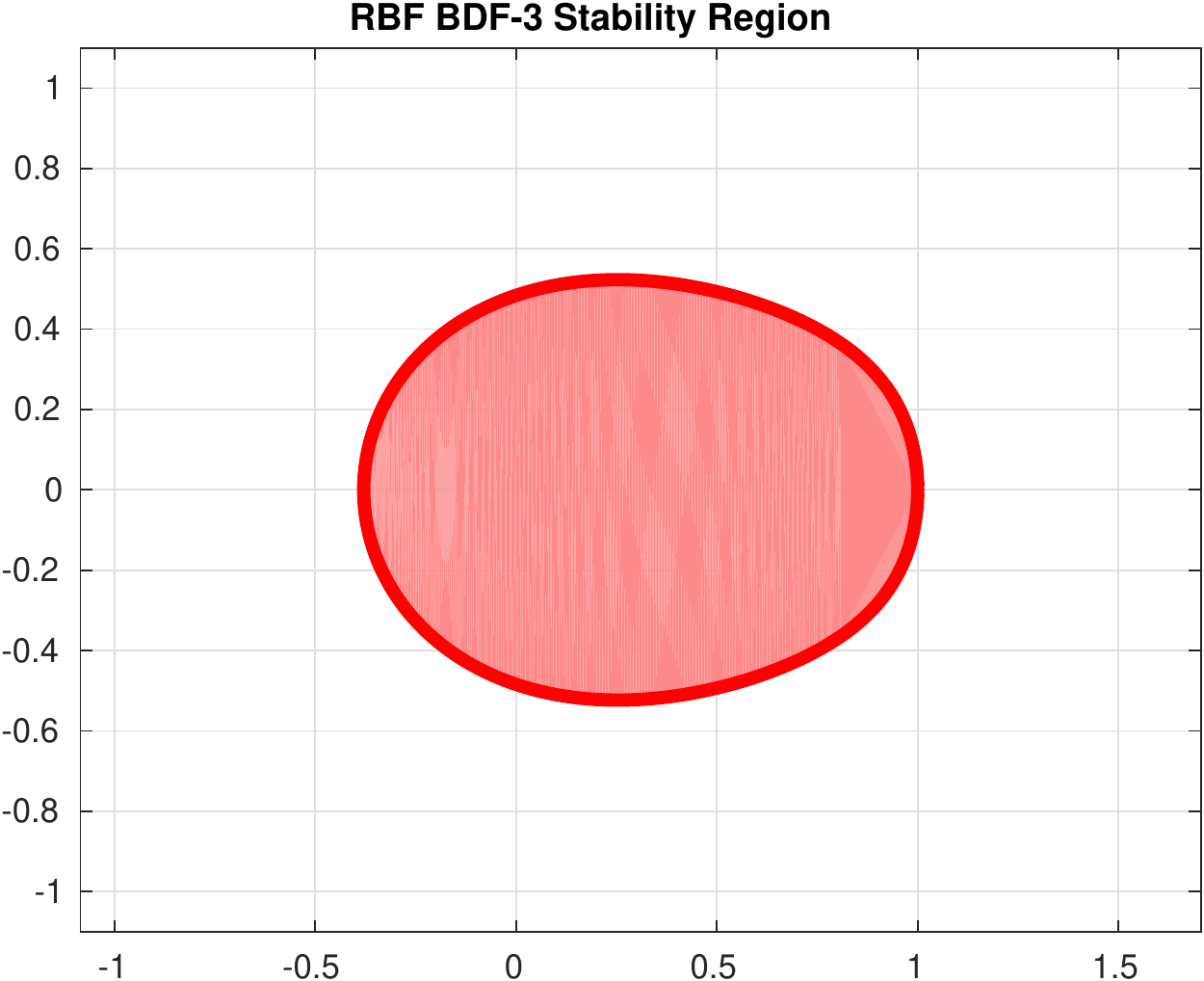}}
	\end{center}
	\par
	\caption{Stability regions using backward finite differences (BDF): a) nivel one, b) two levels and c) three levels.}
	\end{figure}
	

\subsection{\normalsize{The method of lines: global collocation}}\label{subsection.line.method}
	From the ansatz for the steady Stokes system, \eqref{Anzatdist.Stokes.system}, we develop a 
	numerical approach for the unsteady Stokes system with either Navier--slip or Dirichlet boundary conditions
	based in the method of lines to separate the variables in our approximation. Then, following the scheme 
	carried out in the previous section, our ansatz is defined by:
	\begin{equation*}
		(\hat{\mathbf y},\hat{p})(\boldsymbol{x,t})=\sum_{i=1}^{d}\sum_{j=1}^{n_{b}}\left(B_{ij}\right)
		\Phi(\boldsymbol{x-z})\alpha_{ij}^{B}(t)+\sum_{i=1}^{d}\sum_{j=1}^{n_{in}}\left(L_{ij}\right)
		\Phi(\boldsymbol{x-z})\alpha_{ij}^{L}(t),
	\end{equation*}
	with certain coefficients 
	$\overline{\alpha}(t)=(\alpha_{ij}^B(t),\alpha_{ij}^L(t))\in\mathbb{R}^{d(n_{b}+n_{in})}$, 
	$L_{ij}$ is  given by
	\begin{equation*}
		L_{ij}=\delta_{\boldsymbol{x_{j}}}\circ -\mu \Delta e_{i}+\delta_{\boldsymbol{x_{j}}}\circ\partial_{i}e_{d+1},
	\end{equation*}
	meanwhile $B_{ij}$ is defined as in \eqref{DFG.ansatz.slip.functionals2}. In a similar way to Section  
	\ref{subsection.operator.discre}, we develop this method in 2D and consider the vector value function 
	$(\phi^{B_{1}},\phi^{B_{2}},\phi^{L_{1}},\phi^{L_{2}})$ satisfying 
	
	\begin{eqnarray*} \label{eq:Anzat estacionario global}
		(\hat{\mathbf y},\hat{p})(\boldsymbol{x,t}) 
		=&&\hspace{-0.6cm} \sum_{j=1}^{n_{b}}\phi^{B_{1}}(\boldsymbol{x}-\boldsymbol{x_{j}^{b}})
		\alpha_{1j}^{B}(t)+\sum_{j=1}^{n_{b}}\phi^{B_{2}}(\boldsymbol{x}-\boldsymbol{x_{j}^{b}})\alpha_{2j}^{B}(t)\\
 		&& +  \sum_{j=1}^{n_{in}}\phi^{L_{1}}(\boldsymbol{x}-\boldsymbol{x_{j}^{in}})\alpha_{1j}^{L}(t)+\sum_{j=1}
 		^{n_{in}}\phi^{L_{2}}(\boldsymbol{x}-\boldsymbol{x_{j}^{in}})\alpha_{2j}^{L}(t). \nonumber
	\end{eqnarray*}
	
	\noindent Putting the ansatz above in the Stokes system,  we obtain the following scheme
	\begin{equation}\label{eq:stokes-evolutivo}
		\left\{ \begin{array}{lll}
		\mathbf{\hat{y}}_{t}+L(\mathbf{\hat{y}},\hat{p})=\mathbf{f} & \text{ in } & Q,\\
		\mathbf{B\hat{y}}=\mathbf{g} & \text{ on } & \Sigma,\\
		\hat{y}(\cdot,0)=\mathbf{y}_{0}(\cdot) & \text{ in } & Q,
		\end{array}\right.
	\end{equation}
	which implies  the system of ODEs
	\begin{align*}\label{eq:ode global line system} 
		M_{\phi}\frac{\overline{\alpha}(t)}{dt}+M_{L\phi}\overline{\alpha}(t)=\mathbf{f},\\
		M_{B_{\phi}}\overline{\alpha}(t)	=\mathbf{g},\\ 
	\end{align*}
	where, for $\phi=[\phi_1,\phi_2]$,\, we have that $M_{\phi},M_{L_{\phi}}\in \mathcal{M}((2n_{in})\times 
	(2n_{b}+2n_{in});\mathbb{R})$ and are defined as follows:
	\begin{equation*}
		M_{\phi}=\left(
		\begin{array}{ccccc}
		\phi_{1}^{B_{1}} & \phi_{1}^{B_{2}} &  \phi_{1}^{L_{1}} & \phi_{1}^{L_{2}}\\
		\phi_{2}^{B_{1}} & \phi_{2}^{B_{2}} &  \phi_{2}^{L_{1}} & \phi_{2}^{L_{2}}
		\end{array}\right),
		\,\,\,
		M_{L_{\phi}}=\left(\begin{array}{ccccc}
		L_{1}\phi^{B_{1}} & L_{1}\phi^{B_{2}} & L_{1}\phi^{L_{1}} & L_{1}\phi^{L_{2}}\\
		L_{2}\phi^{B_{1}} & L_{2}\phi^{B_{2}} & L_{2}\phi^{L_{1}} & L_{2}\phi^{L_{2}}
		\end{array}\right)
	\end{equation*}
	and $M_{B_{\phi}}\in  \mathcal{M}(2n_{b}\times(2n_{b}+2n_{in});\mathbb{R})$
	\begin{equation*}
		M_{B_{\phi}}=\left(\begin{array}{ccccc}
		B_{1}\phi^{B_{1}} & B_{1}\phi^{B_{2}} & B_{1}\phi^{L_{1}} & B_{1}\phi^{L_{2}}\\
		B_{2}\phi^{B_{1}} & B_{2}\phi^{B_{2}} & B_{2}\phi^{L_{1}} & B_{2}\phi^{L_{2}}
		\end{array}\right).
	\end{equation*}
	
	\noindent Finally, we use BDF methods for solving the former system of ODEs, thus,  at each step we solve the following 
	linear system:

	\begin{equation}\label{eq:bf2 gbobal}
	\left(\begin{array}{c}
		M_{\phi}+\Delta t\,\beta_{s}\,M_{L\phi}\\
		M_{B_{\phi}}
	\end{array}\right)\overline{\alpha}^{n+s}=\left(\begin{array}{c}
	\Delta t\,\beta_{s}\,\mathbf{f}^{n+s}+\sum\limits_{k=0}^{s-1}\sigma_{k}M_{\phi}\overline{\alpha}^{n+k}\\
	\mathbf{g}^{n+s}
	\end{array}\right).
	\end{equation}

	\noindent Since $M_{\phi}\overline{\alpha}(t_n)=	\mathbf{y}^n$, the equation \eqref{eq:bf2 gbobal} is equivalent to

	\begin{equation}\label{eq:bf2_gbobal2}
		\left(\begin{array}{c}
		M_{\phi}+\Delta t\,\beta_{s}\,M_{L\phi}\\
		M_{B_{\phi}}
		\end{array}\right)\overline{\alpha}^{n+s}=\left(\begin{array}{c}
		\Delta t\,\beta_{s}\,\mathbf{f}^{n+s}+\sum\limits_{k=0}^{s-1}\sigma_{k}\mathbf{y}^{n+k}\\
		\mathbf{g}^{n+s}
		\end{array}\right).
	\end{equation}
	
	\begin{Obs}\label{obs.stability.line.method}
	As in subsection \ref{subsection.stability.global1},  we proceed to present an analysis of the stability
	for the method introduced in this subsection. From system \eqref{eq:bf2_gbobal2},  we can deduce 

	\begin{equation*}
		\boldsymbol{y}^{n+s}	=AH_{+}^{-1}H_{-}A^{-1}\sum_{k=0}^{s-1}\sigma_{k}\boldsymbol{y}^{n+s}+AH_{+}^{-1}
		\left(\begin{array}{c}
		\Delta t\,\boldsymbol{f}^{n+s}\\
		\boldsymbol{g}^{n+s}
		\end{array}\right).
	\end{equation*}

	\noindent
Where $A$ is defined as in equation (\ref{eq:interporlation matrix}) and : 	

\[
H_{+}=\left(\begin{array}{c}
M_{\phi}+\Delta t\,\beta_{s}\,M_{L\phi}\\
M_{B_{\phi}}
\end{array}\right)\quad H_{-}=\left(\begin{array}{c}
M_{\phi}\\
0
\end{array}\right)
\]

		\noindent
	Therefore, the method will be stable as long as the eigenvalues of $K=AH_{+}^{-1}H_{-}A^{-1}$ are in the stability region of \eqref{eq:region.stability}. 
	\end{Obs}


	\subsection{\normalsize{Numerical examples: line method BDF2}} \label{sec:global.numerical.results.line}
	In the sequel, we evaluate the accuracy of the previous scheme through BDF2. The objective is test the feasibility 
	of the schemes by considering either nonhomogeneous Dirichlet conditions or nonhomogeneous Navier--slip conditions 
	on the system \eqref{Stokes.system}. In order to generate the divergence 
	free kernel we use inverse multi quadric (MQ), with a shape parameter $c=0.1$. Since this type of kernel are very 
	ill conditioned we have used the Matlab package ADVANPIX for multi precision calculus and set the number of digits 
	to 50. In this subsection, all the computations are done in the programing languages Matlab and FreeFemm++. 
	We introduce a uniform mesh generate by means of the package DISTMESH, where the total number of nodes is 362. 
	
\begin{figure}[H]
	\begin{center}
		\includegraphics[scale=0.18]{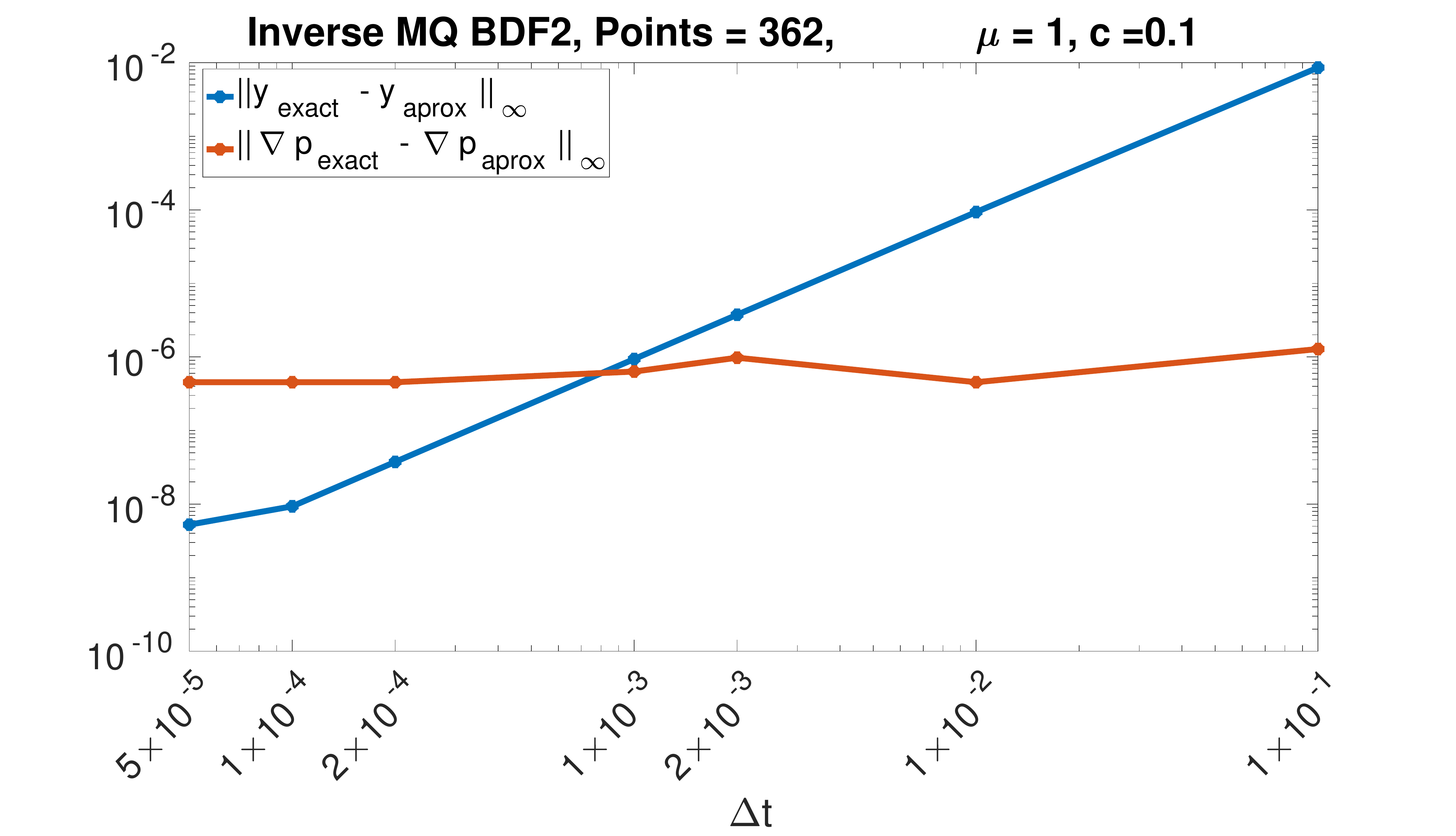}
		\hspace{0.01cm} 
		\includegraphics[scale=0.18]{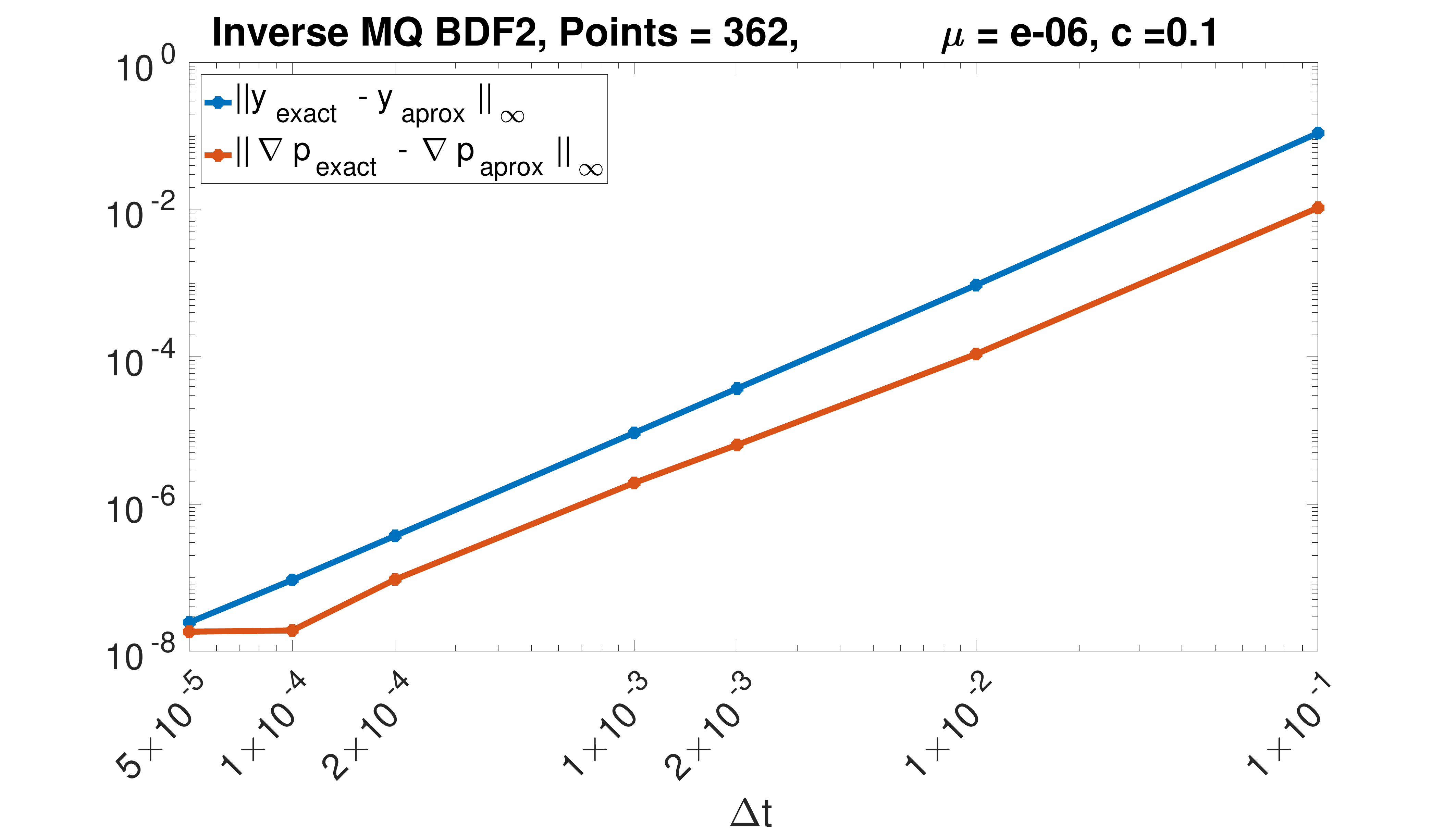}
	\end{center}
		\par
	\caption{Convergence behaviour of $\left\|\boldsymbol{\epsilon(y)}\right\|_{\infty}$ (blue) and 
	$\left\|\boldsymbol{\nabla\epsilon(p)}\right\|_{\infty}$ (red) in the 2D test,  Dirichlet boundary conditions.
	The left--hand side with $\mu=1$, meanwhile in the right--hand side $\mu=1\times 10^{-6}$.}
	\label{fig.time_vs_error.globalmethod}
	\end{figure}

	Let $\Omega$ be the unit circle, i.e.,
	$\Omega=\{(x,y)\in\mathbb{R}^2: x^2+y^2<1\}$, the analytical solution of \eqref{Stokes.system} given by
	$$(\mathbf{y_1}(x,y,t),\mathbf{y_2}(x,y,t))=\Bigl(-\pi\,y\,\sin\left(\frac{\pi}{2}(x^{2}+y^{2})\right)\sin(\pi t), 
	\pi\,x\,\sin\left(\frac{\pi}{2}(x^{2}+y^{2})\right)\sin(\pi t)\Bigr)$$ 
	and pressure $p(x,y,t)=\sin(x-y+t)$. 
	We compare the error in velocity and pressure in the $L^{\infty}$--norm between the exact and numerical 
	solutions for several time steps $\Delta t$. The errors are denoted by 
	$\boldsymbol{\epsilon_y}=\mathbf{y}_{exact}-\mathbf{y}_{aprox}$ and 
	$\boldsymbol{\nabla\epsilon_p}=\nabla p_{exact}-\nabla p_{aprox}$,  respectively. The results are presented in Table 
	\ref{table.difusion_1} and Figure \ref{fig.time_vs_error.globalmethod}.

\begin{table}[H]
\begin{centering}
\begin{tabular}{ccccc}
\hline 
B.C  & \multicolumn{2}{c}{Dirichlet} & \multicolumn{2}{c}{Navier-slip}\tabularnewline
\hline 
$\Delta t$  &\hspace{1cm} $\left\Vert \boldsymbol{\epsilon_{y}}\right\Vert_{\infty} $  &\hspace{1cm} $\left\Vert \boldsymbol{\nabla\epsilon_{p}}\right\Vert_{\infty} $  &
\hspace{1cm} $\left\Vert \boldsymbol{\epsilon_{y}}\right\Vert_{\infty} $  &\hspace{1cm} $\left\Vert \boldsymbol{\nabla\epsilon_{p}}\right\Vert_{\infty} $\tabularnewline
\hline 
 & \multicolumn{4}{c}{$\mu=1$}\tabularnewline
\cline{2-5} 
0.1  &\hspace{1cm} 8.570E-03  &\hspace{1cm} 1.28E-06  &\hspace{1cm} 9.23E-02 &\hspace{1cm} 1.42E-04\tabularnewline
0.01  &\hspace{1cm}  9.32E-05  &\hspace{1cm} 4.53E-07  &\hspace{1cm} 8.00E-04 &\hspace{1cm} 1.43E-04\tabularnewline
0.002  &\hspace{1cm}  3.74E-06  &\hspace{1cm} 9.76E-07  &\hspace{1cm} 4.39E-05 &\hspace{1cm} 1.43E-04\tabularnewline
0.001  &\hspace{1cm} 9.37E-07  &\hspace{1cm} 6.33E-07  &\hspace{1cm} 2.36E-05 &\hspace{1cm} 1.43E-04\tabularnewline
0.0002  &\hspace{1cm} 3.74E-08  &\hspace{1cm} 4.53E-07  &\hspace{1cm} 1.97E-05 &\hspace{1cm} 1.43E-04\tabularnewline
0.0001  &\hspace{1cm} 9.32E-09  &\hspace{1cm} 4.53E-07  &\hspace{1cm} 1.96E-05 &\hspace{1cm} 1.43E-04\tabularnewline
5e-05  &\hspace{1cm} 5.24E-09  &\hspace{1cm} 4.53E-07  &\hspace{1cm} 1.96E-05 &\hspace{1cm} 1.43E-04\tabularnewline
 & \multicolumn{2}{c}{$\mu=10^{-6}$} & \multicolumn{2}{c}{$\mu=10^{-3}$}\tabularnewline
\cline{2-5} 
0.1  &\hspace{1cm} 1.10E-01  &\hspace{1cm} 1.06E-02  &\hspace{1cm} 1.10E-01 &\hspace{1cm} 3.61E-04\tabularnewline
0.01  &\hspace{1cm} 9.48E-04  &\hspace{1cm} 1.09E-04  &\hspace{1cm} 9.79E-04 &\hspace{1cm} 3.57E-04\tabularnewline
0.002  &\hspace{1cm} 3.72E-05  &\hspace{1cm} 6.38E-06  &\hspace{1cm} 1.48E-04 &\hspace{1cm} 3.81E-04\tabularnewline
0.001  &\hspace{1cm} 9.29E-06  &\hspace{1cm} 1.94E-06  &\hspace{1cm} 1.41E-04 &\hspace{1cm} 3.86E-04\tabularnewline
0.0002  &\hspace{1cm} 3.70E-07  &\hspace{1cm} 9.45E-08  &\hspace{1cm} 1.39E-04 &\hspace{1cm} 3.86E-04\tabularnewline
0.0001  &\hspace{1cm} 9.31E-08  &\hspace{1cm} 1.91E-08  &\hspace{1cm} 1.39E-04 &\hspace{1cm} 3.85E-04\tabularnewline
5e-05  &\hspace{1cm} 2.46E-08  &\hspace{1cm} 1.84E-08  &\hspace{1cm} 1.39E-04 &\hspace{1cm}  3.85E-04\tabularnewline
\hline 
\end{tabular}

\par\end{centering}
\caption{Convergence behaviour of velocity and pressure in the $L^{\infty}$\textendash norm
using inverse MQ, 362 nodes and different viscosity coefficients $\mu$.}\label{table.difusion_1} 
\end{table}

	In the left--hand side of figure \ref{fig.time_vs_error.globalmethod}  can be observed that the pressure remains almost constant, this result was obtained by Fuselier 
	in their experiments, \cite{fuselier2016high}. A posible explanation to this behavior is that the convergence of the pressure, which we think it should depend on the fill 
	distance, the shape parameter and the viscosity is lower bounded. In fact,  in the right--hand side of  figure  \ref{fig.time_vs_error.globalmethod} where the viscosity $\mu$
	is of the order of $10^{-6}$,  the pressure decay with almost the same slope as the velocity and it becomes constant after the time step less than $10^{-4}$. 
	
	\begin{figure}[H]
	\begin{center}
		\includegraphics[scale=0.18]{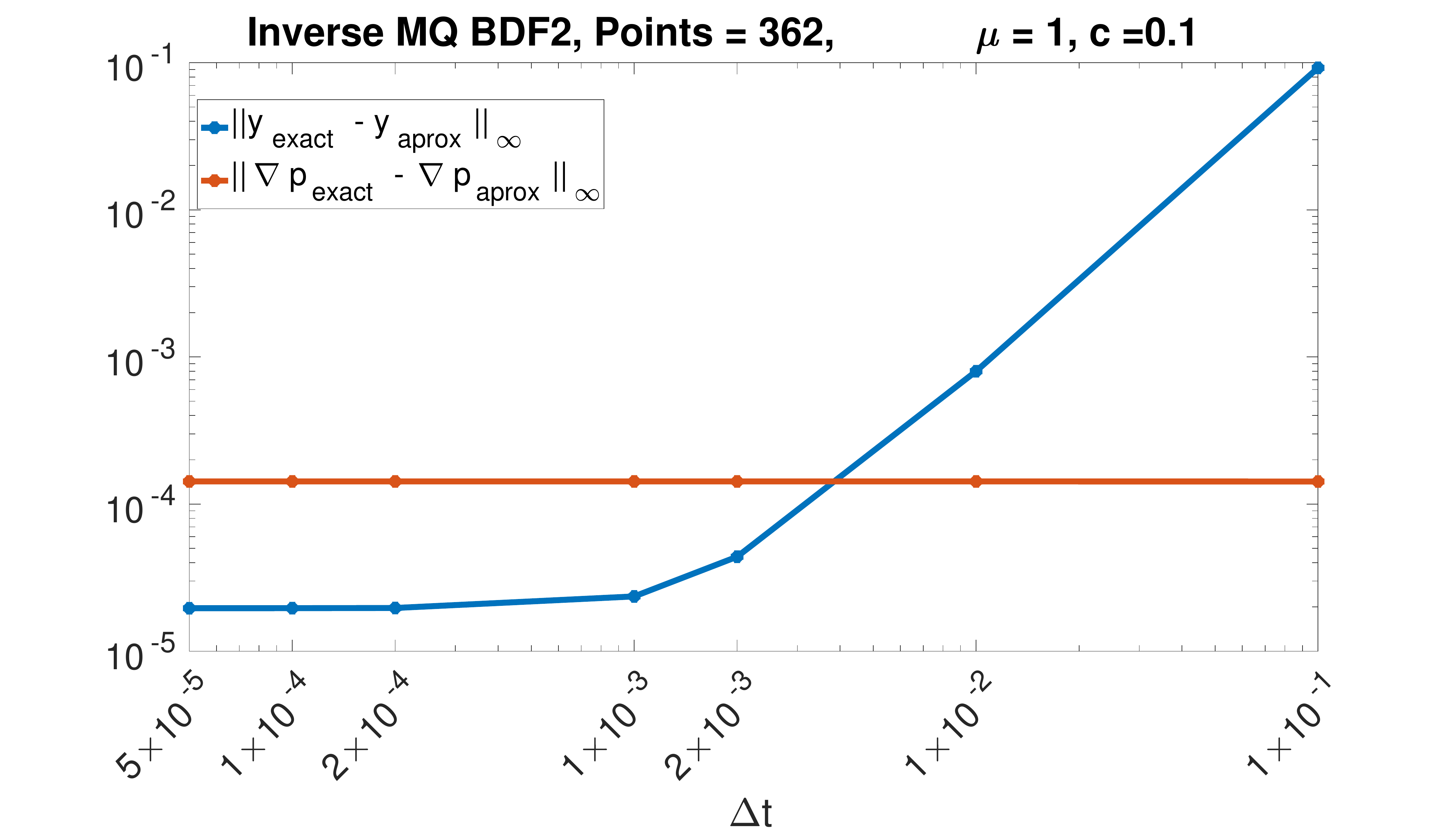}
		\hspace{0.01cm} 
		\includegraphics[scale=0.18]{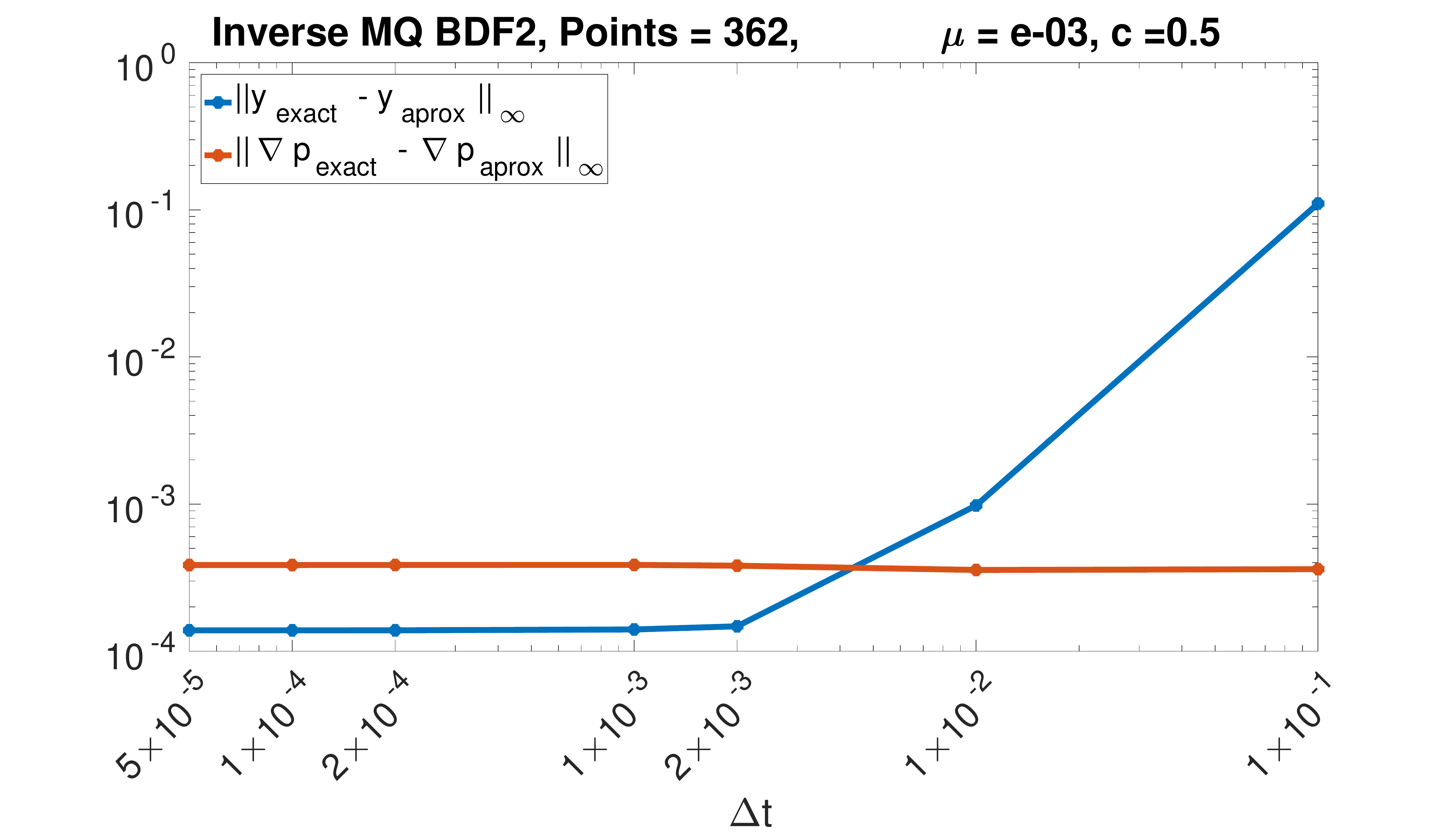}
	\end{center}
		\par
	\caption{Convergence behaviour of $\left\|\boldsymbol{\epsilon(y)}\right\|_{\infty}$ (blue) and 
	$\left\|\boldsymbol{\nabla\epsilon(p)}\right\|_{\infty}$ (red) by considering  Navier--slip boundary conditions and two different values of viscosity $\mu$ and shape 
	parameter $c$.}
	\label{fig.time_vs_error.globalmethod2}
	\end{figure}
\begin{figure}[H]
	\begin{center}
		\includegraphics[scale=0.22]{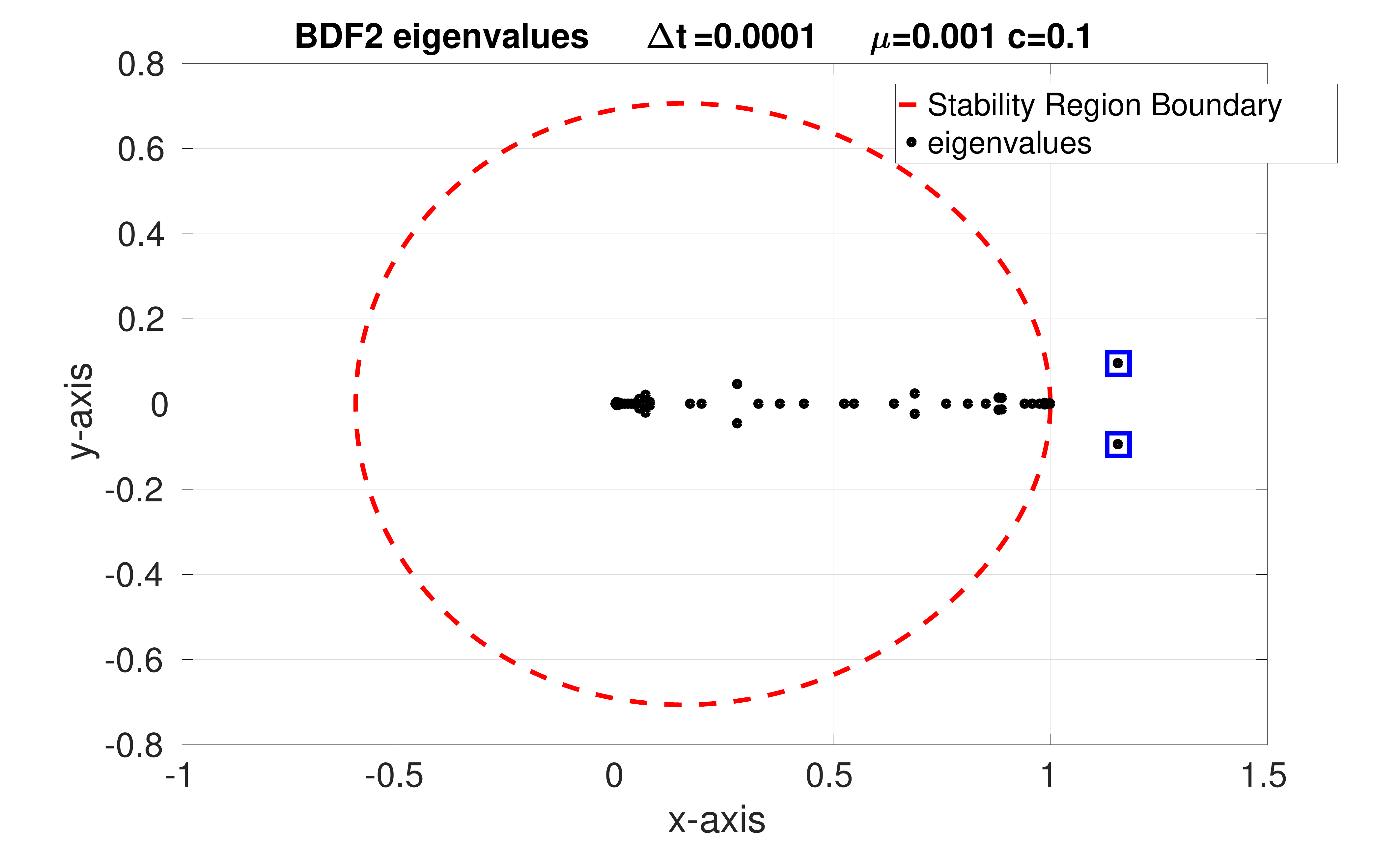}
		\end{center}
	\caption{Eigenvalues of the matrix $K=AH_+^{-1}H_{-} A^{-1}$ defined in subsection \ref{subsection.stability.global1}  concerning stability analysis.The blue squares represent the eigenvalues located outside the stability region. Here, Navier-Slip	boundary conditions were imposed.}	
	\label{fig.eigen_value_plot}
	\end{figure}

	\vskip 0.4cm
	\noindent
	Figure \ref{fig.time_vs_error.globalmethod2} shows the convergence behaviour of  $\left\|\boldsymbol{\epsilon(y)}\right\|_{\infty}$ and 
	$\left\|\boldsymbol{\nabla\epsilon(p)}\right\|_{\infty}$ by putting Navier--slip boundary conditions. Observe that different bounds are obtained in this case,  depending again 
	on the fill distance, the shape parameter and the viscosity. Indeed, if we change the shape parameter for the right figure by  $c=0.1$,  the convergence is not reached.
	It may be a hidden  condition of the style of CFL condition. Numerically, this can be explained by noting that some eigenvalues lies outside the stability region, see figure 		\ref{fig.eigen_value_plot}.

\section{\normalsize{LHI method for the Stokes system}}\label{section_LHI_Stokes}
In this section we formulate a RBF--LHI vectorial technique for the unsteady
Stokes system based on the scalar LHI algorithm, see introduction. First, the numerical
description for the steady Stokes system is formulated and a numerical example 
is presented. Once completed the steady case, we apply the time advancement scheme given in subsection 
\ref{subsection.operator.discre} for obtaining an implicit discretization scheme. A stability analysis is
also carried out.
	   
\subsection{\normalsize{Steady problem}}
	From \cite{wendland2009divergence}, the vectorial LHI-algorithm for the system 
	\begin{equation}
		\left\{ \begin{array}{lll}
		-\mu \Delta\mathbf{y}+\nabla p=\mathbf{F} & \text{ in } & Q,\\
		\nabla\cdot\mathbf{y}=0 & \text{ in } & Q,\\
		+BC & \text{ on } & \Sigma,
		\end{array}\right.\label{eq:Stokes.system}
	\end{equation}
	where $\mathbf{F}$ is known and $+BC$ indicates the boundary conditions (see \eqref{intro.BC.Stokes}),  
	reads as follows: let us define a vector 
	$\mathbf{u}:=(\mathbf{u}_{n},u_{n+1})=(\mathbf{y},p)\in\mathbb{R}^{d+1}$.
	and operators $\boldsymbol{{\cal {L}}}, \boldsymbol{\cal{B}}, \boldsymbol{\cal{B}}_1,\boldsymbol{\cal{B}}_2$ such 
	that the right--hand side of \eqref{eq:Stokes.system} is given by
	\begin{equation*}
		\boldsymbol{{\cal {L}}}\mathbf{u}:=\mathbf{-\mu \Delta\mathbf{u}_{n}}+\nabla u_{n+1}=(-\mu \Delta,\nabla)
		\cdot(\mathbf{u}_{n},u_{n+1});\quad\nabla\cdot\mathbf{u}_{n}=0
	\end{equation*}
	and either
	\begin{equation}\label{boundary.cond.numerics.Stokes}
		\underbrace{\boldsymbol{{\cal {B}}}\mathbf{u}:=\mathbf{u}_{n}}_{\mbox{Dirichlet}}
		\quad\mbox{or}\quad
		\underbrace{\boldsymbol{\cal{B}}_1 \mathbf{u}:=(\sigma(\mathbf{u}_n,u_{n+1})\boldsymbol{\nu})_{tg} 
		=(\sigma(\mathbf{u})\boldsymbol{\nu})_{tg};\quad \boldsymbol{\cal{B}}_2 \mathbf{u}
		=\mathbf{u}_n\cdot \boldsymbol{\nu}.}_{\mbox{Navier--slip}}	
	\end{equation}
	
	\noindent	
    As previously, for each disk $D^{k}$ with center $x_{1}^{k}$, we must define a local sub-system
    (see introduction):
     
	\begin{equation}\label{eq:local.system.steasy}
		\left\{ \begin{array}{lll}
		\mathbf{u}(x_{i})=\overline{h_{i}} & \text{ in } & \Omega_{sc}^{k} \subset \Omega_{sc}\bigcap D^{k},\\
		\\
		\boldsymbol{{\cal {L}}}\mathbf{u}=\mathbf{F} & \text{ in } & \Omega_{pdec}^{k} \subset \Omega_{pdec}\bigcap D^{k},
		\\
		\\
		\eqref{boundary.cond.numerics.Stokes} & \text{ on } & \Omega_{fc}^{k} \subset \Omega_{fc}\bigcap D^{k}.
		\end{array}\right.
	\end{equation}
	In order to solve the systems \eqref{eq:local.system.steasy} we first define the matrix-valued kernel
	
	\begin{equation*}
		\mathbf{\Phi}=\left[{\begin{array}{cc}
		\mathbf{\Phi}_{Div} & 0\\
		0 & \psi
		\end{array}}\right]:\mathbb{R}^{d}\rightarrow\mathbb{R}^{(d+1)\times(d+1)},
	\end{equation*}
	where $\mathbf{\Phi}_{Div}=\Delta\times\nabla\times\psi=\{-\Delta I+\nabla\nabla^{T}\}\psi(\mathbf{x}):
	\mathbb{R}^{d}\rightarrow\mathbb{R}^{d+1}$
	is a divergence-free positive definite kernel, $\Delta$ is the Laplacian,
	$I$ the identity matrix, and  $\phi,\,\psi$ are global $C^{\infty}$ positive
	definite scalar RBFs.

	Since we are choosing free divergence radial kernel, the incompressibility equation is missing, 
	and therefore we lack a differential operator. In order to obtain a fully sparse system from \eqref{eq:calculo de pesos}, we  include the pressure  as unknown variable 
	in the local system \eqref{eq:local.system.steasy}. This leads  to the following local system 
	\begin{equation}\label{lhi:local.system.implicit}
		\left\{ \begin{array}{lll}
		\mathbf{u_n}(x_{i})=\overline{h_{i}} & \text{ in } & \Omega_{sc}^{k} ,\\
		\boldsymbol{{\cal {L}}}\mathbf{u}=\mathbf{F} & \text{ in } & \Omega_{pdec}^{k},
		\\
		\eqref{boundary.cond.numerics.Stokes} & \text{ on } & \Omega_{fc}^{k}.
		\end{array}\right.
	\end{equation}
	
   \noindent 
	Using the generalized interpolation collocation method \cite{wendland2004scattered},  
	the ansatz by considering Navier--slip boundary conditions (see \eqref{boundary.cond.numerics.Stokes}) is given by

	\begin{eqnarray}\label{lhi:Anzatz.stokes.stacionary}
		(\hat{y}^{(k)},\hat{p}^{(k)})(\boldsymbol{x}) & = & 
		\sum_{i=1}^{d}\sum_{j=1}^{n_{sc}^{k}}\left(I_{ij}^{k}\right)
		\mathbf{\Phi}(\boldsymbol{x-z})\alpha_{ij}^{I^{k}}
		 +  \sum_{i=1}^{d+1}\sum_{j=1}^{n_{fc}^{k}}\left(B_{ij}^{k}\right)\mathbf{\Phi}(\boldsymbol{x-z})
		\alpha_{ij}^{B^{k}}\nonumber \\
		&  & \quad +\sum_{i=1}^{d}\sum_{j=1}^{n_{pdec}^{k}}\left(L_{ij}^{k}\right)\mathbf{\Phi}(\boldsymbol{x-z})
		\alpha_{ij}^{L^{k}},\nonumber 
	\end{eqnarray}
	where $\overline{\alpha}=(\alpha_{ij}^{I^{k}},\alpha_{ij}^{B^{k}}, \alpha_{ij}^{L^{k}})
	\in\mathbb{R}^{d(n_{sc}^{k}+n_{fc}^{k}+n_{pdec}^{k})}$,\,   
	$I_{ij}^{k}:=\delta_{\boldsymbol{x_{j}^{sc(k)}}}e_{i}$,\, $B_{ij}^k$
	are vector--valued functionals defined like in \eqref{DFG.ansatz.slip.functionals2}, and 
	$L_{ij}^{k}:=\delta_{\boldsymbol{x_{j}^{sc(k)}}}\circ-\mu\Delta
	e_{i}+\delta_{\boldsymbol{x_{j}^{sc(k)}}}\circ\partial_{i}e_{d+1}$.
	 
	\noindent Now, in order to simplify the notation, we shall proceed our analysis by defining the vector value
	function 
	$(\phi^{I_{1}},\phi^{I_{2}},\phi^{B_{1}},\phi^{B_{2}},\phi^{L_{1}},\phi^{L_{2}})$,
	such that the ansatz above can be rewritten in the form
	
	{\small{}
	\begin{align}\label{lhi:2d.local.anzat}
		(\hat{y}^{(k)},\hat{p}^{(k)})(\boldsymbol{x}) & =
		\sum_{j=1}^{n_{sc}^{k}}\phi^{I_{1}}(\boldsymbol{x}-\boldsymbol{x_{j}^{sc(k)}})\alpha_{1j}^{I}
		+
		\sum_{j=1}^{n_{sc}^{k}}\phi^{I_{2}}(\boldsymbol{x}-\boldsymbol{x_{j}^{sc(k)}})\alpha_{2j}^{I} \\ 
 		& +	
 		\sum_{j=1}^{n_{fc}^{k}}\phi^{B_{1}}(\boldsymbol{x}-\boldsymbol{x_{j}^{fc(k)}})\alpha_{1j}^{B}
		 +
		 \sum_{j=1}^{n_{fc}^{k}}\phi^{B_{2}}(\boldsymbol{x}-\boldsymbol{x_{j}^{fc(k)}})\alpha_{2j}^{B}\nonumber \\
 		 & +
 		\sum_{j=1}^{n_{pdec}^{k}}\phi^{L_{1}}(\boldsymbol{x}-\boldsymbol{x_{j}^{pdec(k)}})\alpha_{1j}^{L}
 		+
 		\sum_{j=1}^{n_{pdec}^{k}}\phi^{L_{2}}(\boldsymbol{x}-\boldsymbol{x_{j}^{pdec(k)}})\alpha_{2j}^{L}.\nonumber 
	\end{align}
	}{\small \par}
	\noindent
	Putting ansatz (\ref{lhi:2d.local.anzat}) in the local Stokes system (\ref{lhi:local.system.implicit}), we obtain the 
	local Gram matrix
	\begin{align}
		A^{(k)}=\left(\begin{array}{ccccccc}
		I_{1}\phi^{I_{1}} & I_{1}\phi^{I_{2}} & I_{1}\phi^{B_{1}} & I_{1}\phi^{B_{2}}  & I_{1}\phi^{L_{1}}
		 & I_{1}\phi^{L_{2}}\\
		I_{2}\phi^{I_{1}} & I_{2}\phi^{I_{2}} & I_{2}\phi^{B_{1}} & I_{2}\phi^{B_{2}}  & I_{2}\phi^{L_{1}} 
		& I_{2}\phi^{L_{2}}\\
		B_{1}\phi^{I_{1}} & B_{1}\phi^{I_{2}} & B_{1}\phi^{B_{1}} & B_{1}\phi^{B_{2}}  & B_{1}\phi^{L_{1}} 
		& B_{1}\phi^{L_{2}}\\
		B_{2}\phi^{I_{1}} & B_{2}\phi^{I_{2}} & B_{2}\phi^{B_{1}} & B_{2}\phi^{B_{2}}  & B_{2}\phi^{L_{1}}
		 & B_{2}\phi^{L_{2}}\\
		L_{1}\phi^{I_{1}} & L_{1}\phi^{I_{2}} & L_{1}\phi^{B_{1}} & L_{1}\phi^{B_{2}}  & L_{1}\phi^{L_{1}}
		& L_{1}\phi^{L_{2}}\\
		L_{2}\phi^{I_{1}} & L_{2}\phi^{I_{2}} & L_{2}\phi^{B_{1}} & L_{2}\phi^{B_{2}} & L_{2}\phi^{L_{1}} 
		& L_{2}\phi^{L_{2}}
	\end{array}\right),
	\end{align}
	which in turn let us to compute the weights by solving the following system
	
	\begin{equation*}
	\begin{cases}
		A^{(k)}W_{L_{1}}^{(k)}\left(x_{1}^{sc(k)}\right)=L_{1}H^{(k)}(x_{1}^{sc(k)}), & \quad k=1,....,n_{sc},\\
		A^{(k)}W_{L_{2}}^{(k)}\left(x_{1}^{sc(k)}\right)=L_{2}H^{(k)}(x_{1}^{sc(k)}), & \quad k=1,....,n_{sc},\\
	\end{cases}
	\end{equation*}
	where $H^{k}(x)$ in this case is given by:	
	\begin{align*}
		H^{(k)}(x)=\left(\begin{array}{c}
		\phi^{I_{1}}(||x-x^{sc(k)}||)\\
		\phi^{I_{2}}(||x-x^{sc(k)}||)\\
		\phi^{B_{1}}(||x-x^{fc(k)}||)\\
		\phi^{B_{2}}(||x-x^{fc(k)}||)\\
		\phi^{L_{1}}(||x-x^{pdec(k)}||)\\
		\phi^{L_{2}}(||x-x^{pdec(k)}||)
	\end{array}\right).
	\end{align*}
	
	Once the weights are known, we can build the sparse global matrix from the following system
	\begin{equation*}
	\begin{cases}
		W_{L_{1}}^{(k)}\left(x_{1}^{sc(k)}\right)d^{(,k)}=F_{1}(x_{1}^{sc(k)}) & \quad k=1,....,n_{sc},\\
		W_{L_{2}}^{(k)}\left(x_{1}^{sc(k)}\right)d^{(k)}=F_{2}(x_{1}^{sc(k)}) & \quad k=1,....,n_{sc}, \\
	\end{cases}
	\end{equation*}
	where 
	\begin{equation*}
		d^{(k)}=\Bigl(y_{1}(x^{sc(k)}),\, y_{2}(x^{sc(k)}),\, g_{1}(x^{fc(k)}),\, g_{2}(x^{fc(k)}),\, 
		{F}_{1}(x^{pdec(k)}),\, {F}_{2}(x^{pdec(k)})\Bigr)^T.
	\end{equation*}
	It is important to highlight that in order to avoid singularity of the sparse system we need that 
	$x_{1}^{sc(k)}\not\in\Omega_{pdec}^{k}$, see framework by \cite{stevens2011local}.	
	\begin{Obs}
		Analogous to Remark \ref{remark.timescheme.glogalscheme}, by considering Dirichlet boundary conditions,  
		the vector $\overline{\alpha}$ belongs to $\mathbb{R}^{d(n_{sc}^{k}+n_{pdec}^{k})}$, since 
		$B_{ij}$ is replaced by $B_{ij}=\delta_{\boldsymbol{x_{j}}} e_i$, and in consequence the 
		local Gram matrix $A^{(k)}$ is akin to the previous one. In the next numerical example we test such boundary
		conditions and postpone numerical experiments concerning the Navier--slip conditions until subsection 
		\ref{subsection.LHI.examples.unsteady}.  
	\end{Obs}

	
	\subsection{\normalsize{Numerical result: stationary case}}
	Here, the numerical example on an unitary square described in \cite{wendland2009divergence} is used to 
	examine the convergence of our method with nonhomogeneous Dirichlet boundary conditions. We have used
	advanpix package to super pass the ill condition Gram matrix	.  Thus, for 
	$\Omega=(0,1)\times(0,1)$ and $\mu=1$, we consider the following analytical solution to \eqref{eq:Stokes.system}: 
	$\mathbf{u}(\boldsymbol{x})=(20x_{1}x_{2}^{3},5x_{1}^{4}-5x_{2}^{4})$,\, 
	$p(\boldsymbol{x})=60x_{1}^2x_{2}-20x_{2}^{3}$. 
	Non--dimensional shape parameters  $c^*$ of values $1.0$ and $0.01$ are used throughout, as well as 25, 50, 80 local nodes 
	for inverse MQ and Gaussian RBFs.  Table \ref{table.LHI.steady} contains the approximation orders of the $L^2$ and 
	$L^\infty$ error of the velocity field, i.e., 
	$\|\boldsymbol{\epsilon_u}\|_{\infty}:=\|\mathbf{u}_{exact}-\mathbf{u}_{aprox}\|_{\infty}$ and
	$\|\boldsymbol{\epsilon_u}\|_{2}:=\|\mathbf{u}_{exact}-\mathbf{u}_{aprox}\|_{2}$.  
	
	\begin{center}
	\begin{table}[H]
		\centering{}%
				\begin{tabular}{cccccc}
		\hline 
		Total nodes  & Local nodes  & $c^*$ & $cond(A^{k})$  & $\|\boldsymbol{\epsilon_u}\|_{\infty}$  
		& $\|\boldsymbol{\epsilon_u}\|_{2}$ \tabularnewline
		\hline 
		400  & 50  & 1.0  & 3.77e+37  &\hspace{0.5cm} 4.78e-08  &\hspace{0.5cm} 2.74e-07 \tabularnewline
		1225  & 25  & 1.0  & 2.94e+32  &\hspace{0.5cm} 5.49e-06  &\hspace{0.5cm} 1.10e-04 \tabularnewline
		1225  & 80  & 1.0  & 2.99e+52  &\hspace{0.5cm} 9.34e-15  &\hspace{0.5cm} 8.31e-14 
		\vspace{0.5cm}
		\tabularnewline
		400  & 50  & 1.0  & 1.13e+26  &\hspace{0.5cm} 0.000485214  &\hspace{0.5cm} 0.00291627 \tabularnewline
		1225  & 25  & 1.0  & 5.15e+25  &\hspace{0.5cm} 1.09e-02  &\hspace{0.5cm} 1.94e-01 \tabularnewline
		1225  & 25  & 0.01  & 2.46e+40  &\hspace{0.5cm} 9.14e-09  &\hspace{0.5cm} 1.63e-07 \tabularnewline
		1225  & 80  & 1.0  & 4.65e+38  &\hspace{0.5cm} 1.24e-07  &\hspace{0.5cm} 1.41e-06 \tabularnewline
		\hline 
		\end{tabular}

		\caption{Approximation errors of the velocity field for the steady Stokes problem using LHI method.  
		Two shape parameters $c^*$ are considered; top: Gaussian RBF; bottom: inverse MQ RBF.}
		\label{table.LHI.steady}
	\end{table}
	\par\end{center}

	\subsection{\normalsize{Evolutionary problem}}
	In this subsection we formulate a RBF--LHI vectorial technique for the evolutionary
	Stokes problem based in Section \ref{subsection.operator.discre}, which create a modified PDE operator using
	the finite difference method (FDM) for the time discretization, as well as BFD methods. Specifically, for the system
	
	\begin{equation}\label{system.lhi-unsteadycase}
	\left\{ \begin{array}{lll}
		\mathbf{y}_{t}+L(\mathbf{y},p)=\mathbf{f} & \text{ in } & Q,\\
		\nabla\cdot\mathbf{y}=0 & \text{ in } & Q,\\
		B\mathbf{y}=\mathbf{g} & \text{ on } & \Sigma,\\
		\mathbf{y}(\cdot,0)=\mathbf{y}_{0}(\cdot) & \text{ in } & Q,
	\end{array}\right.
	\end{equation} 
	where $L(\mathbf{y},p)=-\mu\Delta \mathbf{y}+\nabla p$, we define the operator of left--hand side $\bar{L}$ by 
	\begin{equation*}
		\bar{L}(\mathbf{y}^{n+s},p^{n+s}):=\mathbf{y}^{n+s}+\Delta t\,\beta_{s}\,L(\mathbf{y}^{n+s},p^{n+s})
	\end{equation*} 
	and the source term 	
	\begin{equation*}
		\boldsymbol{F}^{n+s}:=\Delta t\,\beta_{s}\,\mathbf{f}^{n+s}+\sum\limits_{k=0}^{s-1}\sigma_{k}\mathbf{y}^{n+k}.
	\end{equation*} 
	\noindent
	Following the same procedure as in the stationary case, our ansatz  is given by	

	{\small{}
	\begin{align}\label{eq:local anzat bdf2}
		(\hat{y}^{(n+2,k)}(\boldsymbol{x}),\hat{p}^{(n+2,k)}(\boldsymbol{x})) & =
		\sum_{j=1}^{n_{sc}^{k}}\phi^{I_{1}}(\boldsymbol{x}-\boldsymbol{x_{j}^{sc(k)}})\alpha_{1j}^{I}
		+
		\sum_{j=1}^{n_{sc}^{k}}\phi^{I_{2}}(\boldsymbol{x}-\boldsymbol{x_{j}^{sc(k)}})\alpha_{2j}^{I} \\
		 & +
		 \sum_{j=1}^{n_{fc}^{k}}\phi^{B_{1}}(\boldsymbol{x}-\boldsymbol{x_{j}^{fc(k)}})\alpha_{1j}^{B}
		 +
		 \sum_{j=1}^{n_{fc}^{k}}\phi^{B_{2}}(\boldsymbol{x}-\boldsymbol{x_{j}^{fc(k)}})\alpha_{2j}^{B}\nonumber\\
 		& +
 		\sum_{j=1}^{n_{pdec}^{k}}\phi^{\bar{L}_{1}}(\boldsymbol{x}-\boldsymbol{x_{j}^{pdec(k)}})\alpha_{1j}^{L}
		 +
		 \sum_{j=1}^{n_{pdec}^{k}}\phi^{\bar{L}_{2}}(\boldsymbol{x}-\boldsymbol{x_{j}^{pdec(k)}})\alpha_{2j}^{L},\nonumber 
		\end{align}
	}{\small \par}
	which is collocated in the local Stokes system
	\begin{equation}\label{eq:local system}
		\left\{ \begin{array}{lll}
		\mathbf{u_n}(x_{i})=\overline{h_{i}} & \text{ in } & \Omega_{sc}^{k} \subset \Omega_{sc}\bigcap D^{k},\\
		\\
		{\bar {L}}\mathbf{u}=\mathbf{F} & \text{ in } & \Omega_{pdec}^{k} \subset \Omega_{pdec}\bigcap D^{k},
		\\
		\\
		B\mathbf{u}=g & \text{ on } & \Omega_{fc}^{k} \subset \Omega_{fc}\bigcap D^{k}.
		\end{array}\right.
	\end{equation}
	After that, we can obtain the local Gram matrix, which	in turn let us to compute the weights. 
	Once the weights are known, we can build a sparse global system using the system  
	
	\begin{equation}\label{lhi:sparse.stokes.evol}
	\begin{cases}
		W_{L_{1}}^{(k)}\left(x_{1}^{sc(k)}\right)d^{(n+s,k)}=\Delta t\,\beta_{s}\,f^{n+s}_{1}(x_{1}^{sc(k)})
		+\sum\limits_{k=0}^{s-1}\sigma_{k}y^{n+k}_{1}(x_{1}^{sc(k)}), & \quad k=1,\dots,n_{sc},\\
		W_{L_{2}}^{(k)}\left(x_{1}^{sc(k)}\right)d^{(n+s,k)}=\Delta t\,\beta_{s}\,f^{n+s}_{2}(x_{1}^{sc(k)})+\sum\limits_{k=0}
		^{s-1}\sigma_{k}y^{n+k}_{2}(x_{1}^{sc(k)}),  & \quad k=1,\dots,n_{sc},
	\end{cases}
	\end{equation}
	for each $x_{1}^{sc(k)}$ in  $D^{k}$, and where $d^{(n+s,k)}$ is
	\begin{equation*}
		(y_{1}^{(n+s)}(x^{sc(k)}),\, y_{2}^{(n+s)}(x^{sc(k)}),\, g_{1}^{n+s}(x^{fc(k)}),\, g_{2}^{n+s}(x^{fc(k)}),\,
		F^{n+s}_{1}(x^{pdec(k)}),\, F^{n+s}_{2}(x^{pdec(k)}))^T.
	\end{equation*}
	Observe that (\ref{lhi:sparse.stokes.evol}) can be expressed in matrix form as:

	\begin{equation}\label{lhi:matrix.form.sparse}
	\left(\begin{array}{cccccccc}
		W_{L_{1}}^{y_{1}} & W_{L_{1}}^{y_{2}}  & W_{L_{1}}^{B_{1}} & W_{L_{1}}^{B_{2}}  
		& W_{L_{1}}^{L_{1}} & W_{L_{1}}^{L_{2}}\\
		W_{L_{2}}^{y_{1}} & W_{L_{2}}^{y_{2}}  & W_{L_{2}}^{B_{1}} & W_{L_{2}}^{B_{2}} 
		& W_{L_{2}}^{L_{1}} & W_{L_{2}}^{L_{2}}
	\end{array}\right)\left[\begin{array}{c}
		y_{1}^{n+s}\\
		y_{2}^{n+s}\\
		g_{1}^{n+s}\\
		g_{2}^{n+s}\\
		F_{1}^{n+s}\\
		F_{2}^{n+s}
	\end{array}\right]=\left[\begin{array}{c}
		F_{1}^{n+s}\\
		F_{2}^{n+s}
	\end{array}\right],
	\end{equation}
	where the functions $F_{1}^{n+s}, F_{2}^{n+s}$ in the right--hand side of \eqref{lhi:matrix.form.sparse} are 
	evaluated 
	in every $x_{1}^{sc(k)}$ of $D^{k}$, meanwhile $F_{1}^{n+s}$ and  $F_{2}^{n+s}$ in the left--hand side of
	\eqref{lhi:matrix.form.sparse} in $x^{pdec(k)}$. Thus, it allows to compute $(y^{n+s},p^{n+s})$. 
	
	\begin{Obs}
	In a more general setting, we can establish an analysis of stability for the LHI method presented above.
	By simplicity, let us consider $x_{1}^{sc(k)}=x^{pdec(k)}$ 
	on each disk $D^k$. Besides, let $\mathbf{e}:=\mathbf{y}-\hat{\mathbf{y}}$ be the error between the exact and 
	approximate solution in which we have discarded the local truncation error. Thus,  
	(\ref{lhi:matrix.form.sparse}) can be transformed for the error $\mathbf{e}_n$ in the form
	
	\begin{equation*}
	\left(\begin{array}{cc}
		W_{L_{1}}^{y_{1}} & W_{L_{1}}^{y_{2}}\\
		W_{L_{2}}^{y_{1}} & W_{L_{2}}^{y_{2}}
	\end{array}\right)\left[\begin{array}{c}
		e_{1}^{n+s}\\
		e_{2}^{n+s}
	\end{array}\right]=\left[\begin{array}{c}
		\sum\limits_{k=0}^{s-1}\sigma_{k}e_{1}^{n+k}\\
		\sum\limits_{k=0}^{s-1}\sigma_{k}e_{2}^{n+k}
	\end{array}\right]-\left(\begin{array}{cc}
		W_{L_{1}}^{L_{1}} & W_{L_{1}}^{L_{2}}\\
		W_{L_{2}}^{L_{1}} & W_{L_{2}}^{L_{2}}
	\end{array}\right)\left[\begin{array}{c}
		\sum\limits_{k=0}^{s-1}\sigma_{k}e_{1}^{n+k}\\
		\sum\limits_{k=0}^{s-1}\sigma_{k}e_{2}^{n+k}
	\end{array}\right].
	\end{equation*}
	Thus,  $\mathbf{e}^{n+s}=S_{1}^{-1}(I-S_{2})(\sum\limits_{k=0}^{s-1}\sigma_{k}e^{n+k})$, where 
	\[
		S_{1}=\left(\begin{array}{cc}
		W_{L_{1}}^{y_{1}} & W_{L_{1}}^{y_{2}}\\
		W_{L_{2}}^{y_{1}} & W_{L_{2}}^{y_{2}}
		\end{array}\right)\quad\mbox{and}\quad S_{2}=\left(\begin{array}{cc}
		W_{L_{1}}^{L_{1}} & W_{L_{1}}^{L_{2}}\\
		W_{L_{2}}^{L_{1}} & W_{L_{2}}^{L_{2}}
		\end{array}\right).
	\]
	
	Therefore, the LHI method will be stable as long that eigenvalues of $S_{1}^{-1}(I-S_{2})$ lie in the stability 
	region of \eqref{eq:region.stability}.
	\end{Obs}


\subsection{\normalsize{Numerical results: evolutionary case}}\label{subsection.LHI.examples.unsteady}
	From the theoretical description previously developed for the RBF--LHI method for the unsteady 
	Stokes system, the present  subsection shows numerical results of its implementation.  We build the divergence 
	free kernel by using inverse multi quadric (MQ) with a shape parameter $c=0.1$ and again, the Matlab package ADVANPIX  is considered. In this case, 
	we introduce a uniform mesh generate by means of the package DISTMESH with 1312 total nodes. The viscosity coefficient is $\mu=1$. In order to test 
	our numerical approach, we consider the following analytical solution to \eqref{system.lhi-unsteadycase} over the unitary circle. 
	$$(\mathbf{y_1}(x,y,t),\mathbf{y_2}(x,y,t))=\Bigl(-\pi\,y\,\sin\left(\frac{\pi}{2}(x^{2}+y^{2})\right)\sin(\pi t), 
	\pi\,x\,\sin\left(\frac{\pi}{2}(x^{2}+y^{2})\right)\sin(\pi t)\Bigr)$$ 
	and pressure $p(x,y,t)=\sin(x-y+t)$. 
	
	 \noindent
	 As in section \ref{sec:global.numerical.results.line}, we compare the velocity error in the $L^{\infty}$--norm between the
	 exact and numerical solutions, i.e., $\boldsymbol{\epsilon_y}=\mathbf{y}_{exact}-\mathbf{y}_{aprox}$. Here, we decide to omit the pressure error since an extra
	computation is needed to calculate. The results are presented in Table \ref{table.difusion_2}.
	 \vskip 0.4cm
	
	In addition to Table \ref{table.difusion_2} where we note non--convergence cases, Figure \ref{fig.eigen_lhi} gives an answer about it. Figure \ref{fig.eigen_lhi} shows
	the eigenvalues of the matrix $S_{1}^{-1}(I-S_{2})$ described in the previous section and its distribution over the stability region (SR).  Indeed, we can appreciate that 
	for the non-convergence cases the eigenvalues  are located outside the SR of the BDF2 method (\eqref{eq:region.stability}), and therefore the stability condition is broken. 
	That means again that there exists a hidden condition depending on the shape parameter, the diffusion coefficient, time--step, fill-distance and local-fill distance. 
	
	\begin{table}[H]
	\begin{center}
		\begin{tabular}{ccccc}
	\hline 
	B.C  & \multicolumn{2}{c}{Dirichlet} & \multicolumn{2}{c}{Navier-slip}\tabularnewline
	\hline 
	Stencil Size & 30 & 60 & 30 & 60\tabularnewline
	\hline 
	$\Delta t$  & \multicolumn{2}{c}{$\left\Vert \boldsymbol{\epsilon_{y}}\right\Vert_\infty $ } & \multicolumn{2}{c}{$\left\Vert \boldsymbol{\epsilon_{y}}\right\Vert_\infty $ }
	\tabularnewline
	\hline 
	0.02 &\hspace{0.7cm} 3.72E-04 &\hspace{0.7cm} 3.72E-04 &\hspace{0.7cm} 3.07E-03 &\hspace{0.7cm}  3.08E-03\tabularnewline
	0.01 &\hspace{0.7cm} 9.36E-05 &\hspace{0.7cm} 9.36E-05 &\hspace{0.7cm} 7.85E-04 &\hspace{0.7cm} 7.91E-04\tabularnewline
	0.005 &\hspace{0.7cm} 2.35E-05 &\hspace{0.7cm} N.C &\hspace{0.7cm} 1.96E-04 &\hspace{0.7cm} 2.01E-04\tabularnewline
	0.002 &\hspace{0.7cm} N.C &\hspace{0.7cm} N.C &\hspace{0.7cm} N.C &\hspace{0.7cm} 3.38E-05\tabularnewline
	\hline 
	\end{tabular}

	\par
	\caption{Convergence behaviour of velocity in the $L^{\infty}$\textendash norm
	using inverse MQ.  Non--convergence (N.C).}
	\label{table.difusion_2} 
	\end{center}
	\end{table}

	\begin{figure}[H]\label{fig.eigen_lhi}
	\begin{center}
	\subfloat[]{\includegraphics[scale=0.1]{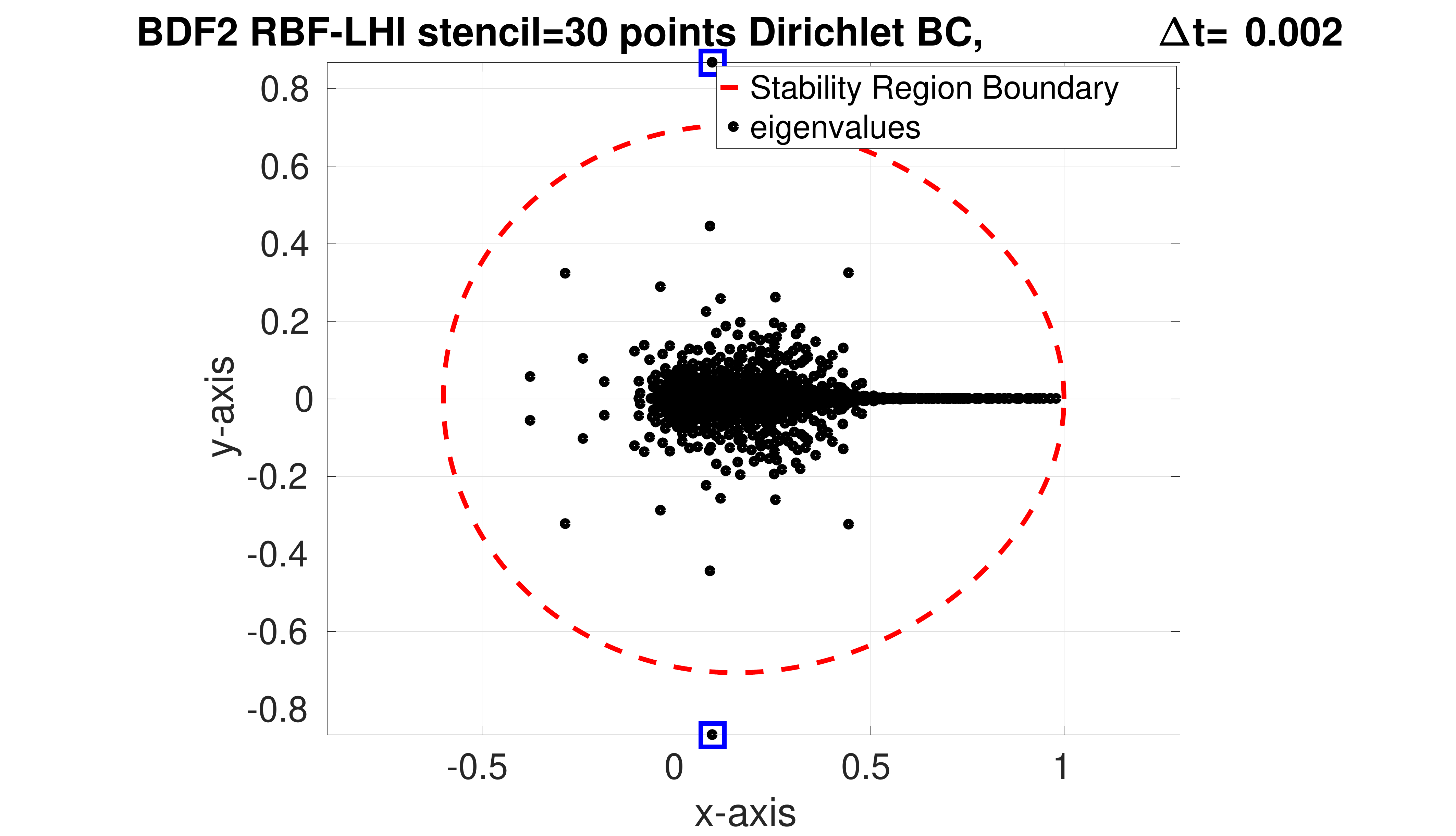}}
	\hspace{0.01cm} 
	\subfloat[]{\includegraphics[scale=0.1]{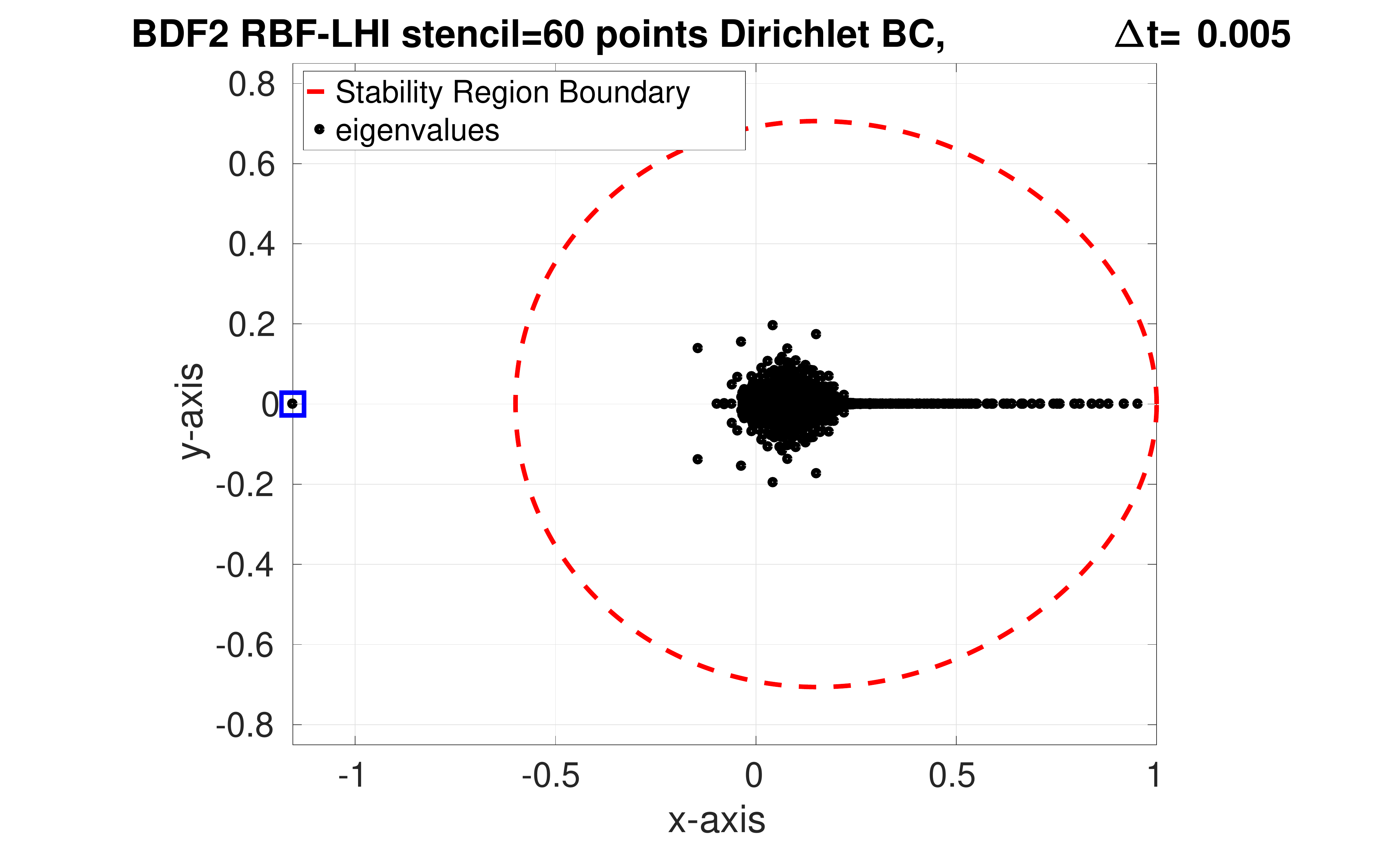}}
	\hspace{0.01cm} 
	\subfloat[]{\includegraphics[scale=0.1]{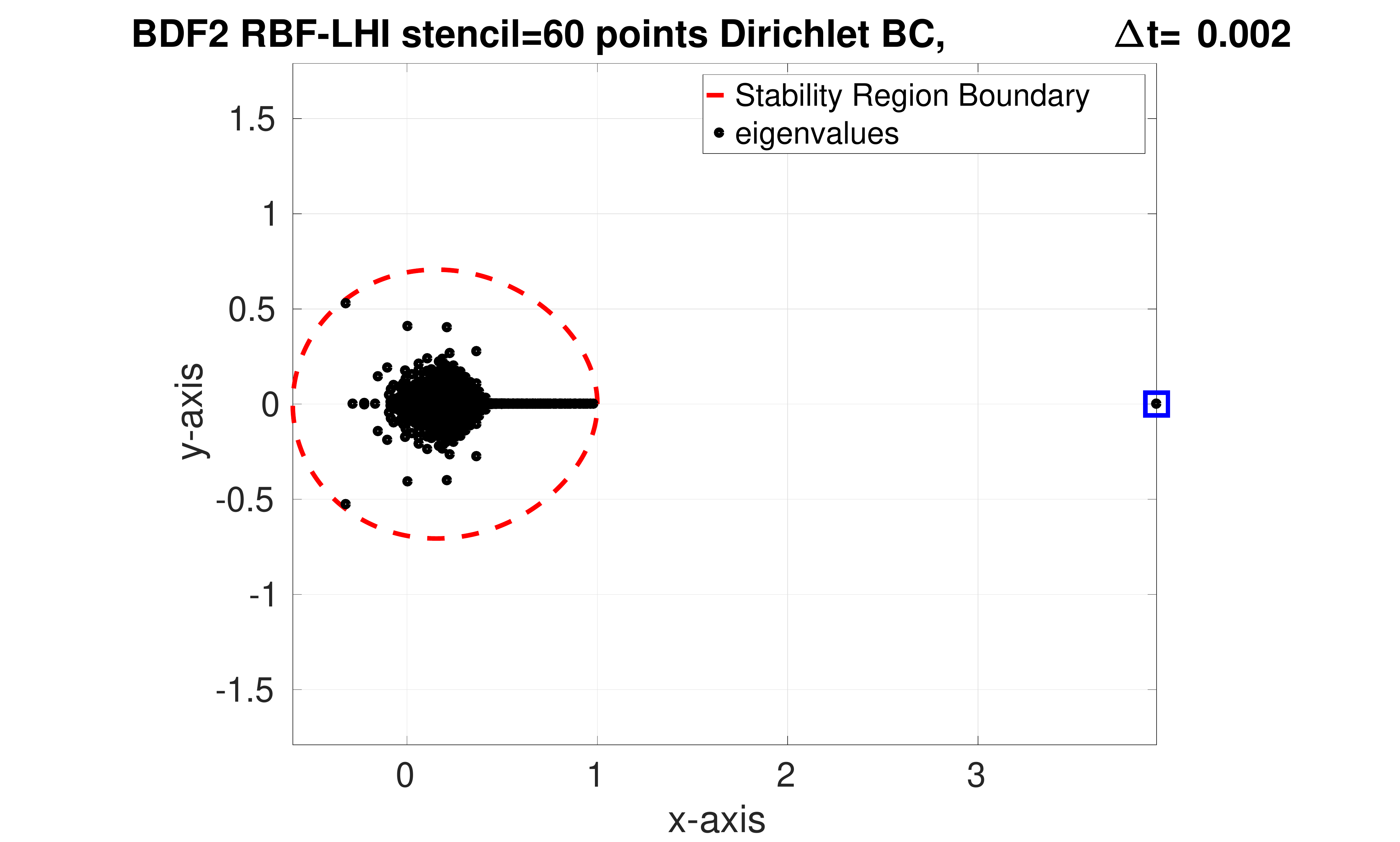}}
	\subfloat[]{\includegraphics[scale=0.1]{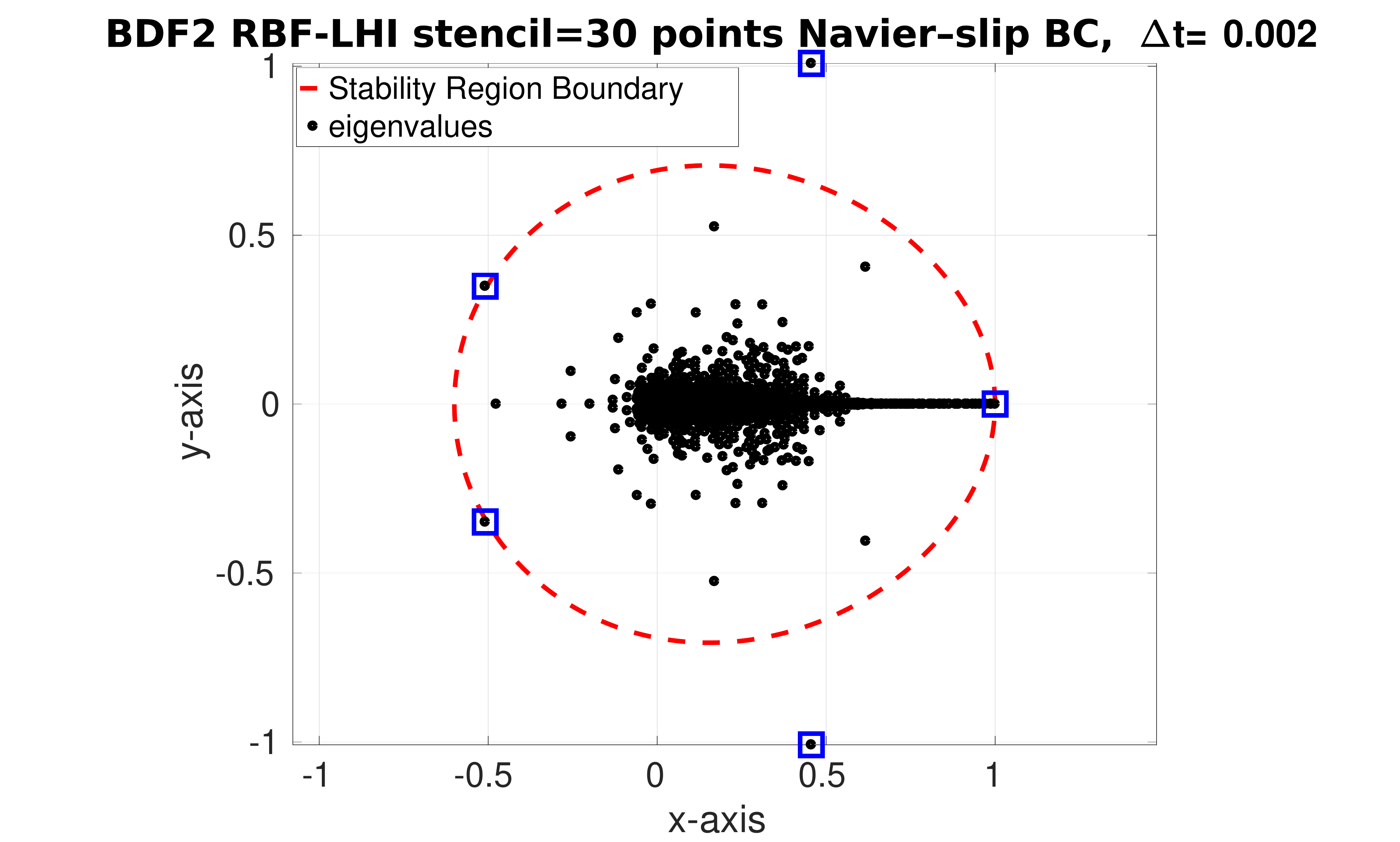}}
	\end{center}
	\par
	\caption{Representation of eigenvalues of the matrix $S_{1}^{-1}(I-S_{2})$.  The blue squares represent the eigenvalues located outside the stability region. 
	The a) b) c) figures are associated to the Dirichlet boundary conditions and d) corresponds to the case of Navier--slip boundary conditions.}
	\label{fig.eigen_lhi}
	\end{figure}


\section{\normalsize{Numerical control problem}}\label{section_control}
	The approximate solution to the null controllability problem for the two-dimensional 
	Stokes system with few scalar controls is carried out in this section. Hence, both
	Dirichlet and Navier--slip boundary conditions are considered.
	As mention before, the numerical implementation follows the method developed above through RBF--LHI, and
	simultaneously, those are compared against a mixed formulation using finite element method (FEM) for only the Dirichlet case, since the implementation of Navier-Slip is not straightforward with FEM. Although this is perfectly possible, the major problem is that both essential and natural boundary conditions should be incorporated within the algorithm to make it reliable. The shape of the domine, in our case a circle, makes this task more difficult.
	\vskip 0.2cm
\noindent
	
	Following \cite{glowinski2008exact},  the CGM is implemented with a stopping criteria of $\epsilon = 10^{-8}$ for solving the dual system \eqref{intro.Stokes.system}, 	
	\eqref{Stokes.adjointsystem1}. We use the  following data: 
	$\Omega=\{(x,y)\in\mathbb{R}^2: x^2+y^2<1\}$, the observation set  ${\omega}{=\{(x,y)\in\mathbb{R}^2: x^2+y^2<0.25\}}$,\, ${T=0.25}$. The reason for this choice of $T$
	is that if we take $T=1$,  the  solution with control force is very close to the solution without control, thus, it is difficult to appreciate the effect of the control in the 
	numerical experiments.
\vskip 0.2cm
\noindent	
	In all cases, we use a uniform mesh of 3512 points generated with FreeFem++, the time step size is $\Delta t = \frac{1}{200}$ and  diffusion coefficient 
	$\mu=1$. For the initial condition we choose $(y^0_1,y^0_2)=(-\pi\,y\,\cos\left(\frac{\pi}{2}(x^{2}+y^{2})\right)^2,\pi\,x\,\cos\left(\frac{\pi}{2}(x^{2}+y^{2})\right)^2)$. 
	Regarding the functional (\ref{functional.Stokes1}),  we set the regularization parameters 
	$c_1=1/300$ , $c_2^{-1}=0$ by having  controls with both non-zeros scalar component ($\mathbf{v}=(v_1,v_2)$) and  $c_1^{-1}=0$, $c_2=1/300$ 
	by considering controls with one scalar control (either $\mathbf{v}=(v_1,0)$ or $\mathbf{v}=(0,v_2)$). 

\begin{Obs}
	For the numerical experiments we use a triangular mesh for two specific reasons, the first one is to have a fair comparison
	 between RBF-LHI and finite element, and the second is because CGM require to calculate integrals over the domain thus for LHI-RBF method it is more efficient 
	 to use the triangulation to calculate integrals with  $\mathcal{P}_1$--type elements. 
\end{Obs}


\subsection{\normalsize{Divergence free RBF-LHI}} 
	To generate the divergence free kernel for the LHI method we use IMQ-RBF with a constant shape parameter $c=0.1$. 
	It's important to underline that  the indicator function $\mathbbm{1}_{\omega}$  given in the system \eqref{intro.Stokes.system} is approximate by the smooth function 
	$1/(1+exp(-2k(0.5-||x||_{2})))$ with $k=20$,  otherwise our solution degenerate, which can be explained due to Gibs phenomena.

	\vskip 0.2cm
	\noindent
	Table \ref{table.iterracion.rbf} shows the number of iterations  to achieve the stopping criteria $\epsilon=10^{-8}$ in the CGM implemented. 
	
	\begin{table}[H]
	\begin{center}
	\begin{tabular}{cccc}
	\hline 
	B.C & $\mathbf{v}=(v_{1},v_{2})$ & $\mathbf{v}=(v_{1},0)$ & $\mathbf{v}=(0,v_{2})$ \tabularnewline
	\hline 
	Navier--slip & 36 & 22 & 25\tabularnewline
	Dirichlet & 22 & 22 & 18\tabularnewline
\end{tabular}
\par
	\caption{Number of iterations for obtaining the convergence criteria of the CGM.}
	\label{table.iterracion.rbf}
	\end{center}
	\end{table}
	
\begin{figure}[H]
	\begin{center}
		\includegraphics[scale=0.16]{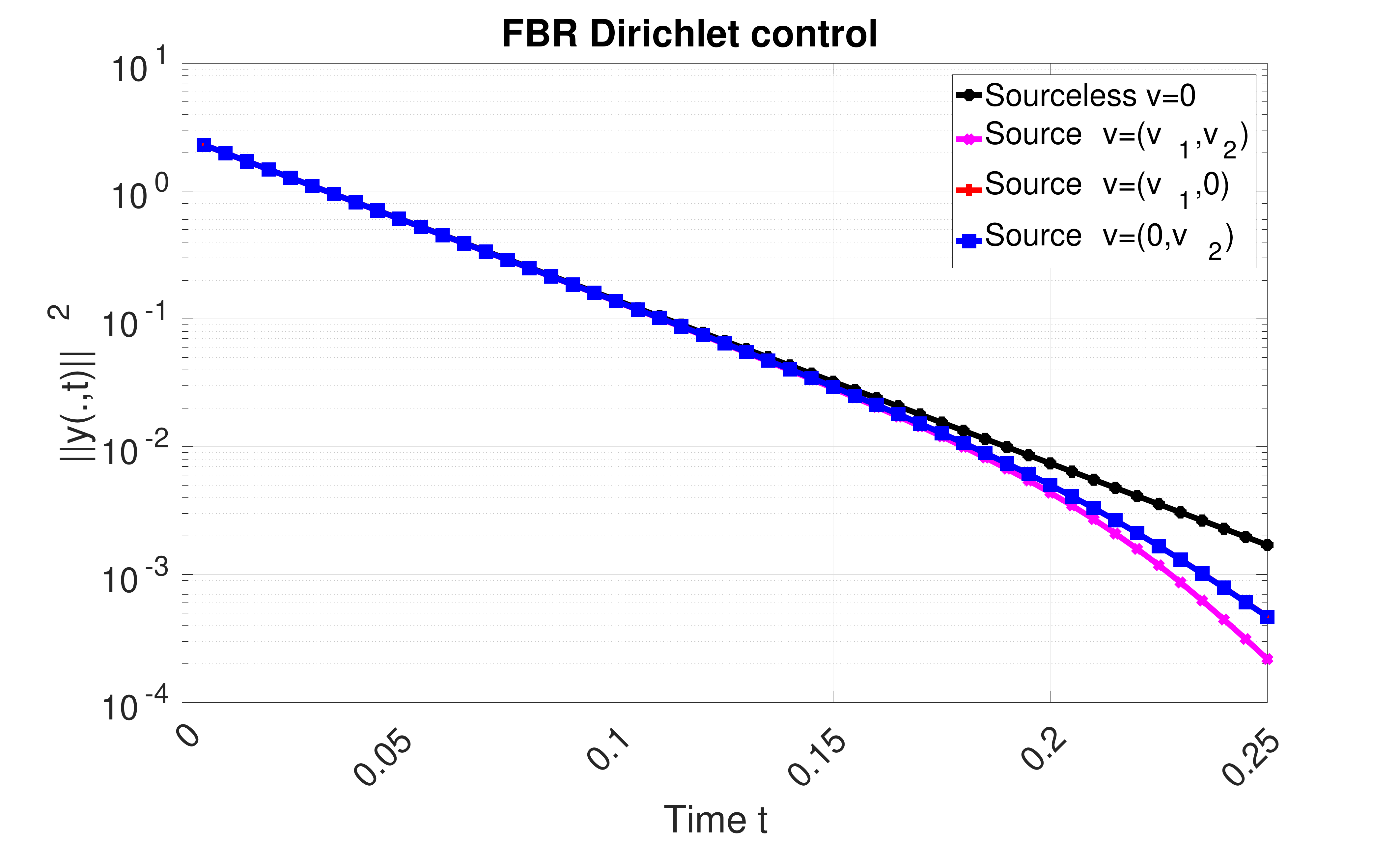}
		\hspace{0.01cm} 
		\includegraphics[scale=0.16]{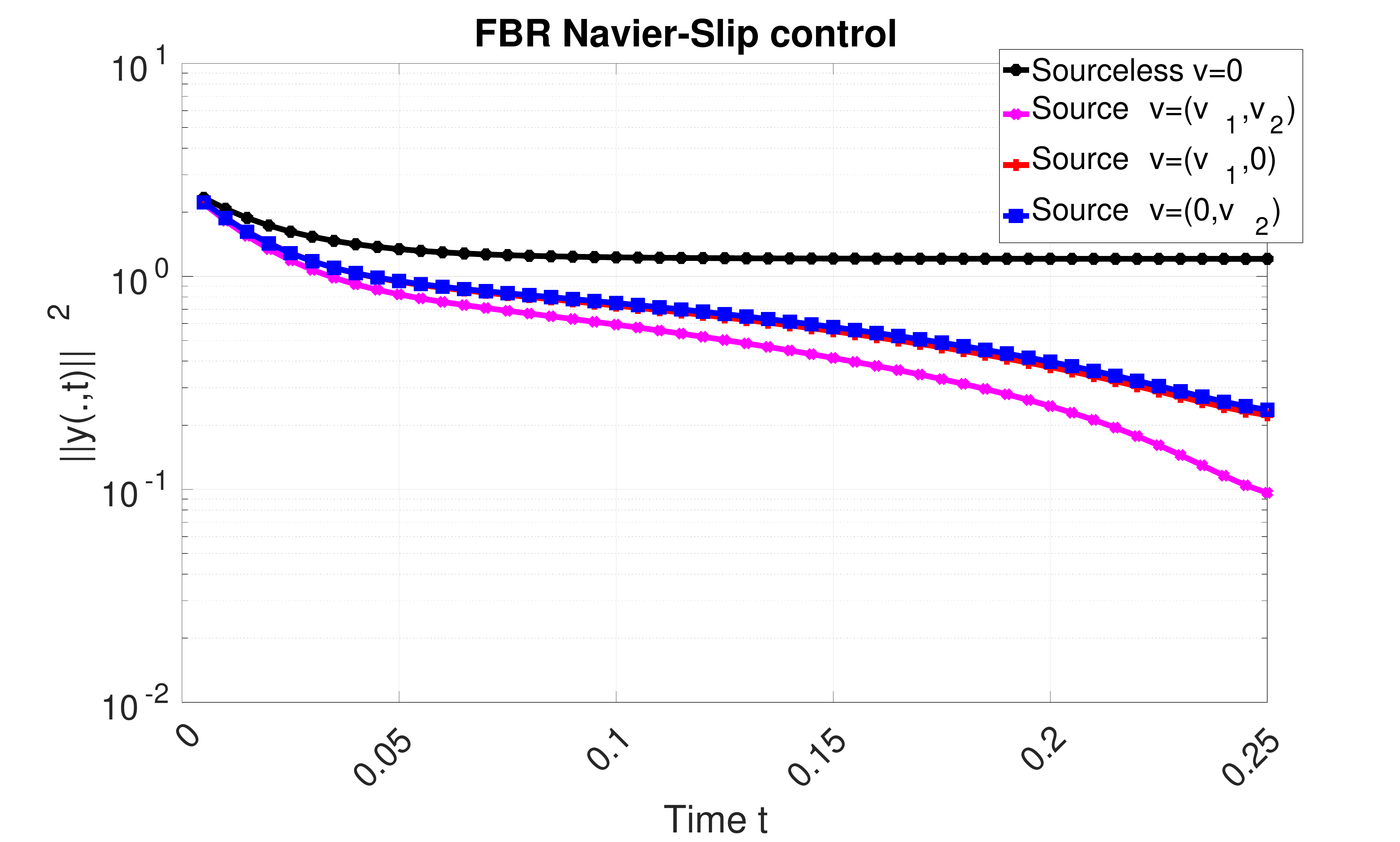}
	\end{center}\par
	\caption{$L^2$--norm square of the velocity field (as time function) of the solution of the null control problem with controls $\mathbf{v}=\mathbf{0}$ (black), 
	$\mathbf{v}=(v_1,v_2)$ (pink), 
	$\mathbf{v}=(v_1,0)$ (red) and $\mathbf{v}=(0,v_2)$ (blue). Dirichlet boundary condition (left) and Navier--slip boundary condition (right).}
	\label{fig.finalstate.fewcontrols}
	\end{figure}

	\noindent
	Table \ref{table.finalstate.fewcontrols} and Figure \ref{fig.finalstate.fewcontrols} show the $L^2$--norm of the velocity vector field for the null control
	problems as time function.  The numerical control function $\mathbf{v}$ has all possible structures, namely, $\mathbf{v}=\mathbf{0}$, $\mathbf{v}=(v_1,v_2)$, 
	$\mathbf{v}=(v_1,0)$ and $\mathbf{v}=(0,v_2)$. 
	
\begin{figure}[H]
		\begin{center}
		\includegraphics[scale=0.12]{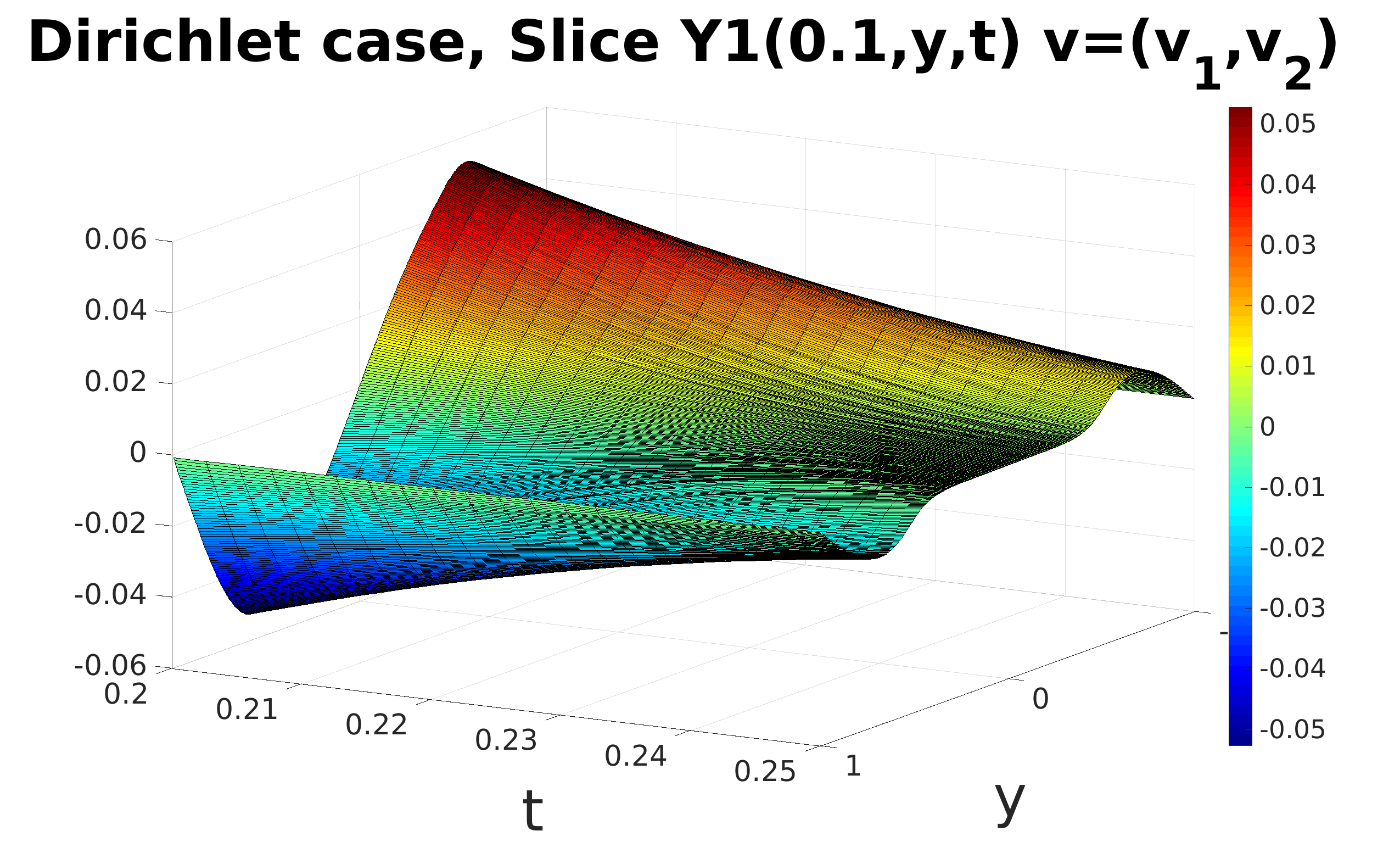} 
   		\hspace{0.01cm}
		\includegraphics[scale=0.12]{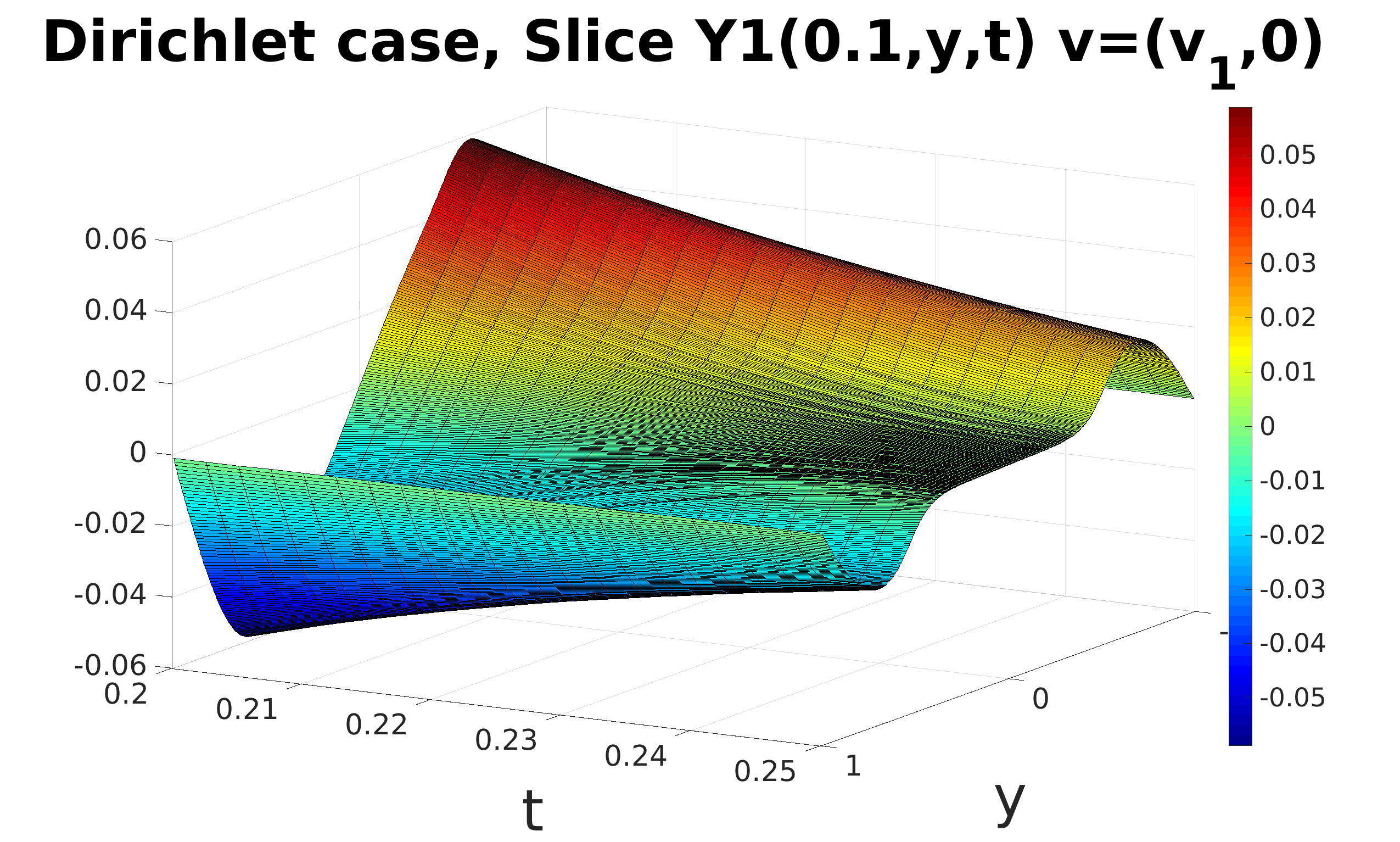} 
    	\hspace{0.01cm}
		\includegraphics[scale=0.12]{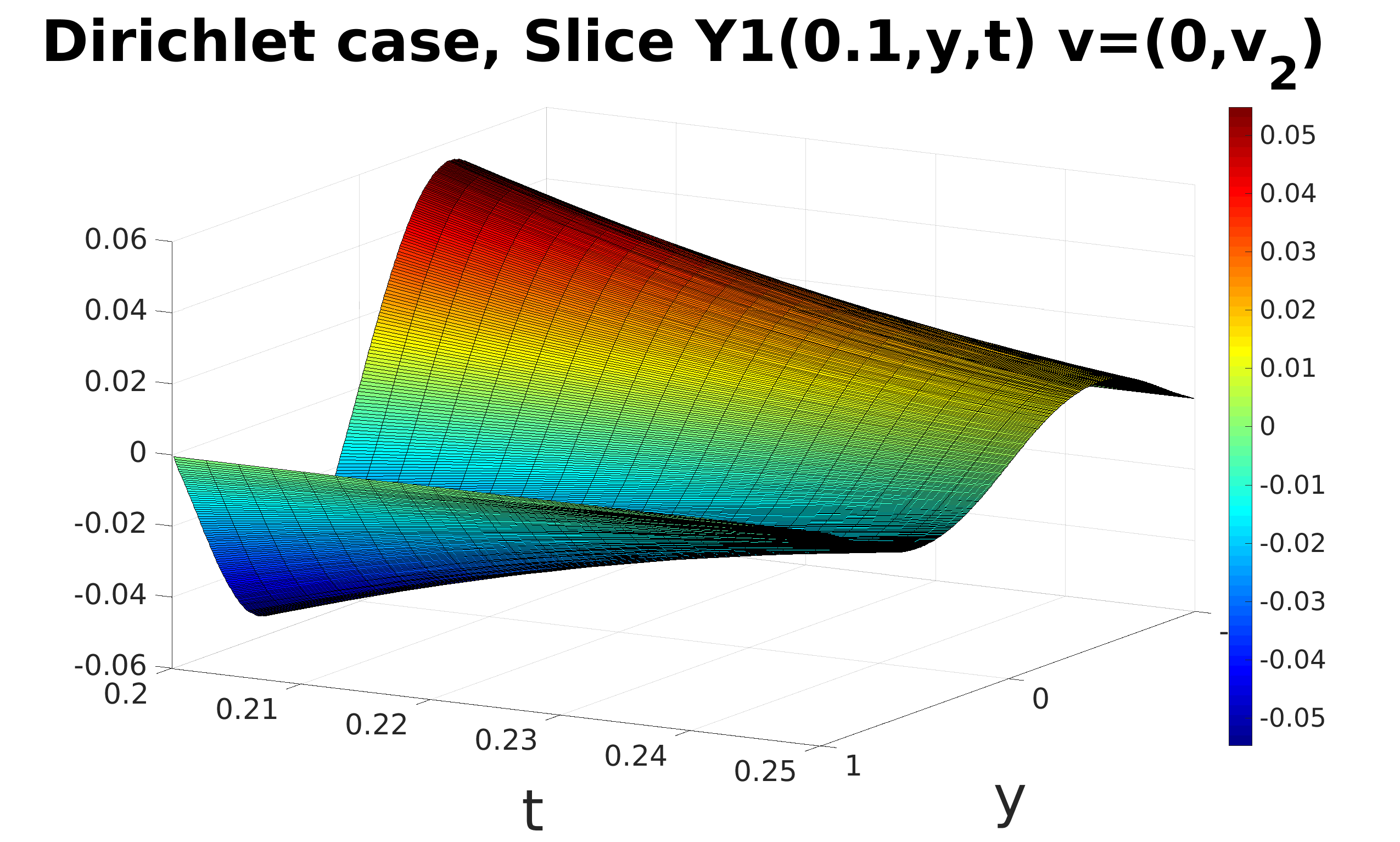} 
   		\vspace{0.01cm}  
		\includegraphics[scale=0.12]{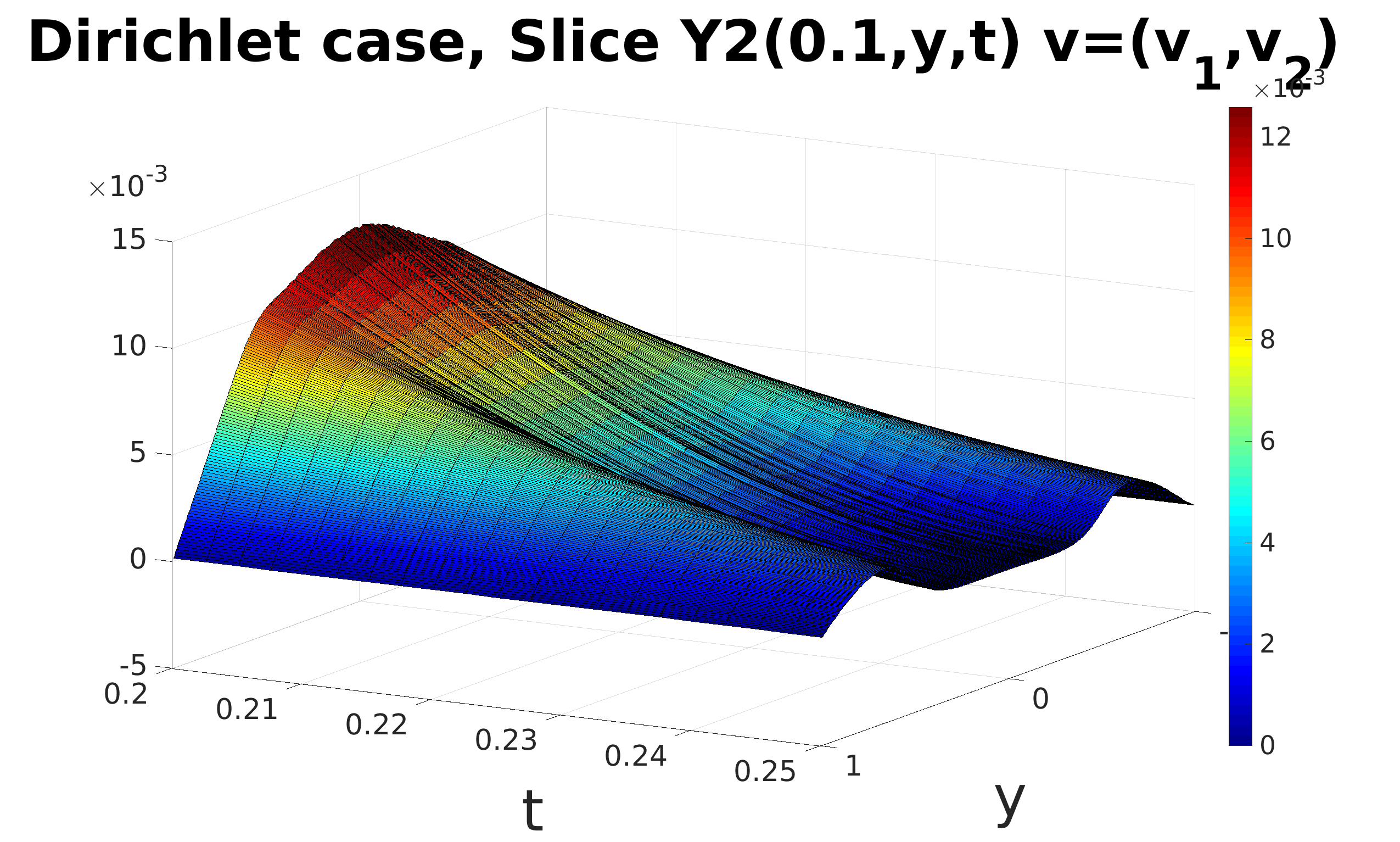} 
		\hspace{0.01cm}  
		\includegraphics[scale=0.12]{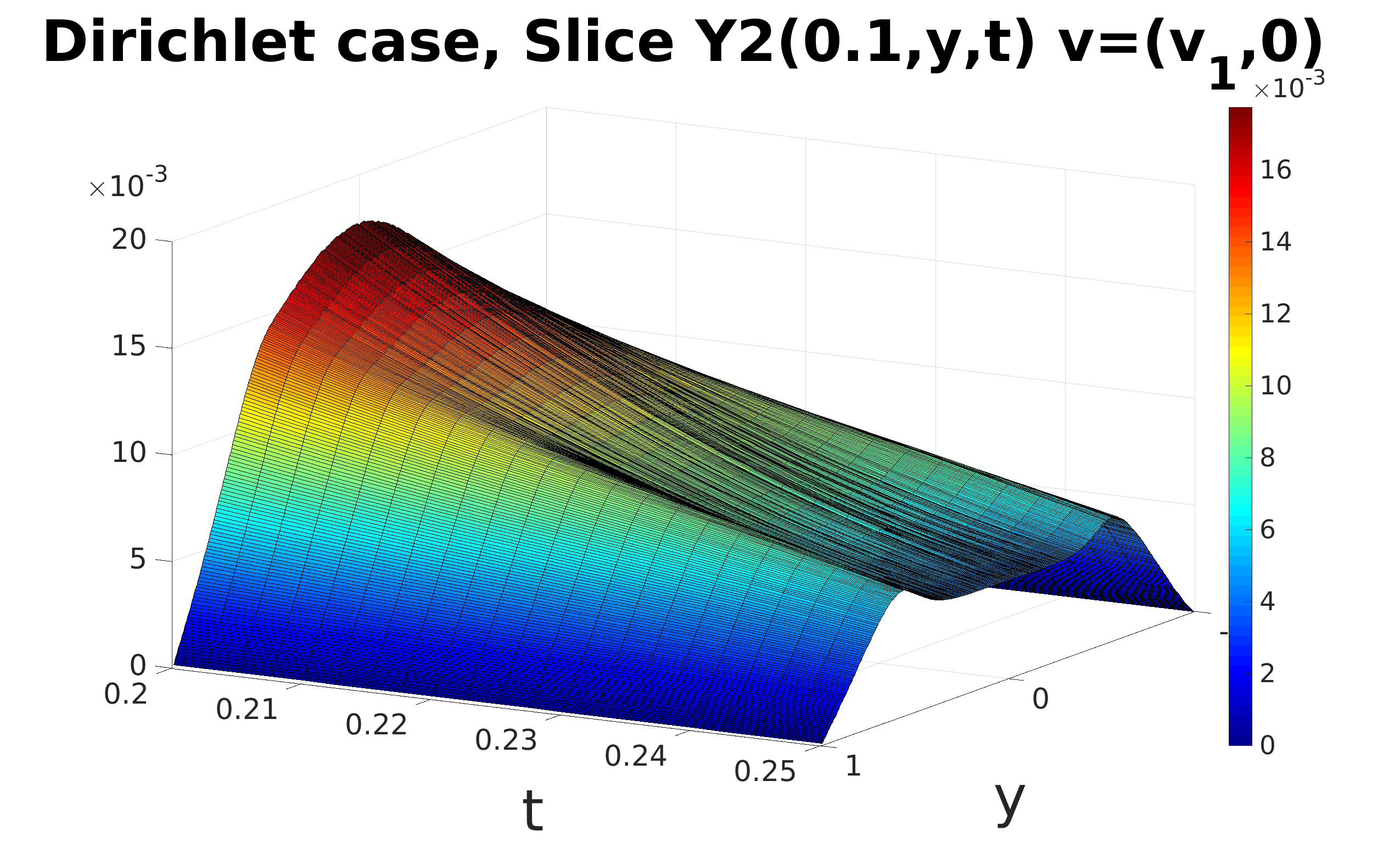} 
   		\hspace{0.01cm}  	
		\includegraphics[scale=0.12]{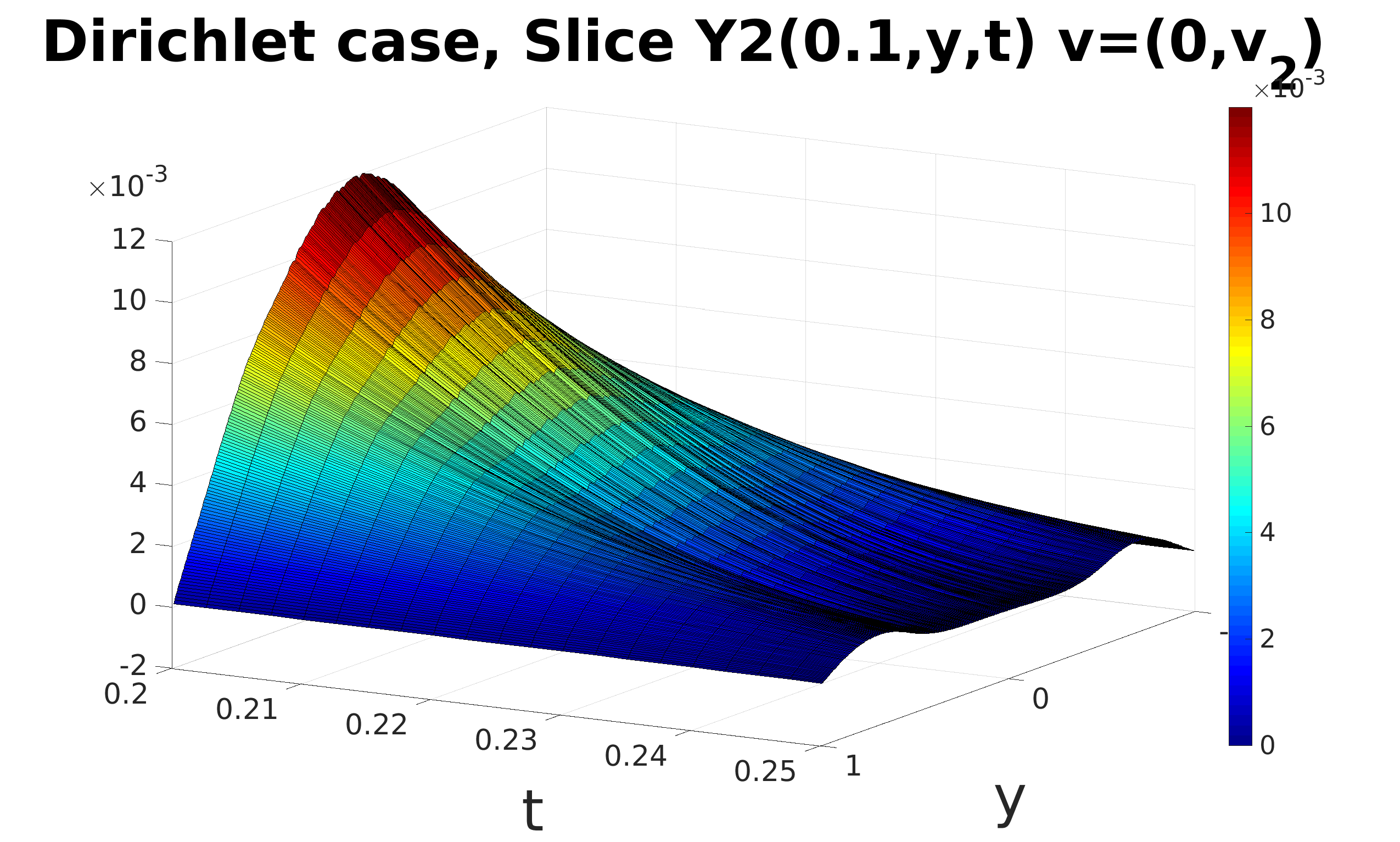} 
   		\hspace{0.01cm}  
		\end{center}\par
	\caption{Dirichlet boundary conditions. Cut sections of the components of the state $\mathbf{y}=(y_1,y_2)$ by the plane $x=0.1$. Cut section of  $y_1$ (first row) and 
	cut section of  $y_2$ (second row). A different control function is represented in each column.}
	\label{fig.cuts.states.RBF.dirich}
	\end{figure}

\begin{figure}[H]
		\begin{center}
		\includegraphics[scale=0.12]{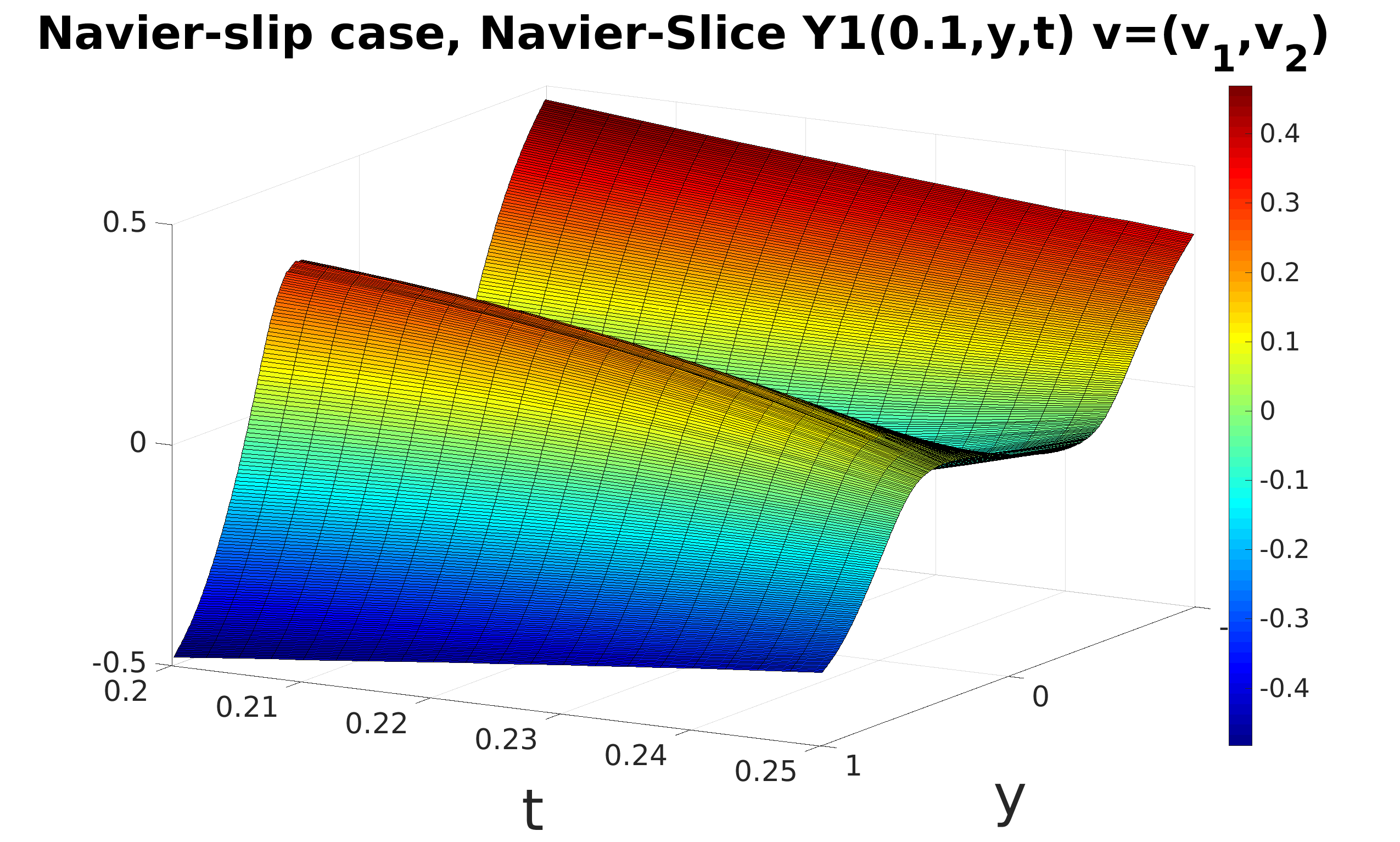} 
   		\hspace{0.01cm} 
		\includegraphics[scale=0.12]{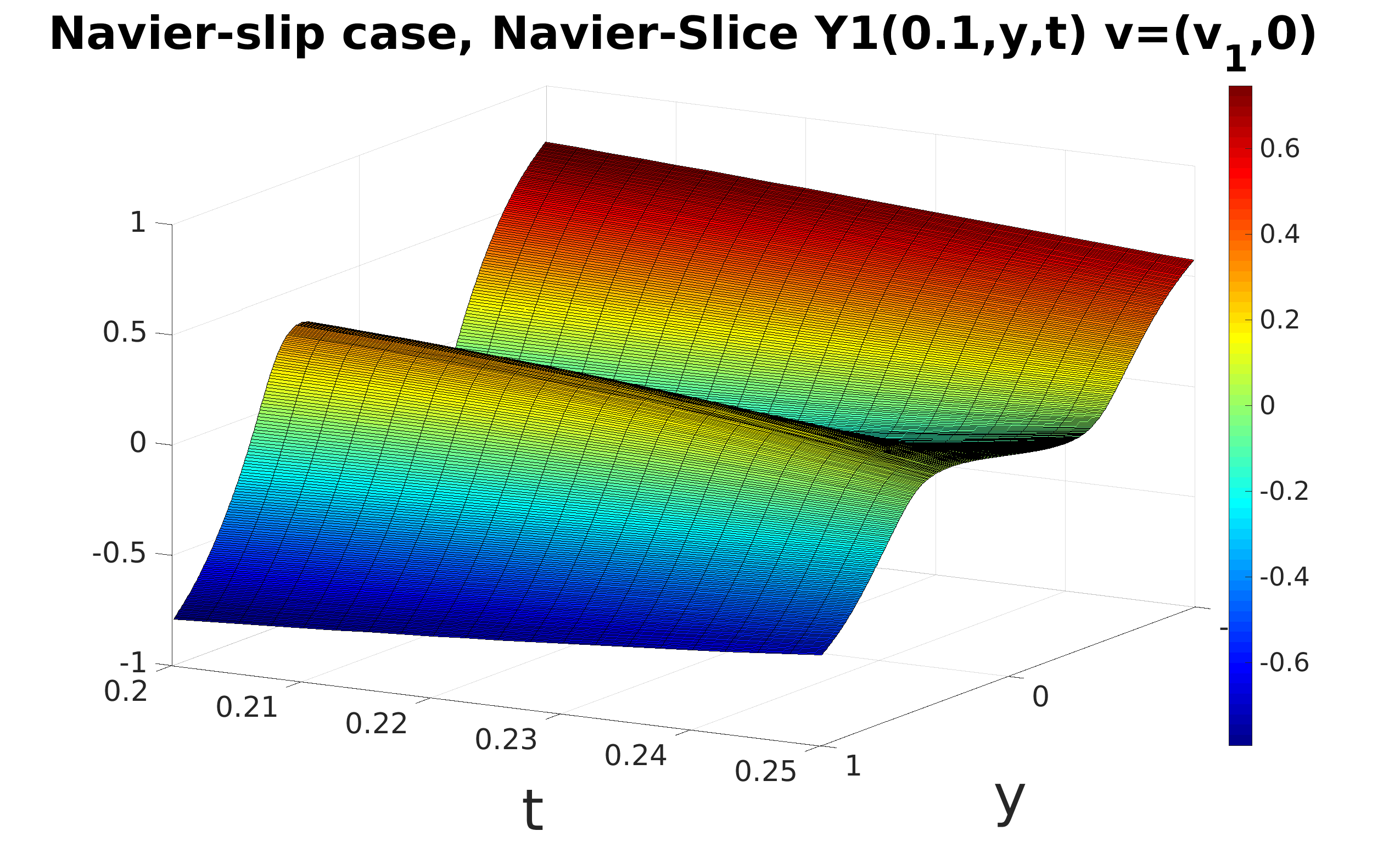} 
   		\hspace{0.01cm} 
		\includegraphics[scale=0.12]{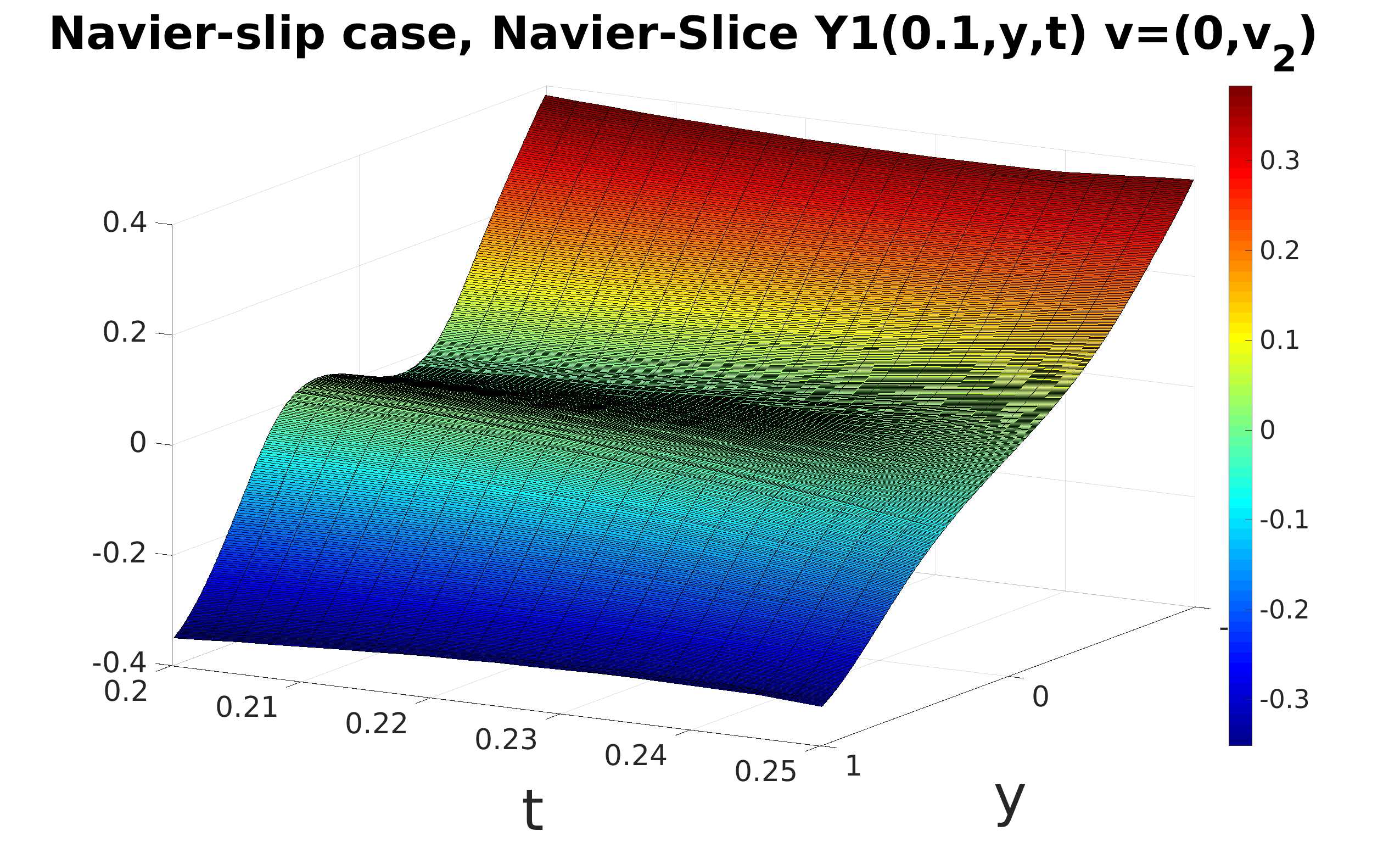} 
   		\hspace{0.01cm} 
		\includegraphics[scale=0.12]{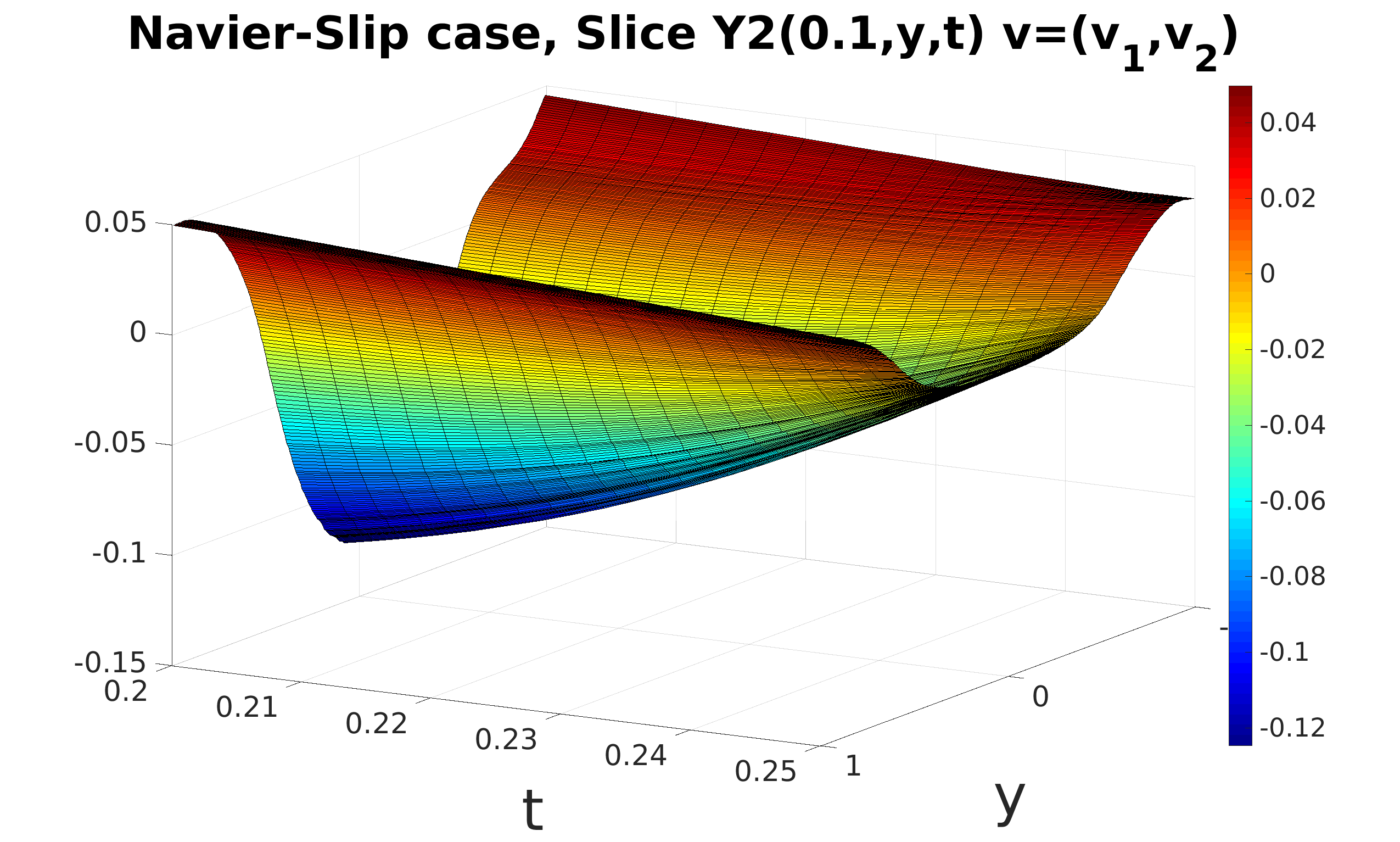} 
   		\hspace{0.01cm}  
		\includegraphics[scale=0.12]{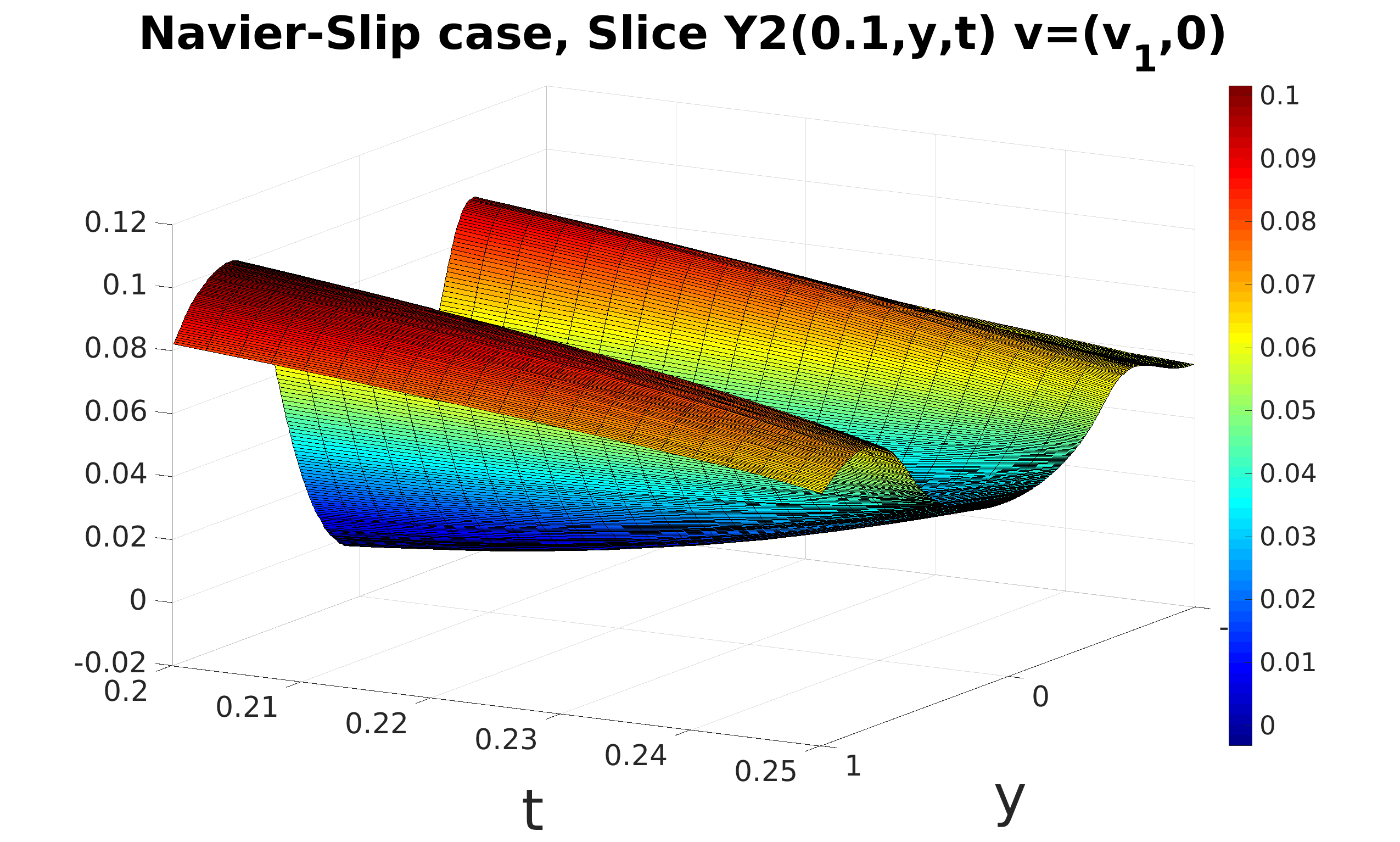} 
   		\hspace{0.01cm}  	
		\includegraphics[scale=0.12]{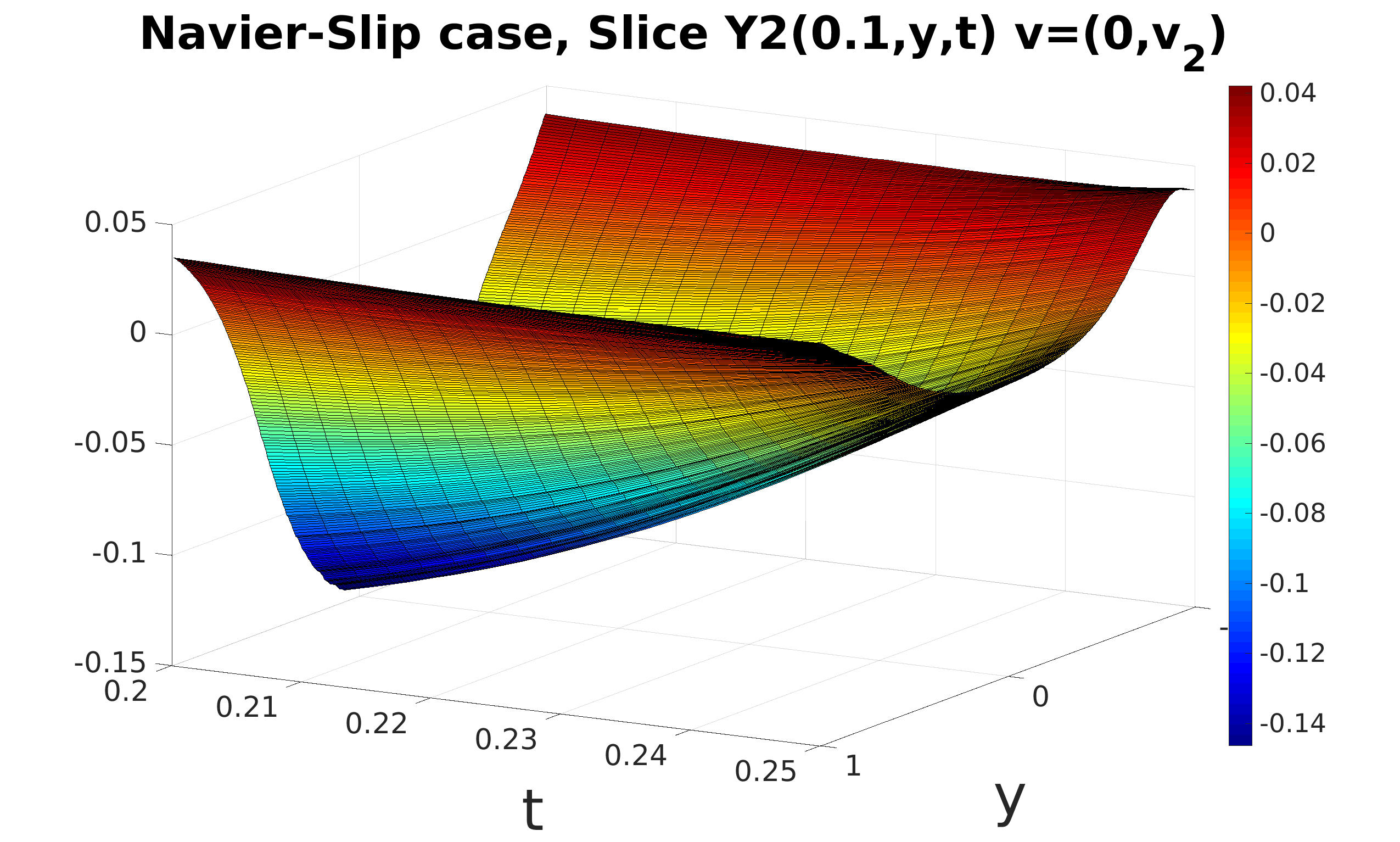} 
   		\hspace{0.01cm}  
		\end{center}\par
	\caption{Navier--slip boundary conditions. Cut sections of the components of the state $\mathbf{y}=(y_1,y_2)$ by the plane $x=0.1$. Cut section of  $y_1$ (first row) and 
	cut section of  $y_2$ (second row). A different control function is represented in each column. }
	\label{fig.cuts.states.RBFs.slip}
	\end{figure}

	On the other hand, for every possible structure of control mentioned above, Figure \ref{fig.cuts.states.RBF.dirich} and Figure \ref{fig.cuts.states.RBFs.slip} show 
	cuts of the first and second component of the state 
	$\mathbf{y}=(y_1,y_2)$ with the plane $x=0.1$.  Every figure displays a type of boundary condition. In order to appreciate a major difference
	in each case (due to scale), we decide to draw the slice at the time interval $[0.2, 0.25]$.

	\begin{table}[H]
	\tiny
	\begin{center}
	\begin{tabular}{ccccc}
	\hline 
	t & $\mathbf{v}=\mathbf{0}$ & $\mathbf{v}=(v_{1},v_{2})$ & $\mathbf{v}=(v_{1},0)$ & $\mathbf{v}=(0,v_{2})$\tabularnewline
	\hline 
	0.005 & 2.30E+00 & 2.30E+00 & 2.30E+00 & 2.30E+00\tabularnewline
	0.025 & 1.28E+00 & 1.27E+00 & 1.27E+00 & 1.27E+00\tabularnewline
	0.05 & 6.11E-01 & 6.08E-01 & 6.09E-01 & 6.09E-01\tabularnewline
	0.075 & 2.93E-01 & 2.89E-01 & 2.90E-01 & 2.90E-01\tabularnewline
	0.1 & 1.40E-01 & 1.37E-01 & 1.38E-01 & 1.38E-01\tabularnewline
	0.125 & 6.73E-02 & 6.36E-02 & 6.44E-02 & 6.44E-02\tabularnewline
	0.15 & 3.22E-02 & 2.87E-02 & 2.94E-02 & 2.94E-02\tabularnewline
	0.175 & 1.55E-02 & 1.20E-02 & 1.28E-02 & 1.28E-02\tabularnewline
	0.2 & 7.40E-03 & 4.37E-03 & 5.01E-03 & 5.01E-03\tabularnewline
	0.225 & 3.55E-03 & 1.18E-03 & 1.67E-03 & 1.67E-03\tabularnewline
	0.25 & 1.70E-03 & 2.18E-04 & 4.68E-04 & 4.66E-04\tabularnewline
	\end{tabular}%
	\quad \quad
	\begin{tabular}{ccccc}
	\hline 
	t & $\mathbf{v}=\mathbf{0}$ & $\mathbf{v}=(v_{1},v_{2})$ & $\mathbf{v}=(v_{1},0)$ & $\mathbf{v}=(0,v_{2})$\tabularnewline
	\hline 
	0.005 & 2.34E+00 & 2.20E+00 & 2.23E+00 & 2.22E+00\tabularnewline
	0.025 & 1.62E+00 & 1.19E+00 & 1.29E+00 & 1.28E+00\tabularnewline
	0.05 & 1.34E+00 & 8.22E-01 & 9.43E-01 & 9.50E-01\tabularnewline
	0.075 & 1.26E+00 & 6.89E-01 & 8.18E-01 & 8.33E-01\tabularnewline
	0.1 & 1.23E+00 & 5.93E-01 & 7.29E-01 & 7.49E-01\tabularnewline
	0.125 & 1.22E+00 & 5.02E-01 & 6.41E-01 & 6.65E-01\tabularnewline
	0.15 & 1.21E+00 & 4.14E-01 & 5.53E-01 & 5.78E-01\tabularnewline
	0.175 & 1.21E+00 & 3.29E-01 & 4.64E-01 & 4.88E-01\tabularnewline
	0.2 & 1.21E+00 & 2.46E-01 & 3.75E-01 & 3.97E-01\tabularnewline
	0.225 & 1.21E+00 & 1.61E-01 & 2.88E-01 & 3.06E-01\tabularnewline
	0.25 & 1.21E+00 & 9.63E-02 & 2.23E-01 & 2.36E-01\tabularnewline
	\end{tabular}

	\end{center} \par
	\caption{$L^2$--norm square of the solution of the null control problem with different internal controls. 
		Dirichlet boundary condition (left) and Navier--slip boundary condition (right).}
	\label{table.finalstate.fewcontrols}	
	\end{table}
	
	\begin{Obs}
	 Observe that Figures \ref{fig.cuts.states.RBF.dirich}, \ref{fig.cuts.states.RBFs.slip} allow us to visualize the impact of every boundary condition in the null control problem 
	 \eqref{intro.Stokes.system}, that is, for every possible control function $\mathbf{v}$ supported in the observatory 
	 ${\omega}{=\{(x,y)\in\mathbb{R}^2: x^2+y^2<0.25\}}$, even with one single scalar control,  the system  \eqref{intro.Stokes.system} with Dirichlet 
	 boundary conditions provides a faster solution in the sense that the $L^2$--norm of the state goes to zero with a higher order than by considering homogeneous 
	 Navier--slip conditions. Numerically, this behaviour was described in Table \ref{table.finalstate.fewcontrols}, however, 
	 to understand the balance between Dirichlet and Navier--slip conditions and the convergence order of the state associated to the null control problem, it is 
	 interesting  to consider friction between the fluid and the boundary (that is, nonhomogeneous Navier--slip conditions) as well as review the stability for the Stokes
	 system with Navier--slip conditions \cite{ding2018stability}. 
	\end{Obs}
		
\subsection{\normalsize{FEM}}	
	Taking as starting point the classical optimal control problem for the Stokes system
	\cite{glowinski2008exact}, we can solve the optimality system given in
	\eqref{intro.Stokes.system}
	\eqref{Stokes1.optimality.condition} and \eqref{Stokes.adjointsystem1} in a similar sense. Here, we only test the Dirichlet boundary conditions. 
	In our case, the time--space discretization of the coupled system  \eqref{intro.Stokes.system}, 
	\eqref{Stokes.adjointsystem1} lies in a mixed 
	finite element formulation in space using $\mathcal{P}_2$--type elements for the velocity and 
	$\mathcal{P}_1$--type elements for the pressure, meanwhile finite differences are used for the time 
	discretization (see \cite{glowinski1992finite,girault2012finite, allaire2005analyse} for a complete review).
	 
\begin{figure}[H]
	\begin{center}
		\includegraphics[scale=0.22]{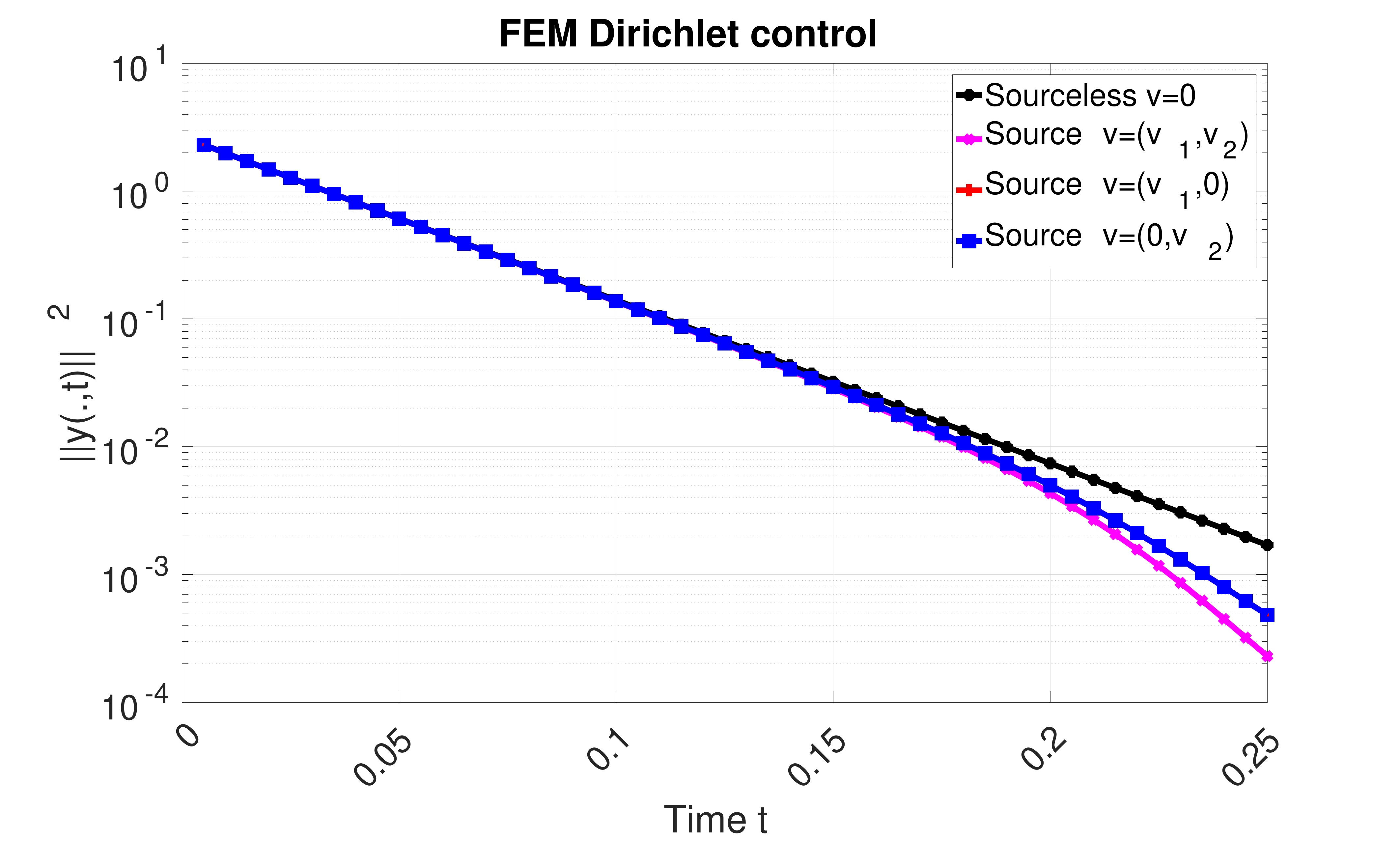}
	\end{center}\par
	\caption{Evolution in time of the $L^2$--norm square for the solution of the null control problem with Dirichlet boundary conditions and  $v=0$ (black), $v=(v_1,v_2)$ (pink), 
	$v=(v_1,0)$ (red) and $v=(0,v_2)$ (blue).}
	\label{fig.finalstate.fewcontrols.fem}
	\end{figure}

	 In this case,  the number of iterations to achieve the stopping criteria $\epsilon=10^{-8}$ in the CGM implemented is 17, for every control function
	 (i.e., $\mathbf{v}=(v_{1},v_{2}), \, \mathbf{v}=(v_{1},0),\, \mathbf{v}=(0,v_{2})$). 
	\noindent Table  \ref{fig.finalstate.fewcontrols.fem} and Figure \ref{fig.finalstate.fewcontrols.fem} display the evolution in time of the $L^2$--norm of the velocity vector field
	$\mathbf{y}=(y_1,y_2)$, which represents the solution to the null control problem \eqref{intro.Stokes.system}, and where  the control function $v$ has different structure, 
	namely, $\mathbf{v}=0$, $\mathbf{v}=(v_1,v_2)$, $\mathbf{v}=(v_1,0)$ and $\mathbf{v}=(0,v_2)$. 

	\begin{table}[H]
	\begin{center}
		\begin{tabular}{ccccc}
	\hline 
	t &\hspace{0.5cm} $\mathbf{v}=\mathbf{0}$ &\hspace{0.5cm} $\mathbf{v}=(v_{1},v_{2})$ &\hspace{0.5cm} $\mathbf{v}=(v_{1},0)$ &\hspace{0.5cm} $\mathbf{v}=(0,v_{2})$
	\tabularnewline
	\hline 
	0.005 &\hspace{0.5cm} 2.30E+00 &\hspace{0.5cm} 2.30E+00 &\hspace{0.5cm} 2.30E+00 &\hspace{0.5cm} 2.30E+00\tabularnewline
	0.025 &\hspace{0.5cm} 1.28E+00 &\hspace{0.5cm} 1.27E+00 &\hspace{0.5cm} 1.27E+00 &\hspace{0.5cm} 1.27E+00\tabularnewline
	0.05 &\hspace{0.5cm} 6.11E-01 &\hspace{0.5cm} 6.08E-01 &\hspace{0.5cm} 6.09E-01 &\hspace{0.5cm} 6.09E-01\tabularnewline
	0.075 &\hspace{0.5cm} 2.93E-01 &\hspace{0.5cm} 2.89E-01 &\hspace{0.5cm} 2.90E-01 &\hspace{0.5cm} 2.90E-01\tabularnewline
	0.1 &\hspace{0.5cm} 1.40E-01 &\hspace{0.5cm} 1.37E-01 &\hspace{0.5cm} 1.37E-01 &\hspace{0.5cm} 1.37E-01\tabularnewline
	0.125 &\hspace{0.5cm} 6.73E-02 &\hspace{0.5cm} 6.36E-02 &\hspace{0.5cm} 6.43E-02 &\hspace{0.5cm} 6.43E-02\tabularnewline
	0.15 &\hspace{0.5cm} 3.22E-02 &\hspace{0.5cm} 2.86E-02 &\hspace{0.5cm} 2.94E-02 &\hspace{0.5cm}2.94E-02\tabularnewline
	0.175 &\hspace{0.5cm} 1.54E-02 &\hspace{0.5cm} 1.20E-02 &\hspace{0.5cm} 1.27E-02 &\hspace{0.5cm} 1.27E-02\tabularnewline
	0.2 &\hspace{0.5cm} 7.40E-03 &\hspace{0.5cm} 4.32E-03 &\hspace{0.5cm} 5.00E-03 &\hspace{0.5cm} 5.00E-03\tabularnewline
	0.225 &\hspace{0.5cm} 3.55E-03 &\hspace{0.5cm} 1.17E-03 &\hspace{0.5cm} 1.67E-03 &\hspace{0.5cm} 1.67E-03\tabularnewline
	0.25 &\hspace{0.5cm} 1.70E-03 &\hspace{0.5cm} 2.29E-04 &\hspace{0.5cm} 4.82E-04 &\hspace{0.5cm} 4.82E-04\tabularnewline
	\end{tabular}

	\end{center}\par
	\caption{Evolution in time of the $L^2$--norm for the solution of the null control problem with Dirichlet boundary conditions and few scalar controls.}\label{table.finalstate.fewcontrols.fem}
	\end{table}
	
	\begin{figure}[H]
		\begin{center}
		\includegraphics[scale=0.12]{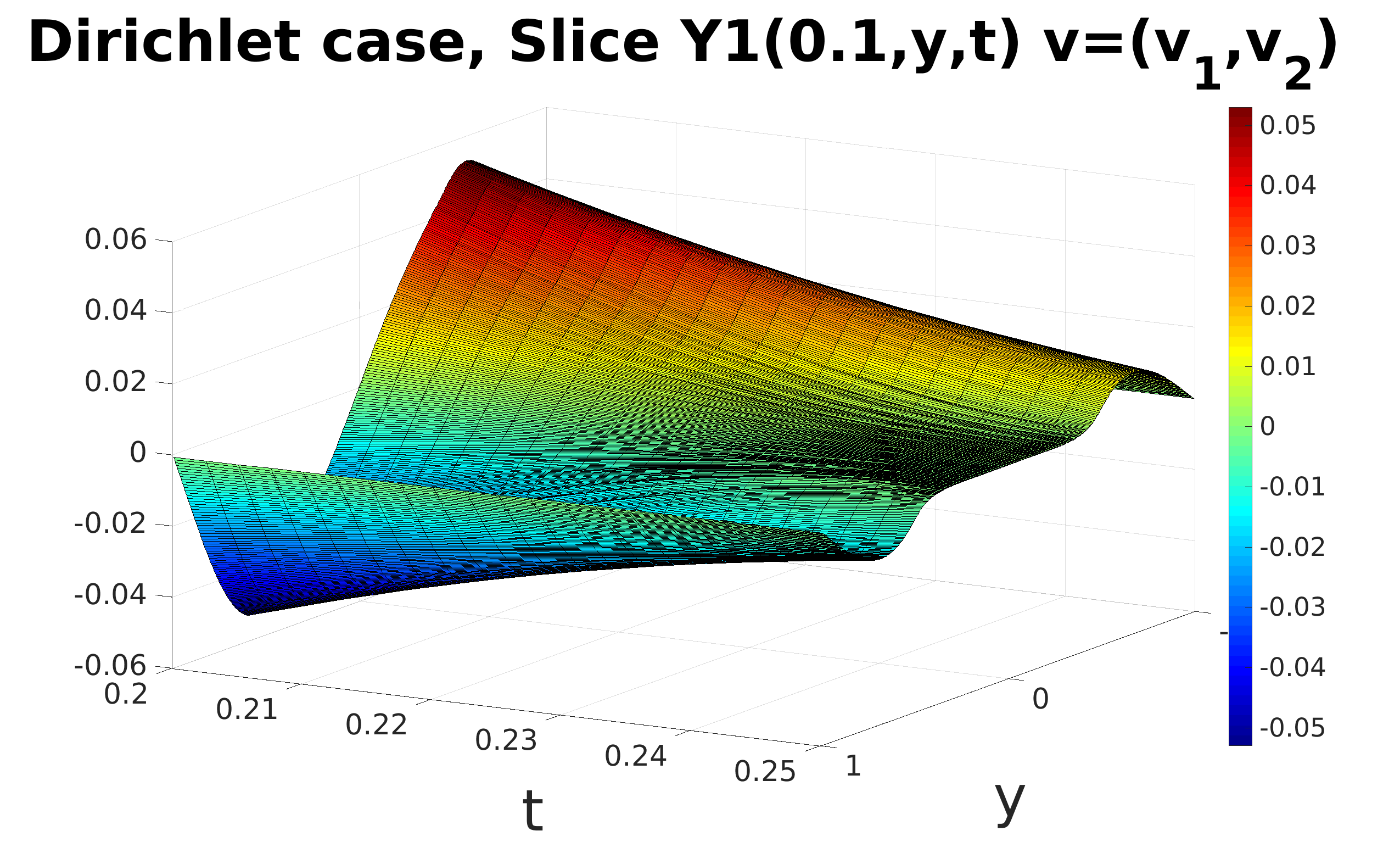} 
   		\hspace{0.01cm} 
		\includegraphics[scale=0.12]{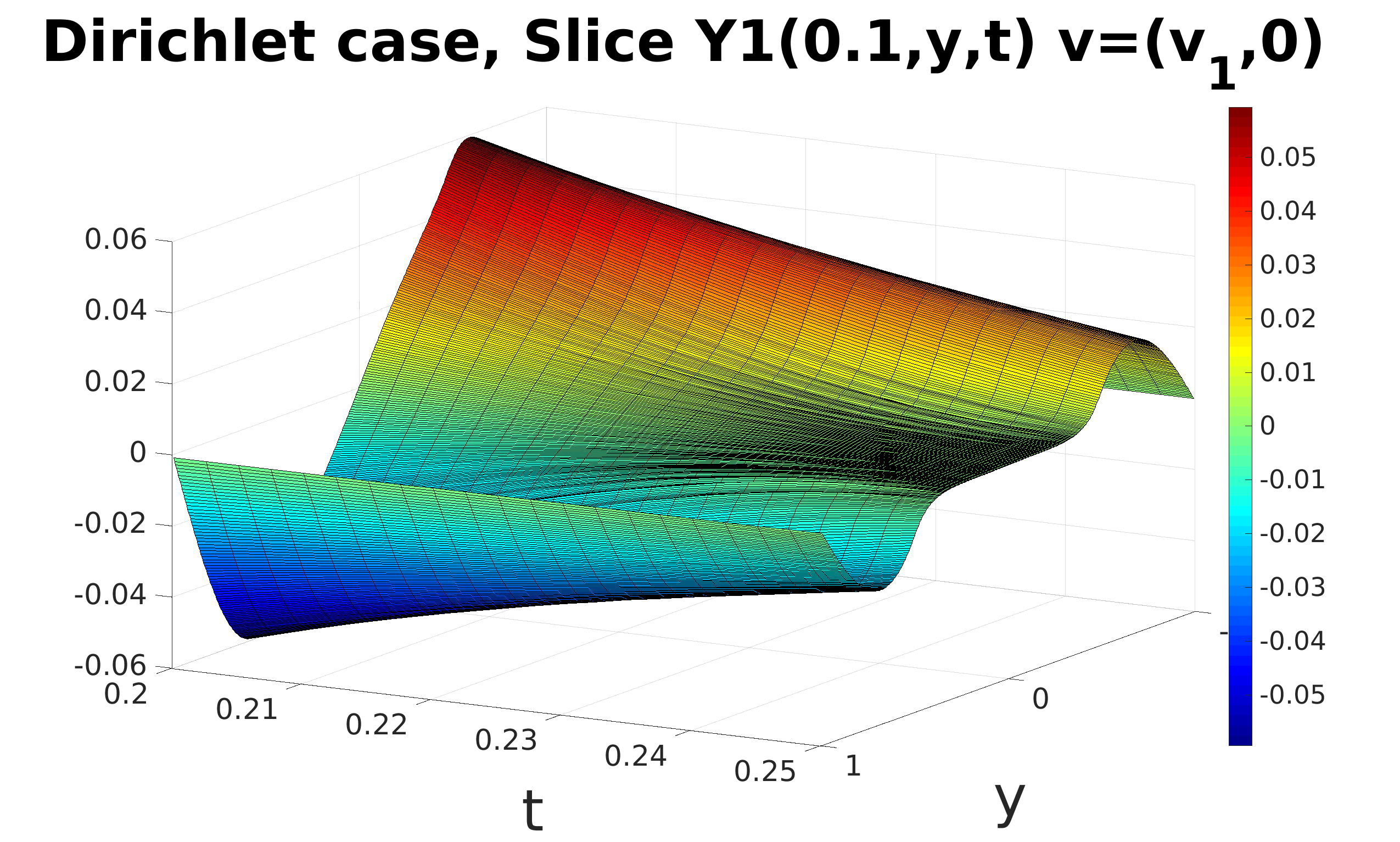} 
   		\hspace{0.01cm} 
		\includegraphics[scale=0.12]{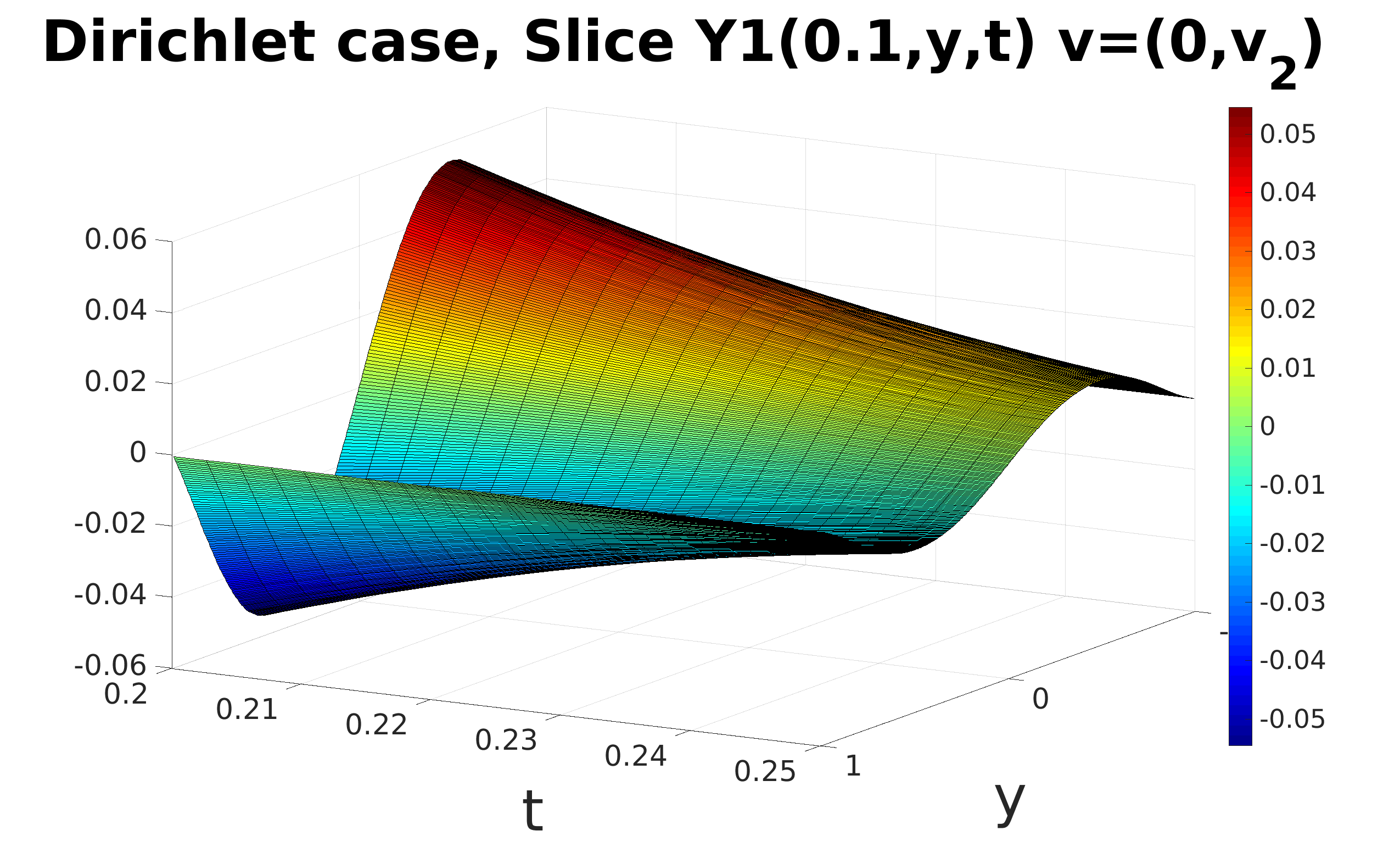} 
   		\hspace{0.01cm} 
		\includegraphics[scale=0.12]{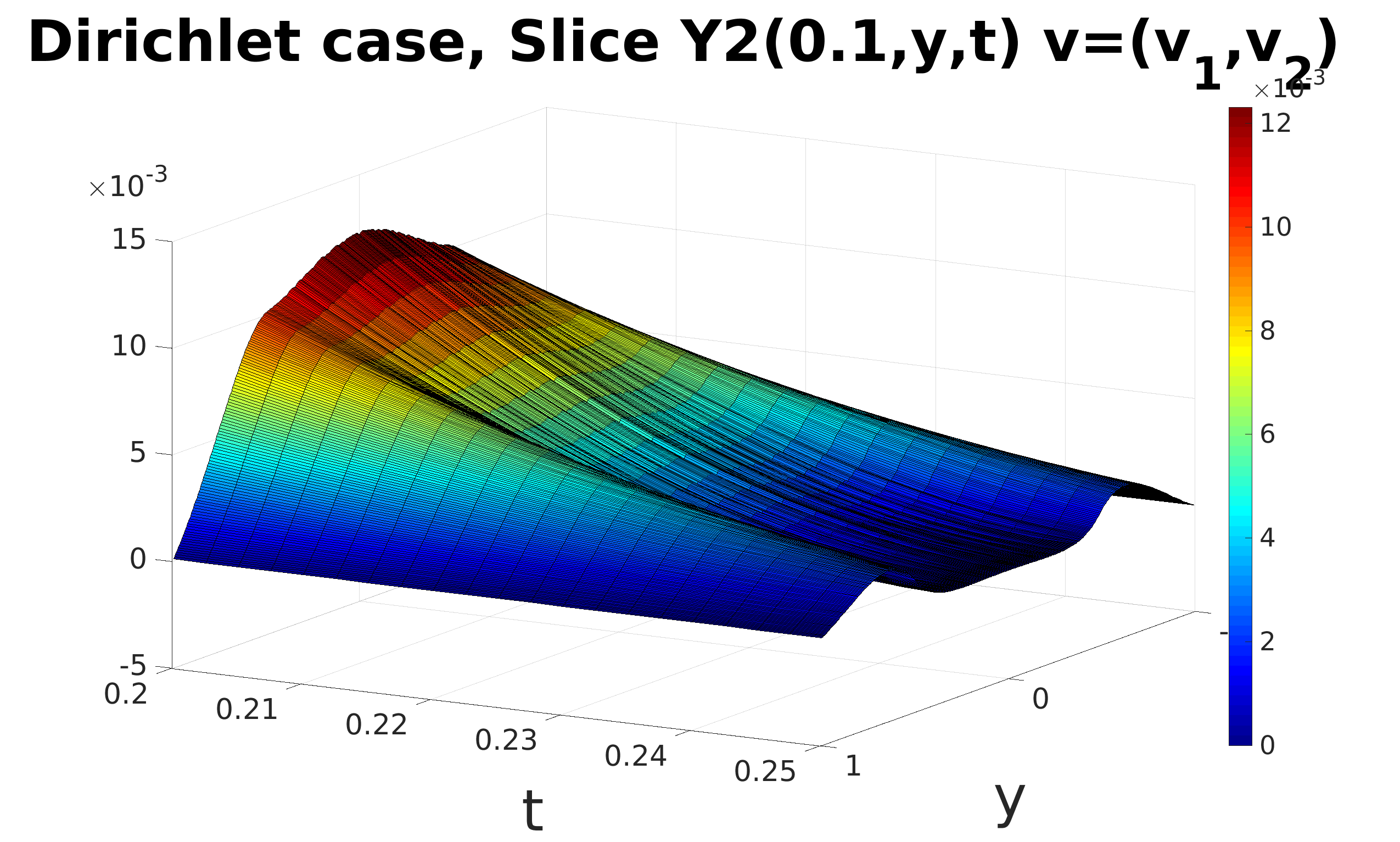} 
   		\hspace{0.01cm} 
		\includegraphics[scale=0.12]{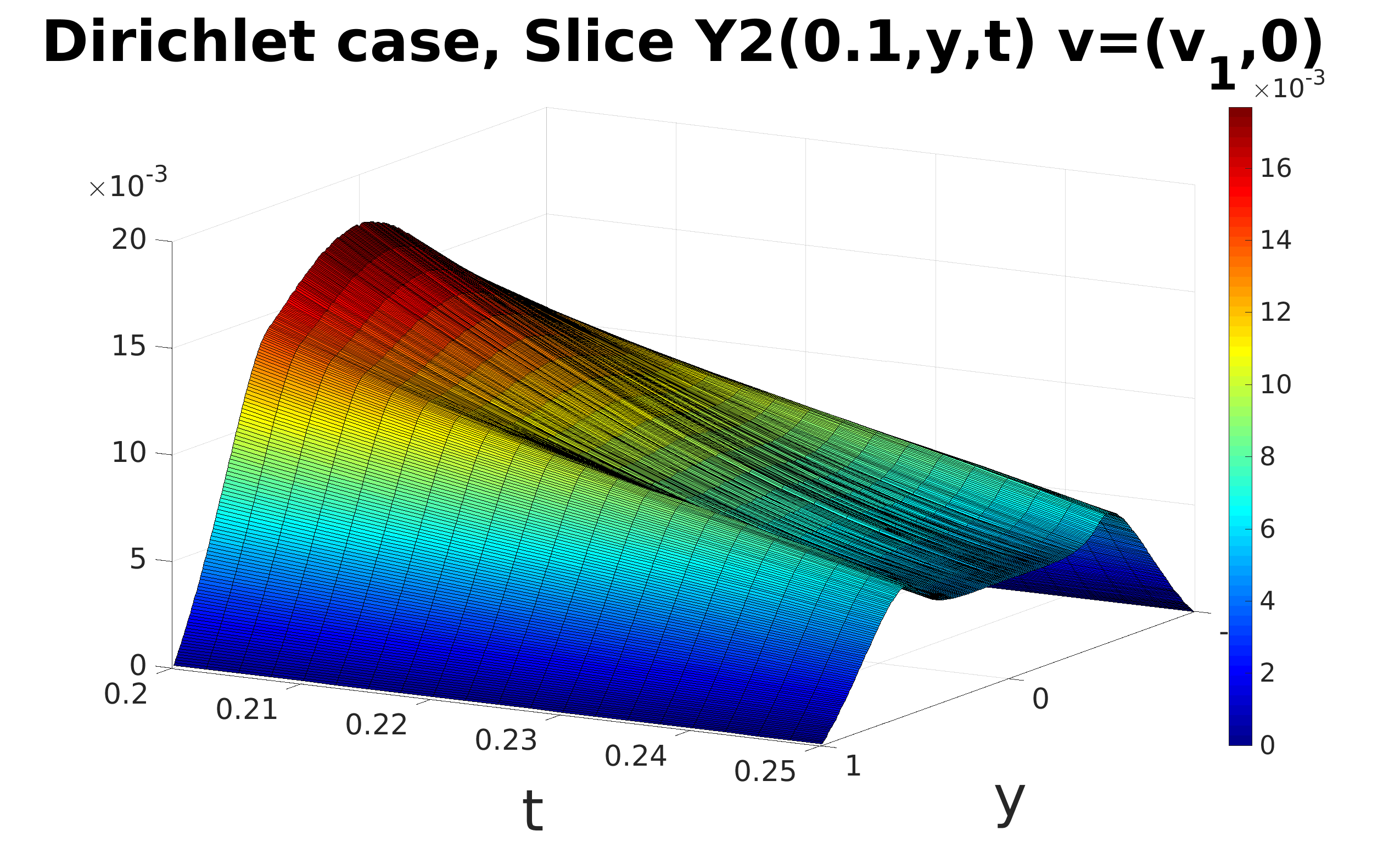} 	 
   		\hspace{0.01cm} 
		\includegraphics[scale=0.12]{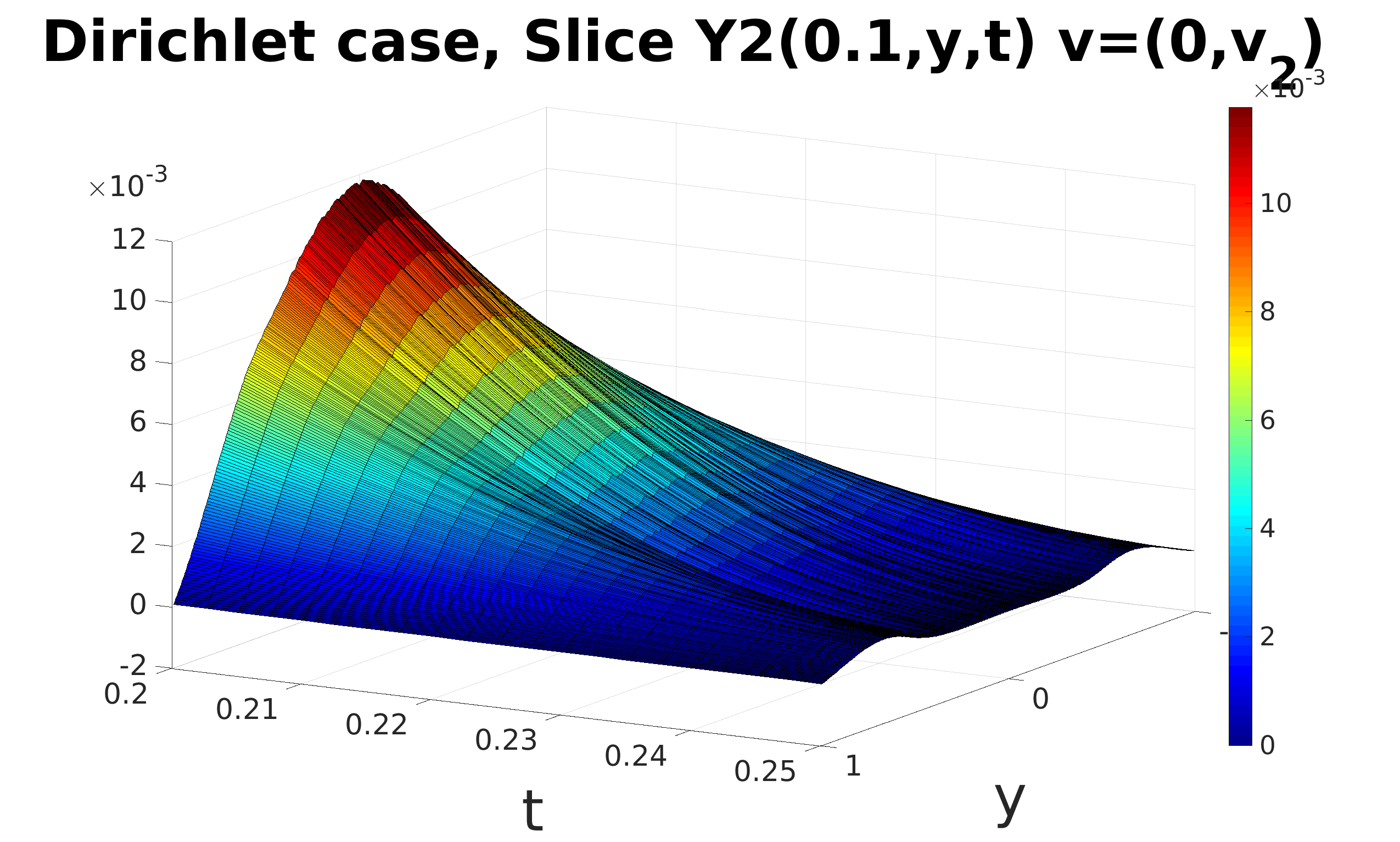} 
		\end{center}\par
	\caption{Dirichlet boundary conditions. Cut sections of the components of the state $\mathbf{y}=(y_1,y_2)$ by the plane $x=0.1$. 
	Cut section of  $y_1$ (first row) and cut section of  $y_2$ (second column). A different control function is represented in each column.}
	\label{fig.cuts.states.fem}
	\end{figure}	

	\begin{Obs}
	As we can see from Tables \ref{table.finalstate.fewcontrols.fem},\ref{table.finalstate.fewcontrols} and Figures \ref{fig.cuts.states.fem},\ref{fig.cuts.states.RBF.dirich} for 
	the Dirichlet case, the control result using RBF-LHI or FEM are practically the same, but nevertheless RBF-LHI method has the advantage of being meshless and 
	having more accuracy by using divergence free kernels.
	
	Unfortunately,  by observing the required number of iterations for the CGM in its convergence  in both methods, RBF--LHI and FEM, the number is 
	is higher for the RBF--LHI method, this might be explain due that in RBF--LHI we are using $\mathcal{P}_1$--type elements approximation to compute the 
	integrals expressions in the CGM, meanwhile for FEM we are using $\mathcal{P}_2$--type elements.
	\end{Obs}
	
\vskip 1cm

\section{\normalsize{Conclusions and final remarks}}\label{conclusiones}	

In this article we have introduced a radial basis function, (RBFs), method to solve null control problems for the Stokes system with few internal scalar controls and 
	Dirichlet or Navier--slip boundary conditions. As far as to our knowledge, this problem has not been treated in the literature through radial basis function methods. 
	The continuous control problem gives rise to a coupled system of Stokes problems, namely the direct and adjoint systems. A conjugate gradient algorithm, adapted to 
	the RBF setting is used to deal with this problem. Direct RBF solvers were designed, based on divergency free global vector RBF to avoid saddle point problems and 
	thus inf/sup conditions. The pressure  and the velocity are incorporated in the solution, in the sense of Wendland  \cite{wendland2009divergence}, so that slip boundary conditions can
	 be implemented. Within the paper, we formulate direct algorithms both for stationary and evolutionary Stokes systems. Stability analysis in the sense of 
	 Chinchapatnam \cite{chinchapatnam2006unsymmetric}, that is by bounding the absolute value of the spectral radius of the discrete system, is performed for each of the techniques formulated within the article. 
	 
We note that recently, Keim and Wendland in \cite{Keim2016AHA} proposed a RBF compact support method for evolutionary Stokes problem with Dirichlet boundary conditions. This is an alternative methodology to the LHI-technique presented in this article.

It is worth pointing that when the stability condition was not satisfied the numerical solution was numerically ill posed.

The treatment of Navier--slip boundary conditions by the LHI method, which gave excellent results of order $10^{-4}$ in the $L^{\infty} $--norm, were obtained by incorporating the slip boundary operator in the ansatz. This was done without increasing the node density near the boundary or incorporating ghost nodes. 

The LHI divergency free technique opens the possibility of treating a large number of data, whenever the fill distance and the shape parameter does not produces a singular value of the corresponding Gram matrix in extended precession. The use of extended precession to deal with the bad condition number has been formerly proposed by E. Kansa, \cite{Kansa2017} and S. Sarra, \cite{Sarra2011}, providing a good possibility to solve many real problems. However, we note that, recently, see Fronberg et al. \cite{Fornberg2015} and references therein, have formulated several techniques to deal with the bad condition problem. These techniques will be 
	the subject of further works to treat null control problems.

 This paper opens the possibility of extending the stability analysis in order to obtain an explicit CLF type condition depending on the shape parameter, fill distance and the diffusion coefficient. Also, it may be possible to extend this analysis to the treatment of Navier--Stokes control problems. This will be considered in further works.

\vskip 1cm
\noindent \textbf{Acknowledgements.}  
    The research was partially supported by Fondecyt grant  3180100 (Cristhian Montoya). The authors acknowledge the CONACYT, project Fordecyt 
    265667. This work also was supported by the National Autonomous University of M\'exico [grant: PAPIIT, IN102116] and by  Network of Mathematics and Development of CONACYT. We wish also to thanks Dr.Jesús Lopéz Estrada for providing some inside in development of the numerical code for this paper.
    
\vskip 1cm

\bibliographystyle{alpha}
\bibliography{biblio}
\end{document}